\documentclass[a4paper,12pt]{article}     

\usepackage{graphicx}		  
\usepackage{amsmath,amssymb}	

\newtheorem{theorem}{Theorem}[section]
\newtheorem{lemma}[theorem]{Lemma}
\newtheorem{cor}[theorem]{Corollary}

\newtheorem{assumption}{Assumption}
\newtheorem{remark}[theorem]{Remark}

\newenvironment{proof}{
\hspace*{-9mm}
{ \it Proof.}}
{\hfill {$\square$}\vspace{1.5em}}

\begin{document}

\begin{center}{
{\Large 
Minimal Charts}
\vspace{10pt}
\\ 
Teruo NAGASE and Akiko SHIMA\footnote
{
The second author was supported by JSPS KAKENHI Grant Number 23540107. 
}
}
\end{center}

\begin{abstract}
In this paper, we give definitions of three kinds of minimal charts, and
we investigate 
properties of minimal charts
and establish fundamental theorems 
characterizing minimal charts.
To classify charts with two or three crossings
we use the fundamental theorems.
In the future paper,
we give an enumeration of the charts with two crossings.
\end{abstract}


\section{Introduction}
\label{s:Intro}

Charts are oriented labeled graphs
in a disk 
with three kinds of vertices
called black vertices, crossings,
and white vertices 
(see page 3 for the precise definition of charts).
From a chart, we can construct an oriented closed surface 
embedded in 4-space ${\Bbb R}^4$ 
 (see \cite[Chapter 14, Chapter 18 and Chapter 23]{BraidBook}). 
A C-move 
is a local modification between two charts
in a disk (see Section~\ref{s:Prel} for C-moves).
A C-move between two charts induces 
an ambient isotopy between oriented closed surfaces 
corresponding to the two charts.
Two charts are said to be {\it C-move equivalent} 
if there exists
a finite sequence of C-moves 
which modifies one of the two charts 
to the other.

We will work in the PL or smooth category. 
All submanifolds are assumed to be locally flat.
A {\it surface link} is a closed surface embedded in 4-space ${\Bbb R}^4$. 
A {\it $2$-link} is a surface link each of whose connected component is a $2$-sphere.
A {\it $2$-knot}
is a surface link which is a 2-sphere.
An orientable surface link is called a 
{\it ribbon surface link}
if there exists an immersion of a $3$-manifold $M$
into ${\Bbb R}^4$ sending the boundary of $M$ onto the surface link
such that each connected component of $M$ is a handlebody and
its singularity
consists of ribbon singularities,
here a ribbon singularity
is a disk in the image of $M$
whose pre-image consists of 
two disks;
one of the two disks is a proper disk of $M$ 
and
the other is a disk in the interior of $M$.
In the words of charts,
a ribbon surface link is
a surface link corresponding to a {\it ribbon chart}, 
a chart C-move equivalent to 
a chart
without white vertices \cite{BraidThree}.
A chart is called a {\it $2$-link chart}
if a surface link corresponding to the chart is a 2-link.

In this paper, 
we denote the closure, the interior, 
the boundary, and the complement of $(...)$ by 
$Cl(...)$, Int$(...)$, 
$\partial(...)$, $(...)^c$ 
respectively.
Also for a finite set $X$, 
the notation $|X|$ denotes 
the number of elements in $X$.

At the end of this paper,
there is the index of new words and notations 
introduced in this paper.

Kamada showed that 
any $3$-chart is a ribbon chart 
\cite{BraidThree}.
Kamada's result was extended by Nagase and Hirota:
Any $4$-chart with at most one crossing
is a ribbon chart \cite{NH}.
We showed that any $n$-chart with at most one crossing is a ribbon chart
\cite{OneCrossing}.
We also showed that any $2$-link chart 
with at most two crossings
 is a ribbon chart  \cite{TwoCrossingI},  
 \cite{TwoCrossingII}.

\begin{figure}
\begin{center}
\includegraphics{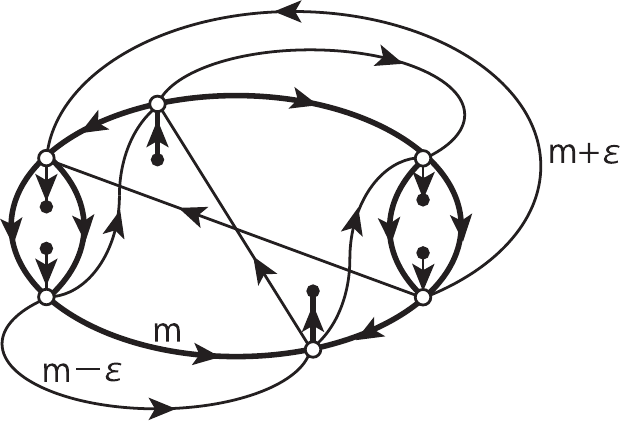}
\end{center}
\caption{ \label{fig01} The letter $m$ is a label and $\varepsilon=\pm1$. }
\end{figure}

\begin{figure}
\begin{center}
\includegraphics{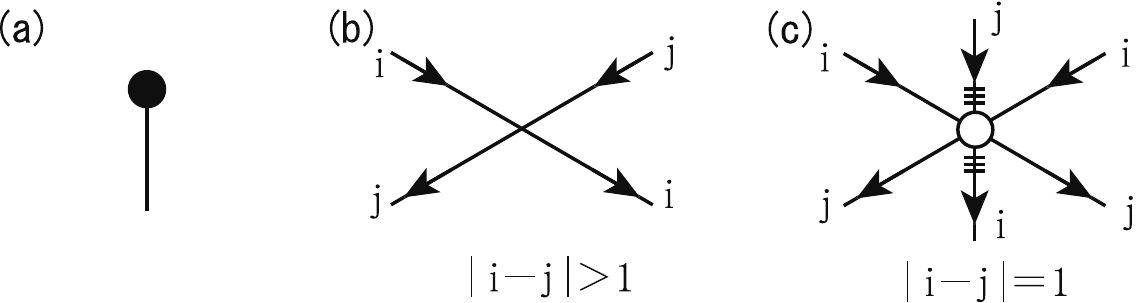}
\end{center}
\caption{ \label{fig02} (a) A black vertex. (b) A crossing. (c) A white vertex. 
Each arc with three transversal short arcs is a middle arc at the white vertex. }
\end{figure}

Charts have strong conditions 
on orientations of arcs around vertices. 
In a small neighborhood of
each white vertex,
there are six short arcs,
three consecutive arcs are
oriented inward 
and
the other three are outward
(see Fig.~\ref{fig02}(c)).
Among six short arcs
in a small neighborhood of
a white vertex,
a central arc of each three consecutive arcs
oriented inward (resp. outward) 
is called a  
{\it middle arc} at the white vertex.
Observing precisely middle arcs, 
orientations of edges, and 
a part of a chart cutting by a disk 
called a {\it tangle},
we shall prove the following theorem \cite{ThreeCrossing}: 
\begin{enumerate}
\item[]Any $2$-link chart with at most three crossings is C-move equivalent to either a ribbon chart, 
or the disjoint union of a ribbon chart and a chart as shown in Fig.~\ref{fig01} or its ``reflection''.
\end{enumerate}
In this paper
we establish fundamental theorems
characterizing $c$-minimal charts, $w$-minimal charts and $cw$-minimal charts.
For the classification theorem above, 
we use the fundamental theorems
obtained in this paper.

For a $4$-chart as shown in Fig.~\ref{fig01},
 we obtain a $2$-twist spun trefoil 
by setting $m=2$ (see \cite[p. 144]{BraidThree}, \cite[p. 170]{BraidBook}).
It is well known that 
the $2$-knot is not a ribbon $2$-knot.  
On the other hand, Hasegawa showed that 
if a non-ribbon chart representing a $2$-knot 
is minimal, 
then the chart must possess at least six white vertices \cite{H1} 
where a minimal chart $\Gamma$ means 
its complexity $(w(\Gamma), -f(\Gamma))$ 
is minimal among 
the charts C-move equivalent to 
the chart $\Gamma$ with respect to 
the lexicographic order 
of pairs of integers, 
here 
$w(\Gamma)$ is the number of white vertices 
in $\Gamma$, 
$f(\Gamma)$ is the number of free edges 
in $\Gamma$.
Here a {\it free edge} is 
an edge of $\Gamma$ containing
two black vertices.
Nagase, Ochiai, and Shima showed that
there does not exist a minimal chart 
with exactly five white vertices \cite{ONS}.
Nagase and Shima show that
there does not exist a minimal chart
with exactly seven white vertices \cite{ChartApp1},\cite{ChartAppII},\cite{ChartAppIII},\cite{ChartApp4},\cite{ChartApp5}.
Ishida, Nagase, and Shima showed that
any minimal chart with exactly four white vertices
 is C-move equivalent to a chart in two kinds of classes \cite{INS}.

Let $n$ be a positive integer.
An $n$-{\it chart} (a braid chart of degree $n$ \cite{KS}
or a surface braid chart of degree $n$ \cite{BraidBook}) 
is 
an oriented labeled graph 
in the interior of a disk,
which may be empty 
or
have closed edges without vertices
satisfying the following four conditions
(see Fig.~\ref{fig02}):
\begin{enumerate}
\item[(i)] 
Every vertex has degree $1$, $4$, or $6$.
\item[(ii)] 
The labels of edges are 
in $\{1,2,\dots,n-1\}$.
\item[(iii)]
In a small neighborhood of
each vertex of degree $6$,
there are six short arcs,
three consecutive arcs are
oriented inward 
and
the other three are outward,
and
these six are labeled $i$ and $i+1$
alternately for some $i$,
where the orientation and label of
each arc are inherited from
the edge containing the arc.
\item[(iv)]
For each vertex of degree $4$,
diagonal edges have the same label
and
are oriented coherently,
and the labels $i$ and $j$ of
the diagonals satisfy $|i-j|>1$.
\end{enumerate}
We call a vertex of degree $1$ a {\it black vertex},
a vertex of degree $4$ a {\it crossing}, and 
a vertex of degree $6$ a {\it white vertex}
respectively.
Among six short arcs
in a small neighborhood of
a white vertex,
a central arc of each three consecutive arcs
oriented inward (resp. outward) 
is called a   
{\it middle arc} at the white vertex
(see Fig.~\ref{fig02}(c)).
For each white vertex $v$, 
there are two middle arcs at $v$ 
in a small neighborhood of
the white vertex $v$.


\begin{figure}
\begin{center}
\includegraphics{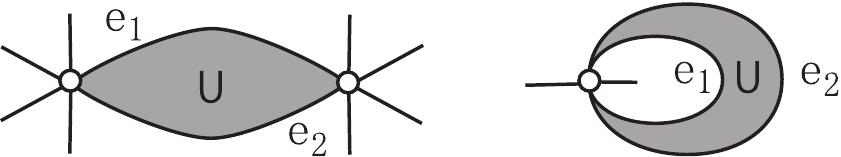}
\end{center}
\caption{ Bigons.\label{fig03} }
\end{figure}

Let $\Gamma$ be a chart. 
Let $e_1$ and $e_2$ be edges of $\Gamma$
which connect two white vertices $w_1$ and $w_2$
where possibly $w_1=w_2$.
Suppose that 
the union $e_1\cup e_2$ bounds 
an open disk $U$.
Then $Cl(U)$ 
is called 
a {\it bigon} of $\Gamma$
provided that
any edge containing $w_1$ or $w_2$ 
does not intersect the open disk $U$
(see Fig.~\ref{fig03}).
Since $e_1$ and $e_2$ are edges of $\Gamma$, 
neither $e_1$ nor $e_2$ 
contains a crossing.

%

Let $\Gamma$ be a chart.
Let $c(\Gamma)$ and
$b(\Gamma)$ be 
the number of crossings,
and
the number of bigons of $\Gamma$
respectively.
We define complexities as follows
(see \cite{BraidThree} 
for the original complexity of charts):\\
$\bullet$ 
The 4-tuple $(c(\Gamma),w(\Gamma),-f(\Gamma),-b(\Gamma))$ is called a 
{\it $c$-complexity} of $\Gamma$.\\
$\bullet$ 
The 4-tuple $(w(\Gamma),c(\Gamma),-f(\Gamma),-b(\Gamma))$ is called a 
{\it $w$-complexity} of $\Gamma$.\\
$\bullet$ 
The 3-tuple $(c(\Gamma)+w(\Gamma),-f(\Gamma),-b(\Gamma))$ is called a 
{\it $cw$-complexity} of $\Gamma$.

A chart $\Gamma$ is said to be 
{\it $c$-minimal $($resp. $w$-minimal or $cw$-minimal$)$} if
its $c$-complexity (resp. $w$-complexity or $cw$-complexity) is minimal among the charts 
which are C-move equivalent to 
the chart $\Gamma$
with respect to 
the lexicographical order of the 
4-tuple (or 3-tuple) of the integers.

{\it In this paper, 
if a chart is $c$-minimal, $w$-minimal or $cw$-minimal,}

{\it then we say that the chart is minimal.}

A {\it hoop} is a closed edge of a chart $\Gamma$
that contains neither crossings nor white vertices.
Therefore a hoop decomposes $\Gamma$ into disjoint pieces:
an inside,  an outside, and itself.
A hoop is said to be {\it simple} 
if one of the complementary domains
of the hoop
does not contain any white vertices.

Let $\Gamma$ be a chart.
For each label $m$, we define
$$\Gamma_m=\text{ 
the union of 
all the edges of label $m$ and 
their vertices in }\Gamma .$$

Let $\Gamma$ be a chart, and $m$ a label of $\Gamma$. 
Let $E$ be a disk 
with $\partial E\subset\Gamma_m$.
Then the disk $E$ is called 
a {\it $3$-color disk} provided that 
\begin{enumerate}
\item[(i)] the disk $E$ does not contain any crossings, and
\item[(ii)] 
$\Gamma\cap E\subset \Gamma_{m-1}\cup\Gamma_m
\cup\Gamma_{m+1}$.
\end{enumerate}
Further a disk $E$ with $\partial E\subset \Gamma_m$
is called a {\it $2$-color disk} 
provided that 
$\Gamma\cap E\subset\Gamma_m\cup\Gamma_{m-1}$ 
or 
$\Gamma\cap E\subset\Gamma_m\cup\Gamma_{m+1}$. 
A $2$-color disk is a special kind of 
$3$-color disk.

Now $3$-color disks and $2$-color disks 
often appear in charts.
For example,
{\it let $m$ be any label of 
a chart $\Gamma$,
and $E$ a disk 
with $\partial E\subset\Gamma_m$ but
without crossings, free edges nor simple hoops.
If $\Gamma$ is a minimal chart,
then we can show that  
$E$ is a $3$-color disk $($see Lemma~\ref{Appendix3ColorDisk} in Section~$\ref{s:Appendix})$.
Further, if $m$ is 
the minimal label or the maximal label 
of the chart, 
then 
$E$ is a $2$-color disk.}
This indicates that 
it is important to investigate $3$-color disks 
and $2$-color disks.

An edge $e$ of $\Gamma$ 
is said to be 
{\it middle at a white vertex $v$} 
if it contains a middle arc at the vertex $v$.

Let $m$ be a label of a chart $\Gamma$.
A simple closed curve in $\Gamma_m$ is 
called a {\it cycle of label $m$}.
Let $C$ be a cycle of label $m$ 
bounding a disk $E$.
An edge $e$ of label $m$ 
is called
{\it an inside $($resp. outside$)$ edge for $C$} 
provided that 
\begin{enumerate}
\item[(i)]
$e\cap C$ consists of 
one white vertex or two white vertices, and
\item[(ii)]
$e\subset E$ (resp. $e\subset Cl(E^c)$).
\end{enumerate}

Let $\Gamma$ be a chart, and 
$m$ a label of $\Gamma$.
For a cycle $C$ of label $m$,
we define\\
$\begin{array}{rl}
{\mathcal{W}}(C)&
\hspace{-2.5mm}= \{ w \ | 
\text{ $w$ is a white vertex in $C$} \},\\
{\mathcal{W}}_I^{{\rm Mid}}(C,m)&
\hspace{-2.5mm}= 
\{ w\in {\mathcal{W}}(C) \ | 
\text{ there exists 
an inside edge for $C$ {\it middle} 
at $w$} \},\\
{\mathcal{W}}_O^{{\rm Mid}}(C,m)&
\hspace{-2.5mm}= 
\{ w\in {\mathcal{W}}(C) \ | 
\text{ there exists 
an outside edge for $C$ {\it middle} 
at $w$} \}.\\
\end{array}
$

\begin{theorem}
\label{Inequation} 
Let $\Gamma$ be a minimal chart, and 
$m$ a label of $\Gamma$.
Let $E$ be a $3$-color disk 
with $\partial E\subset\Gamma_m$
but without free edges nor simple hoops. 
If $\Gamma_m\cap E$ is connected, 
then we have
$$|{\mathcal{W}}_I^{\rm Mid}(\partial E,m)|+2 
\le 
|{\mathcal{W}}_O^{\rm Mid}(\partial E,m)|.$$
\end{theorem}

Let $\Gamma$ be a chart, and 
$D$ a disk.
The pair $(\Gamma\cap D,D)$ is called 
a {\it tangle} provided that 
\begin{enumerate}
\item[(i)]
$\partial D$ does not contain any white vertices,  black vertices nor crossings of $\Gamma$, 
\item[(ii)]
if an edge of $\Gamma$ intersects $\partial D$, 
then the edge intersects $\partial D$ transversely, and
\item[(iii)] $\Gamma\cap D\not=\emptyset$.
\end{enumerate}

Let $\Gamma$ be a chart.
A tangle $(\Gamma\cap D,D)$ is 
called an {\it NS-tangle of label $m$} 
(new significant tangle) 
provided that
\begin{enumerate}
\item[(i)] if $i\neq m$, 
then $\Gamma_i\cap \partial D$ is 
at most one point,
\item[(ii)] 
$\Gamma\cap D$ contains at least one white vertex, and 
\item[(iii)]
for each label $i$, 
the intersection $\Gamma_i\cap D$ contains 
at most one crossing.
\end{enumerate}

\begin{theorem}
\label{NoNS-tangle} 
In a minimal chart, 
there does not exist 
an NS-tangle of any label.
\end{theorem}

The above theorem is 
an extended result of Theorem~$3.5$ in 
\cite{TwoCrossingI}, and
does a significant job for 
the classification of charts 
from the view point of the number of crossings.
The above theorem and 
Consecutive Triplet Lemma (Lemma~\ref{ConsecutiveTriplet}) are our main tools as we see how to use these two tools 
to prove Theorem~\ref{TwoColorGateTangle} below.

Let $\Gamma$ be a chart, and 
$m$ a label of the chart. 
Let $\mathcal W$ be the set of 
all the white vertices of $\Gamma$. 
The closure of a connected component of $\Gamma_m-\mathcal W$ 
is called an {\it internal} edge of label $m$ 
if it contains a white vertex 
but does not contain any black vertex, 
here we consider $\Gamma_m$ as a topological set. 
Thus an internal edge  
begins at a white vertex
passing through several crossings and 
ends at a white vertex.

Let $\Gamma$ be a chart.
A tangle $(\Gamma\cap D,D)$ is said to be
{\it admissible} 
provided that
\begin{enumerate}
\item[(i)] $D$ contains neither free edge nor simple hoop.
\item[(ii)] 
any edge intersecting $\partial D$
is contained in an internal edge, 
\item[(iii)] 
if an internal edge $\overline e$  
intersects $\partial D$,
then each connected component of 
$\overline e\cap D$ contains a white vertex.
\end{enumerate}

Let $\Gamma$ be a chart, 
and $D$ a disk. 
Suppose that an edge $e$ of $\Gamma$ 
transversely intersects $\partial D$.
Let $p$ be a point in 
$e\cap \partial D$, 
and 
$N$ a regular neighbourhood of $p$.
Then the orientation of $e$
induces the one of the arc $e\cap N$.
The edge $e$ is said to be 
{\it locally inward} 
(resp. {\it locally outward}) 
at $p$ 
with respect to $D$
if the oriented arc $e\cap N$ is 
oriented
from a point outside 
(resp. inside) $D$ 
to a point inside 
(resp. outside) $D$.
We often say that 
$e$ is locally inward (resp. outward) at $p$ 
instead of saying that
$e$ is locally inward (resp. outward) at $p$ 
{\it with respect to $D$}, 
if there is no confusion.

For a simple arc $X$, we set\\
\ \ \ $
\begin{array}{rl}
\partial X&= \text{ the set of its two endpoints, 
and}\\
{\rm Int}~X&=~X-\partial X.
\end{array}
$

Let $\Gamma$ be a chart. 
An edge of $\Gamma$ is called 
a {\it terminal edge}
if it contains a white vertex and 
a black vertex.

\begin{figure}
\begin{center}
\includegraphics{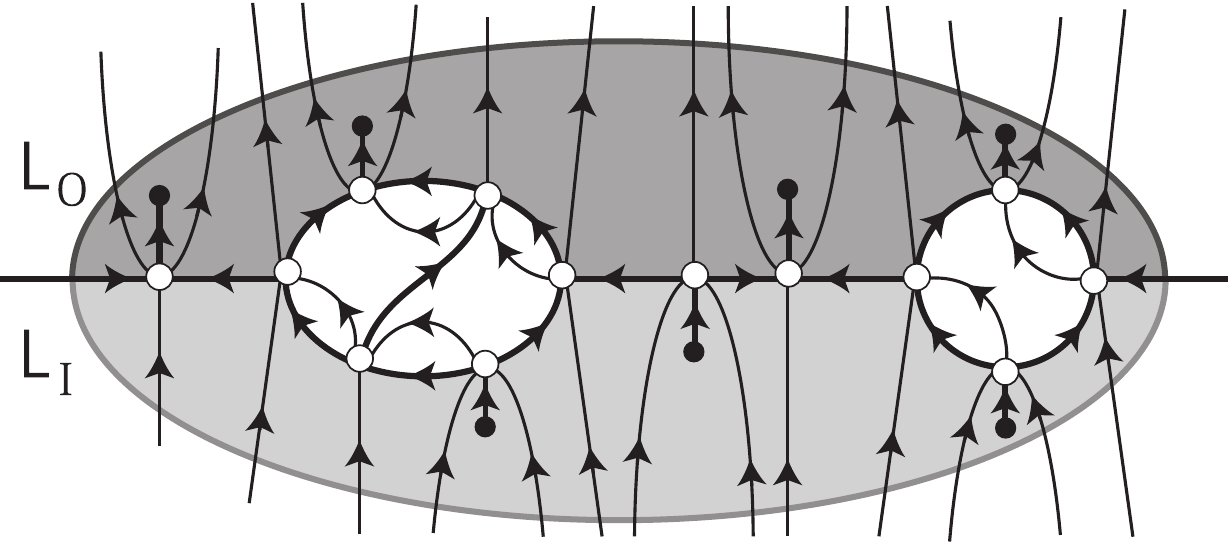}
\end{center}
\caption{\label{fig04}
An IO-tangle of label $m$. 
The thick edges are of label $m$. 
The light gray arc is $L_I$,
and the dark gray arc is $L_O$.
}
\end{figure}

Let $\Gamma$ be a chart, and 
$m$ a label of the chart.
A tangle $(\Gamma\cap D,D)$
is called a 
{\it IO-tangle of label $m$}
provided that 
(see Fig. \ref{fig04})
\begin{enumerate}
\item[(i)] $\partial D$ intersects neither terminal edge nor free edge,
\item[(ii)] 
there exists a label $k$ with $|m-k|=1$ and
$\Gamma\cap D\subset\Gamma_m\cup\Gamma_k$,
\item[(iii)]
there exist two arcs $L_I,L_O$ on $\partial D$
with $L_I\cap L_O=\partial L_I
=\partial L_O=\Gamma_m\cap \partial D$,
\item[(iv)] 
for any point 
$p\in\Gamma\cap$~Int~$L_I$, 
there is an edge 
of label $k$
locally inward at $p$, and\\
for any point 
$p\in \Gamma\cap$~Int~$L_O$, 
there is an edge 
of label $k$ 
locally outward at $p$.
\end{enumerate}
An IO-tangle of label $m$ is 
said to be {\it simple} 
if all the terminal edge in $D$ is of label $m$.

\begin{theorem}
\label{TwoColorGateTangle} 
If $(\Gamma\cap D,D)$ is 
an admissible tangle 
in a minimal chart $\Gamma$ such that 
\begin{enumerate}
\item[{\rm $($a$)$}] 
$\Gamma\cap D\subset \Gamma_m\cup\Gamma_{m-1}$ 
or 
$\Gamma\cap D\subset \Gamma_m\cup\Gamma_{m+1}$ 
for some label $m$,
\item[{\rm $($b$)$}] 
$\Gamma_m\cap \partial D$ 
consists of exactly two points, 
and 
\item[{\rm $($c$)$}] 
$\Gamma_m\cap D$ contains a cycle, 
\end{enumerate}
then
the tangle $(\Gamma\cap D,D)$ is 
a simple IO-tangle of label $m$.
\end{theorem}

\begin{figure}
\begin{center}
\includegraphics{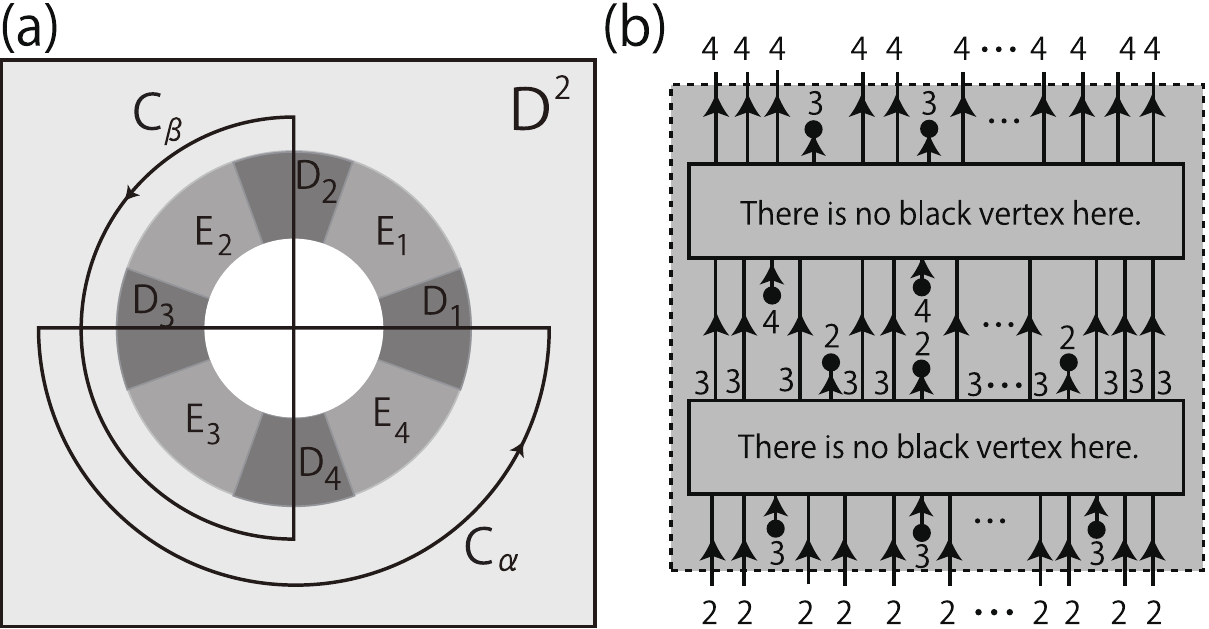}
\end{center}
\caption{\label{fig05}
(b) A tangle $(\Gamma\cap E_i,E_i)$ with $\Gamma\cap E_i\subset \Gamma_2\cup \Gamma_3\cup \Gamma_4$
for the case $\alpha=1$ and $\beta=5$,
here all the free edges and simple hoops are in a regular neighbourhood of $\partial D^2$.
}
\end{figure}

The theorem above determines 
the structure of minimal charts with two crossings
(\cite{TwoCrossingStructureI} and 
\cite{TwoCrossingStructureII}),
and give us an enumeration of the charts 
with two crossings.
The enumeration is much complicated 
than the one of 
$2$-bridge links in ${\Bbb R}^3$, of course. 
We enumerate charts with two crossings as follows
(see \cite{TwoCrossingStructureI},
\cite{TwoCrossingStructureII}):
For any minimal $n$-chart $\Gamma$ 
with two crossings in a disk $D^2$,
there exist two labels 
$1\le \alpha<\beta\le n-1$ such that 
 $\Gamma_\alpha$ and $\Gamma_\beta$
contain cycles $C_\alpha$ and $C_\beta$
with $C_\alpha\cap C_\beta$ the two crossings
and that
for any label $k$ with $k<\alpha$ or $\beta<k$,
the set $\Gamma_k$ does not contain a white vertex.
If $\Gamma_\alpha$ or $\Gamma_\beta$
contains at least three white vertices,
then
after shifting all the free edges and simple hoops
into a regular neighbourhood of $\partial D^2$
by applying C-I-M1 moves and C-I-M2 moves, 
we can find 
an annulus $A$ 
containing all the white vertices of $\Gamma$ 
but not intersecting 
any hoops nor free edge
such that (see Fig.~\ref{fig05}(a))
\begin{enumerate}
\item[(1)] 
each connected component of $Cl(D^2-A)$ 
contains a crossing,
\item[(2)]
$\Gamma\cap \partial A=
(C_\alpha\cup C_\beta)\cap \partial A$, and 
$\Gamma\cap \partial A$ 
consists of eight points.
\end{enumerate}
We can show 
the annulus $A$ can be split 
into 
mutually disjoint four disks 
$D_1,D_2,D_3,D_4$ and 
mutually disjoint four disks 
$E_1,E_2,E_3,E_4$ such that
\begin{enumerate}
\item[(3)] 
for each $i=1,3$ (resp. $i=2,4$)
the tangle $(\Gamma\cap D_i,D_i)$ 
is an IO-tangle of label $\alpha$ 
(resp. label $\beta$),
\item[(4)] 
for each $i=1,2,3,4$, 
$(\Gamma\cap E_i,E_i)$
is a tangle with $\Gamma\cap E_i\subset \cup^{\beta-1}_{j=\alpha+1}\Gamma_j$
as shown in Fig.~\ref{fig05}(b).
\end{enumerate}
We count the number of edges between terminal edges 
in Fig.~\ref{fig05}(b) 
to enumerate charts with two crossings.
As important results, 
from the enumeration 
we can calculate the fundamental group of 
the exterior of the surface link represented by $\Gamma$, 
and the braid monodromy of the surface braid 
represented by $\Gamma$.

Our paper is organized as follows:
In Section~\ref{s:Prel}, 
we introduce 
the definition of C-moves and its related words.
In Section~\ref{s:AdmissibleConsecutiveTriplet}, 
we introduce Consecutive Triplet Lemma 
(Lemma~\ref{ConsecutiveTriplet}).
For a label $m$ of 
a minimal chart $\Gamma$, 
we also investigate 
a complementary domain $U$ 
of $\Gamma_m$ with 
$\Gamma\cap U\subset
\Gamma_{m-1}\cup\Gamma_{m+1}$.
In Section~\ref{s:PfTheorem1}, 
we give a proof of Theorem~\ref{Inequation}.
In Section~\ref{s:Dichromatic}, 
for a chart $\Gamma$ and a label $m$, 
we investigate a neighbourhood $N$ of 
an arc in $\Gamma_m$ with 
$N\cap\Gamma\subset\Gamma_m\cup\Gamma_k$ 
for some label $k$. 
This situation occurs  
an arc in the boundary of a 2-color disk and
an arc in an IO-tangle of label $m$.
In Section~\ref{s:SuspiciousCycle},
for a tangle $(\Gamma\cap D,D)$ and 
each label $m$ with $\Gamma_m\cap D\neq\emptyset$,
we obtain an equation and 
search for conditions for the existence 
of a special cycle which 
never bounds a $2$-color disk 
in a minimal chart.
In Section~\ref{s:PfTheorem2}, 
we give a proof of 
Theorem~\ref{NoNS-tangle}.
In Section~\ref{s:IO-Path}, 
we investigate the boundary of a 2-color disk.
In Section~\ref{s:Bridge}, 
we investigate an arc in $\Gamma_m$ such that 
each white vertex in the arc is contained in 
a terminal edge. 
This situation occurs an arc in 
$\Gamma_m$ connecting vertices in the boundary of 
$2$-color disks.
In Section~\ref{s:PfTheorem3}, 
we give a proof of 
Theorem~\ref{TwoColorGateTangle}.
In Section~\ref{s:Appendix}, 
we give complementary lemmata 
to make our paper self-contained.

\section{Preliminaries}
\label{s:Prel}

In this section we give the definitions of C-moves
 and its related words.

\begin{figure}
\begin{center}
\includegraphics{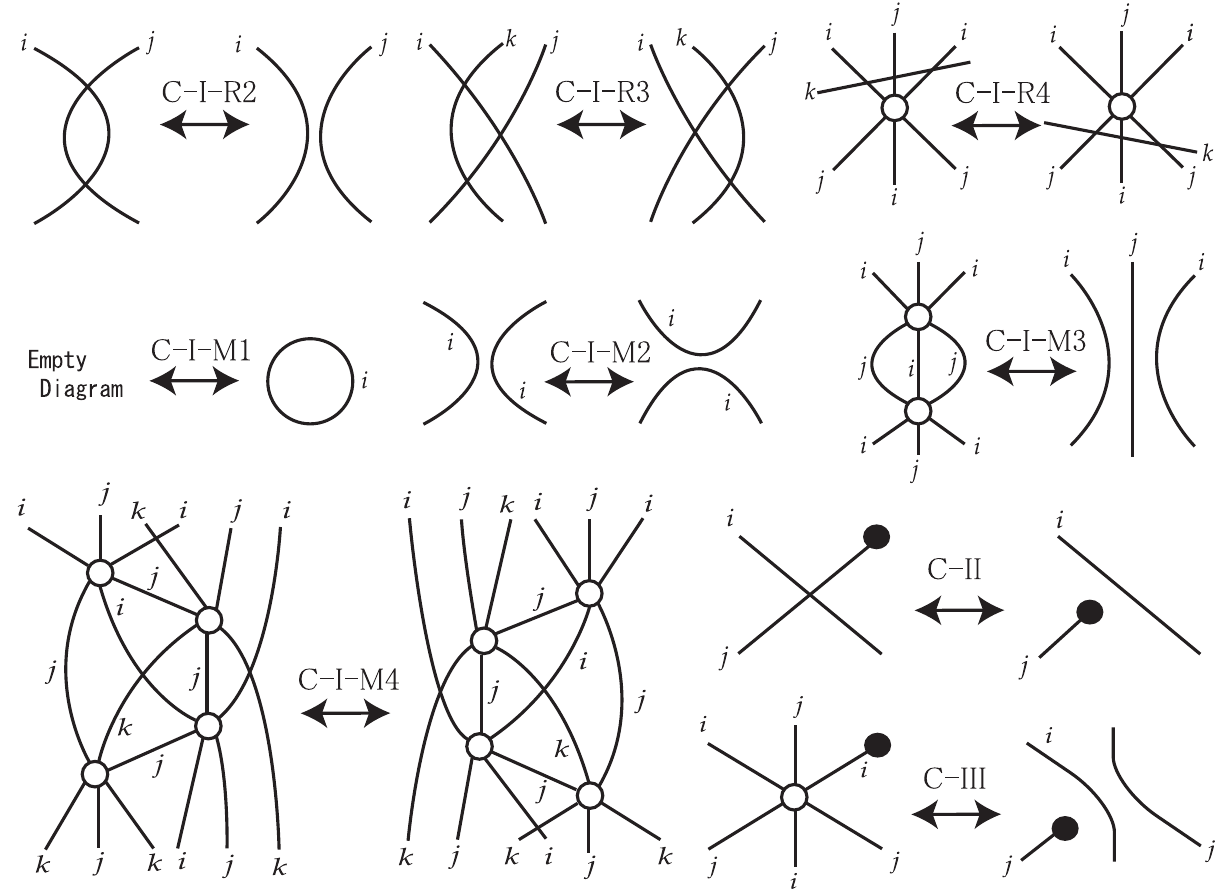}
\end{center}
\caption{ \label{fig06} 
For the C-III move, the edge containing 
the black vertex does not contain a middle arc at
a white vertex in the left figure. }
\end{figure}

Now {\it C-moves} are local modifications 
of charts as shown in Fig.~\ref{fig06}
(cf. \cite[p. 117]{KS}, 
\cite[p. 142--143]{BraidBook} and \cite{Tanaka}).
As one of C-moves,
Kamada originally defined CI-moves
as follows: 
A chart $\Gamma$ is obtained from
a chart $\Gamma'$ in a disk $D^2$
by a {\it CI-move},
if there exists a disk $E$ 
in $D^2$ such that 
\begin{enumerate}
\item[(i)] 
the two charts $\Gamma$ and $\Gamma'$
intersect the boundary of $E$ transversely
or
do not intersect the boundary of $E$, 
\item[(ii)] 
$\Gamma\cap E^c=\Gamma'\cap E^c$, and
\item[(iii)]
neither $\Gamma\cap E$ nor 
$\Gamma'\cap E$ 
contains a black vertex,
\end{enumerate}
where $E^c$ is 
the complement of $E$ in the disk $D^2$.

\begin{remark}
Any CI-move is realized by a finite sequence of 
seven types: C-I-R2, C-I-R3, C-I-R4, 
C-I-M1, C-I-M2, C-I-M3, C-I-M4.
\end{remark}

Let $D_1^2,D_2^2$ be disks,
and $pr_2:D_1^2\times D_2^2\to D_2^2$
the projection defined by $pr_2(x,y)=y$.
Let $Q_n$ be a set of $n$ interior points of $D_1^2$.
A {\it surface braid} $S$ is an oriented surface 
embedded properly in $D_1^2\times D_2^2$
such that
the map $pr_2|_S:S\to D_2^2$ is a branched covering of degree $n$
and $\partial S=Q_n\times\partial D_2^2$
\cite[Chapter 14]{BraidBook}.

\begin{figure}
\begin{center}
\includegraphics{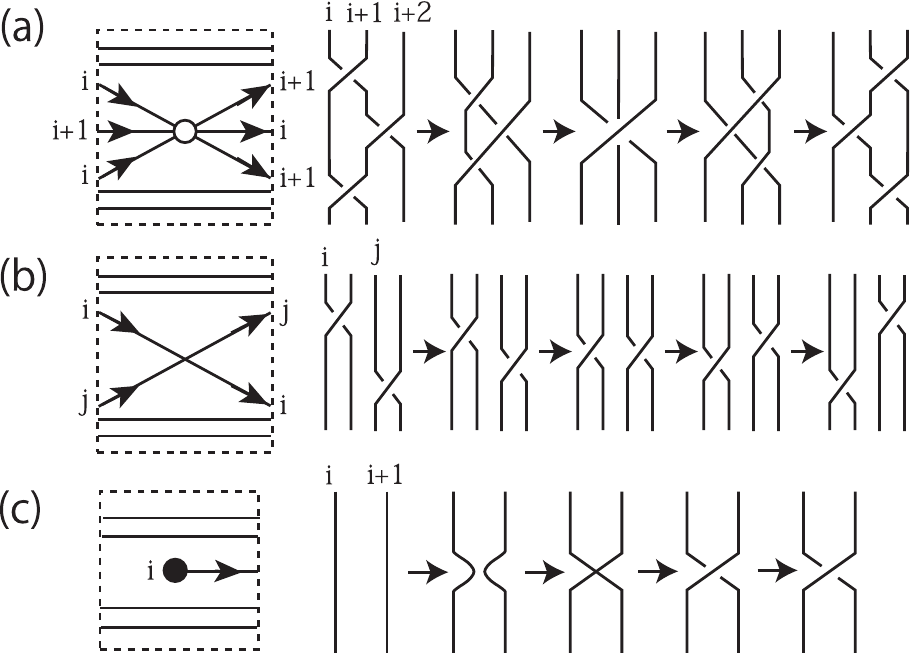}
\end{center}
\caption{ \label{fig07} 
(a) A motion picture for a white vertex.
(b) A motion picture for a crossing.
(c) A motion picture for a black vertex.}
\end{figure}

A surface braid $S$ can be represented by 
a motion picture method as follows:
We identify the disk $D_2^2$ the product of 
the unit intervals $I_3,I_4$.
For each $t\in I_4=[0,1]$,
we define the subset $b_t$ in $D_1^2\times I_3$
by $b_t\times\{t\}=S\cap (D_1^2\times I_3\times\{t\})$
such that $b_t$ is a geometric $n$-braid in 
$D_1^2\times I_3$ except for 
a finite number of values $t_1,t_2,\cdots,t_m\in I_4$. 
Thus the surface braid $S$ is represented by 
the one-parameter family $\{b_t\}_{t\in[0,1]}$
called a {\it motion picture}
 (see Fig.~\ref{fig07}) 
(cf. \cite{BraidBook}).

\begin{figure}
\begin{center}
\includegraphics{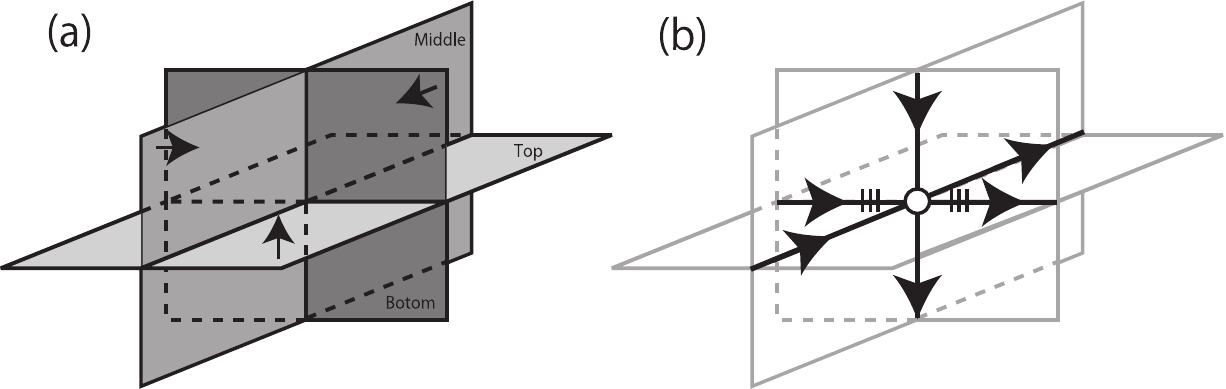}
\end{center}
\caption{ \label{fig08} }
\end{figure}

A chart can be constructed from a surface braid 
as follows:
We identify the disk $D_1^2$ the product of 
the unit intervals $I_1,I_2$.
Let $\pi:D^2_1\times D_2^2=I_1\times I_2\times D_2^2\to
 I_1\times D_2^2$ be the projection 
defined by $\pi(x_1,x_2,y)=(x_1,y)$.
Then the image $\pi(S)$ of a surface braid $S$
has double points, triple points and branch points.
The union of these points is a graph $G$ 
in $I_1\times D_2^2$.
We construct a chart from the surface braid $S$ 
by projecting the graph $G$ into $D_2^2$.
The orientation of each edge is determined as follows:
In a neighbourhood of each double point,
there are two sheets; an upper sheet and a lower sheet.
Let $v_1,v_2$ be the normal vectors at the double point 
in an upper sheet and a lower sheet 
determined from the orientation of $S$ respectively.
Define the tangent vector $v_3$ of the graph $G$ 
at the double point such that
the orientation of $I_1\times D_2^2$ matches 
the triplet $(v_1,v_2,v_3)$.
Then we have the orientated graph $G$;
 the chart each of whose edge is oriented.
In a neighbourhood of each triple point,
there are three sheets;
a top sheet, a middle sheet and a bottom sheet.
A white vertex is corresponding to a triple point,
and middle arcs are recognizable as the intersection 
of a top sheet and a bottom sheet (see Fig.~\ref{fig08}).
The label of each edge is determined as follows:
Let $\sigma_1,\sigma_2,\cdots,\sigma_{n-1}$ be 
the standard generators of 
the classical braid group $B_n$.
Any geometric braid $b$ can be represented by a product
of $\sigma_1,\sigma_2,\cdots,\sigma_{n-1}$ and their inverse $\sigma_1^{-1},\sigma_2^{-1},\cdots,\sigma_{n-1}^{-1}$. A double point in an edge of the graph $G$
is corresponding to a crossing
of a geometric braid $b_t$ for some $t\in I_4$,
and the crossing is corresponding to
some generator $\sigma_i$ or $\sigma_i^{-1}$.
Then we define the label of the edge of $G$ 
by the label $i$.
Thus we have the labeled graph $G$;
 the chart each of whose edge is labeled.

Now for any chart $\Gamma$ 
a disk $D^2$ is assigned 
so that the chart is contained in the disk $D^2$ 
by the definition of charts. 
If we want to emphasize that 
a domain $X$ is contained in Int~$D^2$, 
then 
we say that $X$ is {\it finite}.

\begin{figure}
\begin{center}
\includegraphics{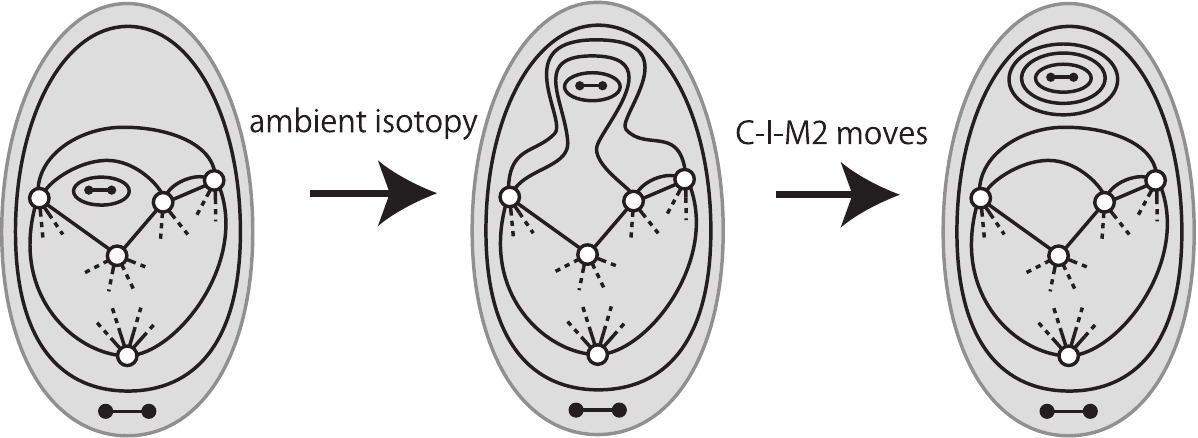}
\end{center}
\caption{ \label{fig09} }
\end{figure}

For any chart in a disk $D^2$
we can move free edges and simple hoops into 
a regular neighbourhood of $\partial D^2$ 
by C-I-M2 moves and ambient isotopies of $D^2$
as shown in Fig.~\ref{fig09}.
Even during argument,
if free edges or simple hoops appear, 
we immediately move them 
into a regular neighbourhood of $\partial D^2$.
Thus we assume the following.
\begin{assumption}
\label{AssumptionFreeEdge}
{\it For any chart in a disk $D^2$, 
all the free edges and 
simple hoops  
are in a regular neighbourhood of 
$\partial D^2$.}
\end{assumption}

Let $\Gamma$ be a chart in a disk $D^2$, and 
$X$ the union of all the free edges 
and simple hoops.
By Assumption~\ref{AssumptionFreeEdge},
the set $X$ is in a regular neighbourhood $N$ of 
$\partial D^2$ in $D^2$. 
Define
$${\rm Main}(\Gamma)=\Gamma-X.$$
Let $\widehat D=Cl(D^2-N)$. 
Then $\Gamma\cap\widehat D=$Main$(\Gamma)$.
Hence $(\Gamma\cap\widehat D,\widehat D)$ is 
a tangle without free edges 
and simple hoops.
\begin{assumption}
\label{AssumptionFreeEdgeSimpleHoop}
{\it In this paper, our arguments are done 
in the disk $\widehat D$, 
otherwise mentioned.}
\end{assumption}

\begin{figure}
\begin{center}
\includegraphics{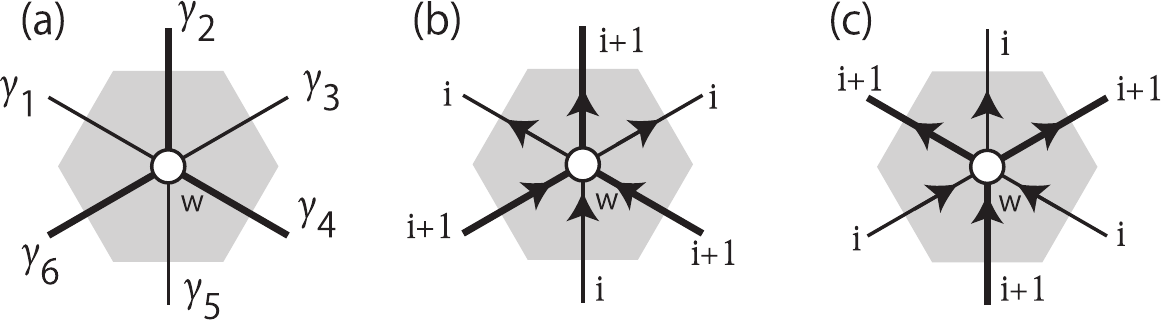}
\end{center}
\caption{ \label{fig10} }
\end{figure}

\begin{remark}
\label{EdgesAroundVertex}
{\rm Let $\Gamma$ be a chart.
Let $\gamma_1,\gamma_2,\cdots,\gamma_6$ be six arcs 
around a white vertex $w$ lying
clockwise in this order (see Fig.~\ref{fig10}(a)).
Then we have the following 
(see Fig.~\ref{fig10}(b) and (c)).
\begin{enumerate}
\item[$(1)$]
For each $j=1,2,\cdots,6$,
one of the two arcs $\gamma_j,\gamma_{j+1}$ 
is not a middle arc at $w$.
\item[$(2)$] For each $j=1,2$,
one of the three arcs $\gamma_j,\gamma_{j+2},\gamma_{j+4}$ 
of the same label is middle at $w$ but 
the others are not middle at $w$.
\end{enumerate}
Here $\gamma_{j+6}=\gamma_j$ for each integer $j$.}
\end{remark}

Let $m$ be a label of a chart $\Gamma$.
A simple closed curve in $\Gamma_m$ 
is called a {\it ring}, 
if it contains a crossing 
but does not contain 
a white vertex nor a black vertex.
An arc is said to be {\it internal} 
if it is contained in an internal edge or a ring.

\begin{figure}
\begin{center}
\includegraphics{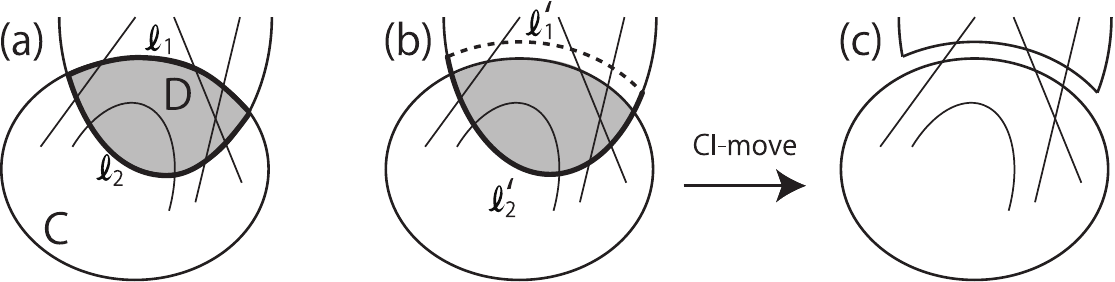}
\end{center}
\caption{ \label{fig11} 
(a) The gray area is the disk $D$. 
(b) The dotted arc is $\ell'_1$ and 
the thick arc is $\ell'_2$.}
\end{figure}

\begin{remark}
\label{Assumption0}
{\rm Let $\Gamma$ be a minimal chart. 
Then we have 
the following:
\begin{enumerate}
\item[(1)]
{\it If an edge of $\Gamma$ contains a black vertex, 
then it is a terminal edge or a free edge.}
For, if the edge contains a crossing, 
then we can eliminate  the crossing 
on the edge by 
a C-II move. 
This contradicts that the chart is minimal.
\item[(2)]
{\it Any terminal edge of $\Gamma$ contains a middle arc 
at its white vertex.}
For, if not, 
we can eliminate the white vertex by a C-III move.
\item[(3)]
{\it Each complementary domain of
any ring must contain 
at least one white vertex} 
(cf. \cite[Assumption 4]{ChartApp1}). 
For, 
suppose that there exists a ring $C$
such that a complementary domain of $C$ 
does not contain any white vertices.
Let $F$ be the closure of the complementary domain.
By (1), 
the ring $C$ does not intersect any terminal edge 
nor free edge. 
Thus any crossing on $C$ 
is contained in a proper internal arc of $F$. 
Since $F$ is a disk or an annulus, 
the domain $F$ contains a disk $D$ bounded by 
an arc $\ell_1$ on $C$ and 
a proper internal arc $\ell_2$ of $F$  
such that any crossing on Int $\ell_1$ is
contained in a proper internal arc $\ell$ of $D$ 
intersecting $\ell_2$ by a crossing
(see Fig.~\ref{fig11}(a)). 
Let $\ell_2'$ 
be an internal arc 
with $\ell'_2\supset\ell_2$ 
such that 
$\ell_2'-\ell_2$ does not contain a crossing.
Let 
$\ell_1'$ be an arc 
outside $F$ parallel to $\ell_1$ 
with $\partial\ell_1'=\partial\ell_2'$
(see Fig.~\ref{fig11}(b)). 
Thanks to Assumption~\ref{AssumptionFreeEdge} and Assumption~\ref{AssumptionFreeEdgeSimpleHoop}, 
there does not exist any black vertex 
in the disk bounded by 
$\ell_1'\cup\ell_2'$. 
Thus 
we can shift the arc $\ell_2'$ 
to the arc $\ell_1'$
by a CI-move so that 
the number of crossings decreases at least two 
(see Fig.~\ref{fig11}(c)). 
This contradicts that the chart is minimal. 
\end{enumerate}
}
\end{remark}

\section{Admissible consecutive triplets}
\label{s:AdmissibleConsecutiveTriplet}
In this section we introduce 
one of our main tools, 
Consecutive Triplet Lemma, 
which minimal charts satisfy. 
Also an original idea for 
Theorem~\ref{Inequation} is given in 
Lemma~\ref{NumberMiddleEdge}.

Let $\Gamma$ be a chart,
and $m$ a label of $\Gamma$.
An internal edge of label $m$ is 
called a {\it loop}
if it contains exactly one white vertex.

\begin{figure}
\begin{center}
\includegraphics{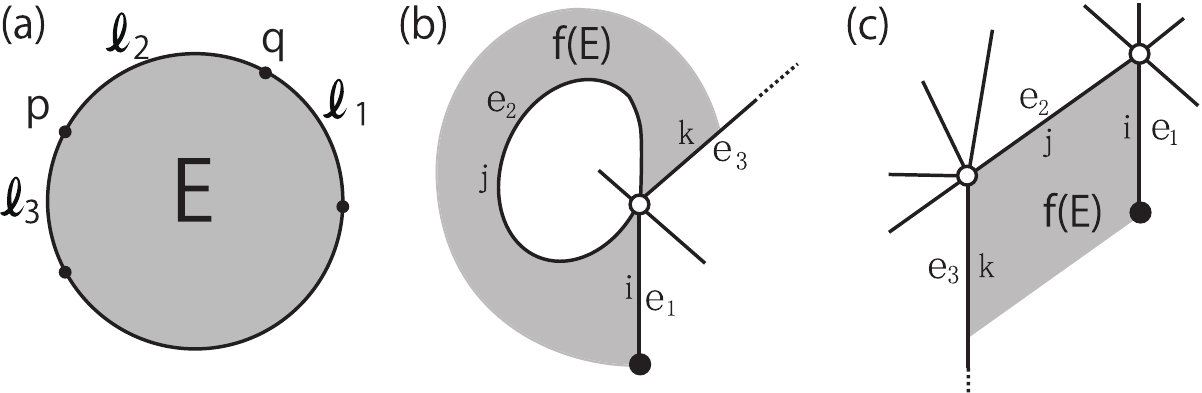}
\end{center}
\caption{ \label{fig12} }
\end{figure}

Let $E$ be a disk, and
$\ell_1,\ell_2,\ell_3$ three arcs on $\partial E$
such that each of $\ell_1\cap \ell_2$ and $\ell_2\cap \ell_3$ is one point and $\ell_1\cap \ell_3=\emptyset$
(see Fig.~\ref{fig12}(a)),
say $p=\ell_1\cap \ell_2$,
$q=\ell_2\cap \ell_3$.
Let $\Gamma$ be a chart in a disk $D^2$.
Let $e_1$ be a terminal edge of 
 $\Gamma$. 
A triplet $(e_1,e_2,e_3)$ of 
mutually different edges of $\Gamma$
is called 
a {\it consecutive triplet}
if there exists
a continuous map $f$ from the disk $E$ 
to the disk $D^2$ such that (see Fig.~\ref{fig12}(b) and (c))
\begin{enumerate}
\item[(i)] the map $f$ is injective on $E-\{p,q\}$,
\item[(ii)] 
$f(\ell_3)$ is an arc in $e_3$, and $f({\rm Int}~E)\cap\Gamma=\emptyset$,
$f(\ell_1)=e_1$,
$f(\ell_2)=e_2$,
\item[(iii)]
$f(p)$ and $f(q)$ are white vertices.
\end{enumerate}
If the label of $e_3$ is different
from the one of $e_1$ 
then the consecutive triplet is said to be
{\it admissible}.

\begin{lemma}[Consecutive Triplet Lemma]
\label{ConsecutiveTriplet}
{\rm $($\cite[Lemma 1.1]{OneCrossing}$)$}
{\it Any consecutive triplet 
in a minimal chart is admissible.}
\end{lemma}

The above lemma was proven by 
the maximality of bigons in a minimal chart
\cite{OneCrossing}. 
In Section~\ref{s:Appendix}, 
we shall give a proof of Consecutive Triplet Lemma 
to make our paper to be self-contained.

In the proof of Lemma~\ref{NumberMiddleEdge}, 
we must be careful of the next remark.

\begin{remark}
\label{NoteBiMiddle} 
{\rm 
Let $w$ be the white vertex of a loop.
In a small neighborhood of the white vertex $w$,
the loop contains two short arcs $\gamma,\gamma'$ with $\gamma\cap\gamma'=w$. 
One of the two arcs $\gamma,\gamma'$
is a middle arc at $w$,
but the other is not a middle arc at $w$.
}
\end{remark}

Let $\Gamma$ be a chart.
Let $\gamma_1,\gamma_2,\cdots,\gamma_6$ be 
six short arcs 
around a white vertex $w$ lying
clockwise in this order (see Fig.~\ref{fig10}(a)).
For each $j=1,2,\cdots,6$
let $e_j$ be the edge containing $\gamma_j$ 
(possibly $e_j=e_{j+2}$), 
here $e_{j+6}=e_j$ for each integer $j$. 
Let $e=e_k$ for some $k\in\{1,2,\cdots,6\}$.
Then the two edges $e_{k-1},e_{k+1}$ are 
called the $(e,w)$-{\it edges}.


\begin{lemma}
\label{NumberMiddleEdge}
Let $\Gamma$ be a minimal chart,
and $m$ a label of $\Gamma$.
Let $U$ be a finite complementary domain 
of $\Gamma_m$ with
$\Gamma\cap U\subset \Gamma_{m-1}\cup \Gamma_{m+1}$.
Suppose that 
$U$ does not contain any crossing.
Then we have the following:
\begin{enumerate}
\item[{\rm $($a$)$}] The component $U$ does 
not contain any white vertex.
\item[{\rm $($b$)$}] $Cl(U)$ does not contain 
any terminal edge of label $m$.
\item[{\rm $($c$)$}]  If $U$ is an open disk
and if $Cl(U)$ contains a white vertex in $\Gamma_m$,
then there exist at least two middle arcs 
of label $m\pm1$ in $Cl(U)$.
\end{enumerate}
\end{lemma}

\begin{figure}[htb]
\begin{center}
\includegraphics{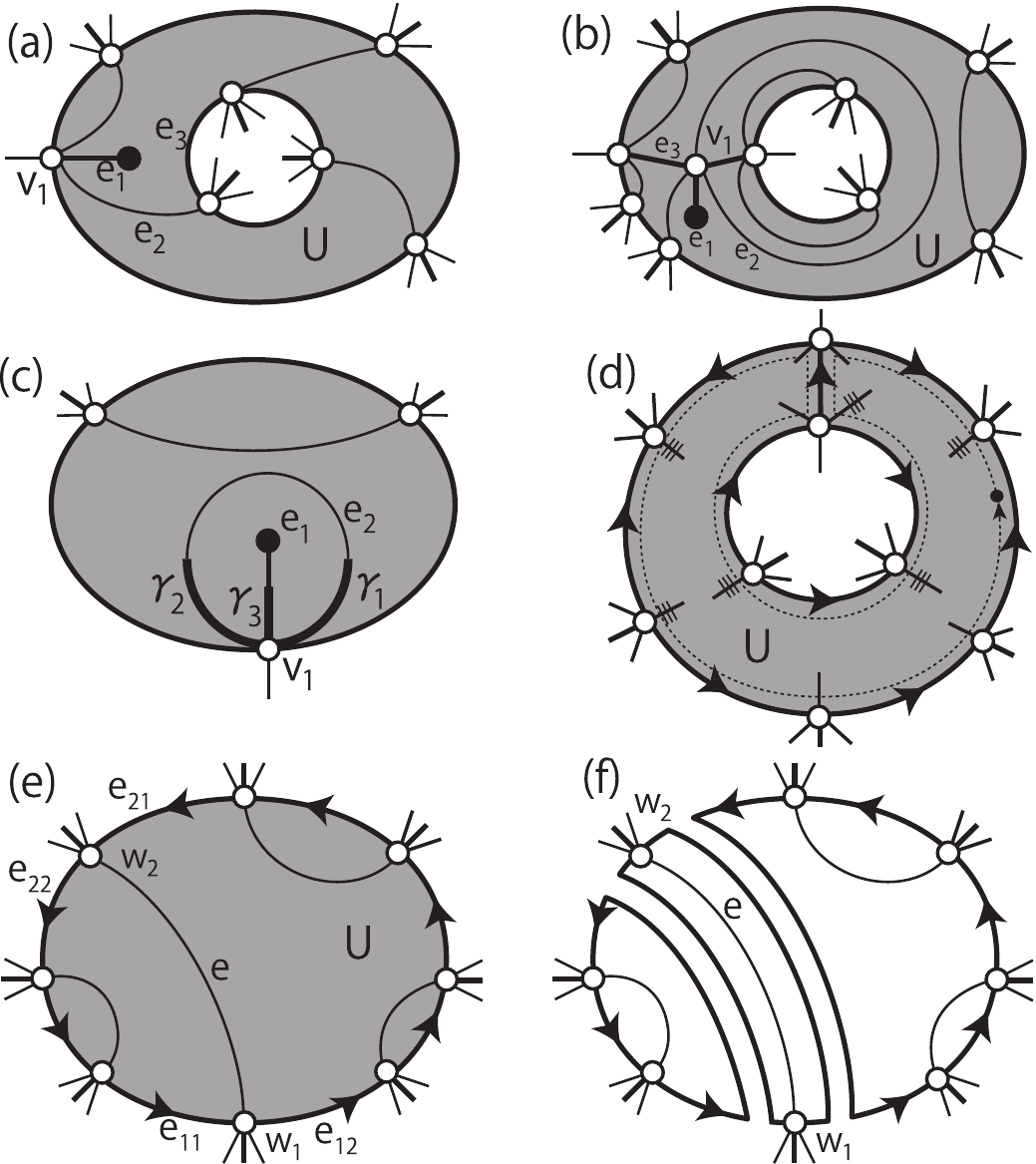}
\end{center}
\caption{ \label{fig13} 
In (a),(b),(c),(d),(e), each gray region is 
a finite complementary domain  $U$ of $\Gamma_m$.
In (c), the thick lines are the short arcs 
$\gamma_1,\gamma_2,\gamma_3$.
In (d), each arc with three transversal short arcs 
is a middle arc in $Cl(U)$ of label $m\pm 1$.}
\end{figure}

\begin{proof}
Since there is no white vertex contained 
in $\Gamma_{m-1}\cap\Gamma_{m+1}$, 
Statement (a) holds.

We show Statement (b). Suppose that 
$Cl(U)$ contains a terminal edge $e_1$ 
of label $m$.
Let $v_1$ be the white vertex of $e_1$. 
Then $v_1\in \partial U$ by Statement (a).
Let $e_2$ be an $(e_1,v_1)$-edge.
Then $e_2$ is of label $m\pm1$ and 
$e_2\cap U\not=\emptyset$.
The edge $e_2$ is not a terminal edge.
For, if $e_2$ is a terminal edge, 
then by Remark~\ref{EdgesAroundVertex}$(1)$
one of terminal edges $e_1$ and $e_2$ 
does not contain a middle arc at $v_1$.
This contradicts 
Remark~\ref{Assumption0}$(2)$.

If $e_2$ is not a loop, 
then there exists 
an edge $e_3$ of label $m$ 
in $\partial U$ such that 
the consecutive triplet $(e_1,e_2,e_3)$ 
is not admissible 
(see Fig.~\ref{fig13}(a)).
This contradicts Consecutive Triplet Lemma 
(Lemma~\ref{ConsecutiveTriplet}).

Suppose that $e_2$ is a loop.
In a regular neighborhood $N$ of $v_1$,
the edge $e_2$ contains 
two short arcs $\gamma_1,\gamma_2$
with $\gamma_1\cap\gamma_2=v_1$.
Now $e_1$ contains 
a short arc $\gamma_3$ 
in $N$ with $\gamma_3\ni v_1$.
If $\gamma_1,\gamma_2,\gamma_3$ are not
consecutive 
around $v_1$ 
(see Fig.~\ref{fig13}(b)),
then $v_1\in$~Int$(Cl(U))$.
Hence again 
there exists an edge $e_3$ of label $m$ 
in $\partial U$ such that 
the consecutive triplet $(e_1,e_2,e_3)$ 
is admissible 
(see Fig.~\ref{fig13}(b)).
This contradicts 
Consecutive Triplet Lemma 
(Lemma~\ref{ConsecutiveTriplet}).

Suppose that $\gamma_1,\gamma_2,\gamma_3$ are 
consecutive around $v_1$.
Then 
$\gamma_1,\gamma_3,\gamma_2$ are 
consecutive arcs
situated around $v_1$ in this order 
(see Fig.~\ref{fig13}(c)). 
By Remark~\ref{Assumption0}$(2)$,
the arc $\gamma_3$ is middle at $v_1$.
Since $e_2$ is a loop,
by Remark~\ref{NoteBiMiddle}
 one of $\gamma_1,\gamma_2$ is 
middle at $v_1$, say $\gamma_1$.
Then the two consecutive arcs $\gamma_1,\gamma_3$ are 
middle at $v_1$.
This contradicts Remark~\ref{EdgesAroundVertex}$(1)$.
Hence Statement (b) holds.

We show Statement (c). 
Let
 $\mathcal W=\{v~|~v$ is a white vertex 
contained in a middle arc 
of label $m\pm 1$ in $Cl(U)\}.$ 
We shall show $\mathcal W\neq\emptyset$.
Suppose $\mathcal W=\emptyset$.
By the assumption, $Cl(U)$ contains 
a white vertex $w$ in $\Gamma_m$.
Then $w\in \partial U=Cl(U)-U$ by 
Statement (a).
Let $e$ be an edge of label $m\pm 1$
intersecting $U$ and 
containing the white vertex $w$. 
Since $\mathcal W=\emptyset$,
the edge $e$ does not contain a middle arc at $w$.
By Remark~\ref{Assumption0}$(2)$,
the edge $e$ is not a terminal edge.

We claim that the edge $e$ is not a loop.
For, suppose that the edge $e$ is a loop.
Since there is no crossing 
in $\Gamma_{m\pm1}\cap\Gamma_m$,
the edge $e$ must be contained in $Cl(U)$.
By Remark~\ref{NoteBiMiddle},
the edge $e$ contains a middle arc at $w$.
Thus $w\in \mathcal W$.
This contradicts $\mathcal W=\emptyset$. 
Hence the edge $e$ is not a loop.

Thus
the edge $e$ contains two white vertices 
in $\partial U$, say $w_1=w$ and $w_2$.
Then for $i=1,2$, 
there exist two edges 
$e_{i1}$ and $e_{i2}$ 
of label $m$ in $\partial U$
containing $w_i$.
Go around $\partial U$ starting from $e_{11}$
and next pass through $e_{12}$ 
and so on.
Without loss of generality we can assume that
$e_{11}$ is oriented inward at $w_1$ 
and that
we pass $e_{11},e_{12},\cdots,e_{21},
e_{22}$ in this order.
Since $Cl(U)$ does not contain 
a terminal edge of label $m$ 
by Statement (b), 
when we go around 
$\partial U$ 
(see Fig.~\ref{fig13}(d)), 
we have that

$(1)$ the orientation of the edge in 
$\partial U$ 
changes at a vertex in $\mathcal W$.\\
Namely, 
on the way going around 
$\partial U$, 
if an edge and the next edge are
inward (or outward) at a white vertex,
then the white vertex is contained 
in $\mathcal W$.
Since $\mathcal W=\emptyset$,
Statement $(1)$ implies that
the edge $e_{12}$ is 
oriented outward at $w_1$,
the edge $e_{21}$ is 
oriented inward at $w_2$, and
the edge $e_{22}$ is 
oriented outward at $w_2$
(see Fig.~\ref{fig13}(e)).
Thus applying two C-I-M2 moves
between $e_{11}$ and $e_{22}$ 
and
between $e_{12}$ and $e_{21}$
(see Fig.~\ref{fig13}(f)),
we can eliminate the two white vertices
$w_1$ and $w_2$ by a C-I-M3 move.
This contradicts the fact that
the chart $\Gamma$ is minimal.
Thus $\mathcal W\neq\emptyset$.

Again Statement $(1)$ implies that 
the set $\mathcal W$ consists of exactly 
even number of vertices.
Thus there exist at least two middle arcs 
of label $m\pm 1$ in $Cl(U)$.
\end{proof}

Lemma~\ref{NumberMiddleEdge}(c) 
is our start point for 
Theorem~\ref{Inequation}.

\begin{cor}
\label{Cor3ColorDiskNoTerminalEdge}
Let $\Gamma$ be a minimal chart,
and $m$ a label of $\Gamma$.
Then we have the following:
\begin{enumerate}
\item[{\rm $($a$)$}] 
If $E$ is a $3$-color disk with
$\partial E\subset\Gamma_m$, 
then $E$ does not contain any terminal edge of label $m$.
\item[{\rm $($b$)$}] 
If $E$ is a $2$-color disk with
$\partial E\subset\Gamma_m$, 
then $E$ does not contain any terminal edge.
\end{enumerate}
\end{cor}

\begin{proof}
We show Statement (a).
Since the $3$-color disk $E$ does not contain any crossing, each component of $E-\Gamma_m$ does not contain any crossing.
By Lemma~\ref{NumberMiddleEdge}(b),
the closure of each component of $E-\Gamma_m$ does not contain any terminal edge of label $m$,
and so does $E$.
Thus Statement (a) holds.

We show Statement (b).
Since a $2$-color disk is a $3$-color disk,
the $2$-color disk $E$ does not contain any terminal edge of label $m$ by Statement (a).

Let $k$ be the label with 
$|m-k|=1$ and  
$\Gamma\cap E\subset\Gamma_m\cup\Gamma_k$.
Suppose that there exists a terminal edge of label $k$ 
in $E$.
Since any white vertex in $\Gamma\cap E$
is contained in $\Gamma_m\cap \Gamma_k$,
we can find a non-admissible consecutive triplet
by a similar way to 
Lemma~\ref{NumberMiddleEdge}(b).
This contradicts 
Consecutive Triplet Lemma 
(Lemma~\ref{ConsecutiveTriplet}).
Hence $E$ does not contain any terminal edge of label $k$. Therefore $E$ does not contain any terminal edge.
\end{proof}


\section{A proof of Theorem~1.1}
\label{s:PfTheorem1}

In this section
we shall prove Theorem~\ref{Inequation}
by using Lemma~\ref{NumberMiddleEdge}(c), 
Corollary~\ref{Cor3ColorDiskNoTerminalEdge}(a),
Lemma~\ref{3ColorDisk} and Lemma~\ref{DecompositionM}.

Let $\Gamma$ be a chart, and
$m$ a label of $\Gamma$.
A simple closed curve in $\Gamma_m$
is called a {\it cycle of label $m$}.
{\it We consider hoops, rings and loops as cycles.}

Let $\Gamma$ be a chart, and $m$ a label of $\Gamma$.
Let $C$ be a cycle of label $m$ bounding a disk $E$.
An edge $e$ of label $m$ 
is called
{\it an inside $($resp. outside$)$ edge for $C$} 
provided that 
\begin{enumerate}
\item[(i)]
$e\cap C$ consists of 
one white vertex or two white vertices, and
\item[(ii)]
$e\subset E$ (resp. $e\subset Cl(E^c)$).
\end{enumerate}
For a cycle $C$ of label $m$ 
bounding a disk $E$, 
we define\\
$\begin{array}{ll}
\mathcal{W}(E)&
\hspace{-2.5mm}= \{ w \ | 
\text{ $w$ is a white vertex in $E$} \},\\
{\mathcal{W}}({\rm Int}~E)&
\hspace{-2.5mm}= \{ w \ | 
\text{ $w$ is a white vertex in Int $E$} \},\\
\mathcal{W}(C)&
\hspace{-2.5mm}=\{w\ | 
\text{ $w$ is a white vertex in $C$}\},\\
\mathcal{W}_I(C,m)&
\hspace{-2.5mm}= \{ w\in \mathcal{W}(C) \ | 
\text{ $w$ is contained in an inside edge 
for $C$} \},\\
\mathcal{W}_I^{{\rm Mid}}(C,m)&
\hspace{-2.5mm}=\{w\in\mathcal{W}(C)\ |
\text{ there exists an inside edge 
for $C$ middle at $w$} \},\\
\mathcal{W}_O(C,m)&
\hspace{-2.5mm}= \{ w\in \mathcal{W}(C) \ | 
\text{ $w$ is contained in an outside edge for $C$} \},\\
\mathcal{W}_O^{{\rm Mid}}(C,m)&
\hspace{-2.5mm}=\{w\in\mathcal{W}(C)\ |
\text{ there exists an outside edge for $C$ middle at $w$} \}.\\
\end{array}$

\begin{remark}
\label{O(C)I(C)} 
{\rm 
Let $\Gamma$ be a chart, and $m$ a label of $\Gamma$. 
For a cycle $C$ of label $m$ 
bounding a disk $E$, 
we have the following:
\begin{enumerate}
\item[$(1)$] 
Let $w$ be a white vertex
in $\mathcal{W}_O(C,m)$.
Then the vertex $w$ is 
in $\mathcal{W}_O^{{\rm Mid}}(C,m)$ 
if and only if 
there exists an arc of label $m\pm1$ 
in $E$ middle at the vertex $w$.
\item[$(2)$] 
Let $w$ be a white vertex 
in $\mathcal{W}_I(C,m)$.
Then the vertex $w$ is {\it not} 
in $\mathcal{W}_I^{{\rm Mid}}(C,m)$
if and only if 
there exists an arc of label $m\pm1$ 
in $E$ middle at the vertex $w$.
\item[$(3)$]
Let $w$ be a white vertex
in $\mathcal{W}_I(C,m)$.
Then the vertex $w$ is 
in $\mathcal{W}_I^{{\rm Mid}}(C,m)$ 
if and only if
 there exists an arc of label $m\pm1$ 
 in $Cl(E^c)$ middle at the vertex $w$.
\item[$(4)$] 
The set $\mathcal{W}(C)$
splits into disjoint subsets 
$\mathcal{W}_O(C,m)$ and 
$\mathcal{W}_I(C,m)$.
\item[$(5)$] 
The set $\mathcal{W}(E)$
splits into three mutually disjoint subsets 
$\mathcal{W}({\rm Int}~E)$,
$\mathcal{W}_O(C,m)$ and 
$\mathcal{W}_I(C,m)$.
\end{enumerate}
}
\end{remark}

\begin{lemma}\label{3ColorDisk}
Let $\Gamma$ be a minimal chart, 
and $m$ a label of $\Gamma$.
Let $E$ be a $3$-color disk 
with $\partial E\subset \Gamma_m$
but without free edges nor simple hoops. 
If $\Gamma_m\cap E$ is connected, 
then $E$ contains neither hoop nor ring.
\end{lemma}

\begin{proof}
Suppose that $E$ contains a hoop. 
By the assumption,
the hoop is not simple. 
Thus there exists a white vertex $v$
in the interior of the disk 
bounded by the hoop.
Any white vertex in $E$ is in 
$\Gamma_m\cap\Gamma_{m-1}$ or 
$\Gamma_m\cap\Gamma_{m+1}$.
Thus 
$\Gamma_m\cap E$ contains 
at least two connected components; 
one containing the vertex $v$,
and the other containing $\partial E$.
Hence $\Gamma_m\cap E$ is not connected.
This is a contradiction.
Thus $E$ does not contain any hoop.

Since any 3-color disk 
does not contain any crossing,
the disk $E$ does not contain any ring.
\end{proof}

Let $\Gamma$ be a chart $\Gamma$, and 
$m$ a label of $\Gamma$.
For a 3-color disk $E$ 
with $\partial E\subset\Gamma_m$, 
we define

$\mathcal{W}^{{\rm Mid}}(E,m\pm 1) = 
\{ v\in \mathcal{W}(E) \ |
\text{ there is an arc of label $m\pm 1$ 
in $E$ middle at $v$}\}.$


\begin{lemma}
\label{DecompositionM}
Let $\Gamma$ be a minimal chart, and 
$m$ a label of $\Gamma$. 
Let $E$ be a $3$-color disk 
bounded by a cycle $C$ in $\Gamma_m$.
Then we have\\
$|\mathcal{W}^{{\rm Mid}}(E,m\pm 1)|
=|\mathcal{W}({\rm Int}~E)|
+|\mathcal{W}_O^{{\rm Mid}}(C,m)|
+|\mathcal{W}_I(C,m)|
-|\mathcal{W}_I^{{\rm Mid}}(C,m)|.$
\end{lemma}

\begin{proof}
By the definition of $\mathcal{W}^{{\rm Mid}}(E,m\pm 1)$,
we have 
$\mathcal{W}^{{\rm Mid}}(E,m\pm 1) 
\subset 
\mathcal{W}(E)$.
By Remark~\ref{O(C)I(C)}$(4)$ and $(5)$,
the set $\mathcal{W}(E)$ is 
the disjoint union of 
$\mathcal{W}({\rm Int}~E)$ and 
$\mathcal{W}(C)$.
Thus

$\mathcal{W}^{{\rm Mid}}(E,m\pm 1)
-\mathcal{W}({\rm Int}~E)
=\mathcal{W}^{{\rm Mid}}(E,m\pm 1)
\cap\mathcal{W}(C)$.\\
Since $\mathcal{W}(C)$ 
is the disjoint union of 
$\mathcal{W}_O(C,m)$ and 
$\mathcal{W}_I(C,m)$
by Remark~\ref{O(C)I(C)}$(4)$, we have\\
$\begin{array}{l}
|\mathcal{W}^{{\rm Mid}}(E,m\pm 1)-\mathcal{W}({\rm Int}~E)|\\
=
|\mathcal{W}^{{\rm Mid}}(E,m\pm 1)\cap\mathcal{W}(C)|\\
=
|\mathcal{W}^{{\rm Mid}}(E,m\pm 1)\cap 
(\mathcal{W}_O(C,m)\cup
 \mathcal{W}_I(C,m))|\\
=
|(\mathcal{W}^{{\rm Mid}}(E,m\pm 1)\cap
  \mathcal{W}_O(C,m))
\cup 
 (\mathcal{W}^{{\rm Mid}}(E,m\pm 1)\cap
  \mathcal{W}_I(C,m))|\\
=
|\mathcal{W}^{{\rm Mid}}(E,m\pm 1)\cap
 \mathcal{W}_O(C,m)|
+ 
|\mathcal{W}^{{\rm Mid}}(E,m\pm 1)\cap
 \mathcal{W}_I(C,m)|.
\end{array}
$\\
Now by Remark~\ref{O(C)I(C)}$(1)$, we have\\
$\begin{array}{l}
\mathcal{W}^{{\rm Mid}}(E,m\pm 1)\cap 
\mathcal{W}_O(C,m)\\
=\{
v\in \mathcal{W}_O(C,m) ~|~
\text{there exists an arc of $m\pm 1$ in $E$ 
middle at $v\}$}\\
=\mathcal{W}_O^{{\rm Mid}}(C,m).
\end{array}
$\\
Using Remark~\ref{O(C)I(C)}(2), we have\\
$\begin{array}{l}
\mathcal{W}^{{\rm Mid}}(E,m\pm 1)\cap 
\mathcal{W}_I(C,m)\\
=\{
v\in \mathcal{W}_I(C,m) ~|~
\text{there exists an arc of $m\pm 1$ in $E$ 
middle at }v\}\\
=\{ v\in \mathcal{W}_I(C,m) ~|~
v\not\in \mathcal{W}_I^{{\rm Mid}}(C,m)~\}\\
=\mathcal{W}_I(C,m)
-\mathcal{W}_I^{{\rm Mid}}(C,m).
\end{array}
$\\
Therefore\\
$\begin{array}{l}
|\mathcal{W}^{{\rm Mid}}(E,m\pm 1)|\\
=
|\mathcal{W}({\rm Int}~E)|
+|\mathcal{W}^{{\rm Mid}}(E,m\pm 1)
-\mathcal{W}({\rm Int}~E)|\\
=|\mathcal{W}({\rm Int}~E)|
+|\mathcal{W}^{{\rm Mid}}(E,m\pm 1)
\cap \mathcal{W}_O(C,m)|
+|\mathcal{W}^{{\rm Mid}}(E,m\pm 1)
\cap \mathcal{W}_I(C,m)|\\
=|\mathcal{W}({\rm Int}~E)|
+|\mathcal{W}_O^{{\rm Mid}}(C,m)|
+|\mathcal{W}_I(C,m)
-\mathcal{W}_I^{{\rm Mid}}(C,m)|
\end{array}
$\\
$\begin{array}{l}
=|\mathcal{W}({\rm Int}~E)|
+|\mathcal{W}_O^{{\rm Mid}}(C,m)|
+|\mathcal{W}_I(C,m)|
-|\mathcal{W}_I^{{\rm Mid}}(C,m)|.
\end{array}$
\end{proof}

{\it Proof of Theorem~\ref{Inequation}.}
Let $C=\partial E$.
Since the disk $E$ is a 3-color disk,

$(1)$ $E$ does not contain any crossing by Condition (i) of a $3$-color disk. \\
Since $\Gamma_m\cap E$ is connected
by the assumption,

$(2)$ $E$ does not contain a hoop nor a ring
by Lemma~\ref{3ColorDisk}.\\
Thus

\begin{enumerate}
\item[$(3)$] 
for each connected component of $E-\Gamma_m$,
its closure contains a white vertex.
\end{enumerate}
Since $E$ is a $3$-color disk with $\partial E\subset\Gamma_m$,
by Corollary~\ref{Cor3ColorDiskNoTerminalEdge}(a)
the disk $E$ does not contain any terminal edge of 
label $m$.
Hence all the vertices of $\Gamma_m$ in $E$ are
white vertices.
Let 
$\mathcal V$ be the number of 
white vertices of $\Gamma_m$ in $E$,
and
$\mathcal E$ the number of 
edges of $\Gamma_m$ in $E$,
and
$\mathcal F$ the number of connected components 
of $E-\Gamma_m$. 
Since $E$ is a 3-color disk 
 with $\partial E\subset\Gamma_m$,
all the white vertices in $E$ are 
white vertices of $\Gamma_m$.
Thus we have

$(4)$~$\mathcal V=|\mathcal{W}(E)|$.\\
By Remark~\ref{O(C)I(C)}$(5)$,
the set $\mathcal{W}(E)$ splits into 
three mutually disjoint subsets 
$\mathcal{W}({\rm Int}~E)$,
$\mathcal{W}_O(C,m)$ and 
$\mathcal{W}_I(C,m)$.
Thus by $(4)$, we have

(5) $\mathcal V=|\mathcal{W}({\rm Int}~E)|
+|\mathcal{W}_O(C,m)|
+|\mathcal{W}_I(C,m)|$.\\
Since the intersection 
$\Gamma_m\cap E$ is connected
by the assumption,
 
$(6)$ each connected component of $E-\Gamma_m$
is an open disk.

{\bf Claim 1.} 
$2\mathcal F
=2+\mathcal V-|\mathcal{W}_O(C,m)|$.

{\it Proof of Claim $1$.}
In a small neighbourhood of 
a white vertex in $\Gamma_m$,
there are exactly three short arcs of label $m$ 
intersecting each other at the white vertex.
We fix the three short arcs of label $m$ 
for each white vertex in $\Gamma_m$.
Each edge of $\Gamma_m$ in $E$ has 
two short arcs of label $m$.
Thus there are $2\mathcal E$ short arcs in $E$.
On the other hand, 
any white vertex in 
$\mathcal{W}({\rm Int}~E)
\cup \mathcal{W}_I(C,m)$ is 
incident with three short arcs 
of label $m$ in $E$,
and 
any white vertex in 
$\mathcal{W}_O(C,m)$ is 
incident with two short arcs 
of label $m$ in $E$.
Hence we have
$$2\mathcal E =  3(|\mathcal{W}({\rm Int}~E)|+
|\mathcal{W}_I(C,m)|)
+2|\mathcal{W}_O(C,m)|.$$
Thus by the equation (5),
we have 
$$
\begin{array}{ll}
2\mathcal E  &=  
3(|\mathcal{W}({\rm Int}~E)|
+|\mathcal{W}_I(C,m)|
+|\mathcal{W}_O(C,m)|)
-|\mathcal{W}_O(C,m)|\\
&= 3\mathcal V-|\mathcal{W}_O(C,m)|.
\end{array}$$
Since $E$ is a disk, 
we have the equation $\mathcal V-\mathcal E+ \mathcal F=1$ 
by Euler formula.
Hence
$$
\begin{array}{ll}
2\mathcal F
&=2-2\mathcal V+2\mathcal E\\
&=2-2\mathcal V+(3\mathcal V-|\mathcal{W}_O(C,m)|)\\
&=2+\mathcal V-|\mathcal{W}_O(C,m)|.
\end{array}$$ 
This completes the proof of Claim~$1$.

{\bf Claim 2.}
 $|\mathcal{W}^{{\rm Mid}}(E,m\pm 1)|
 -2\mathcal F\ge0$.

{\it Proof of Claim $2$.}  
Let ${\mathcal{M}^*}(E,m\pm 1)$ 
be the set of all the middle arcs 
of label $m\pm1$ in $E$. Now

(7) each white vertex is contained in 
exactly one middle arc of label $m\pm 1$.\\
Hence we have

(8) $|\mathcal{M}^*(E,m\pm 1)|
=|\mathcal{W}^{{\rm Mid}}(E,m\pm 1)|.$\\
By $(1)$, $(3)$ and $(6)$,
Lemma~\ref{NumberMiddleEdge}(c)
assures that
for each connected component $U$ 
of $E-\Gamma_m$,
the closure $Cl(U)$ contains
at least two middle arcs of 
label $m\pm1$ in $\mathcal{M}^*(E,m\pm 1)$.
Hence (7) assures us

$|\mathcal{M}^*(E,m\pm 1)|-2\mathcal F\ge0.$\\
Thus (8) implies 
$|\mathcal{W}^{{\rm Mid}}(E,m\pm 1)|
-2\mathcal F\ge0.$
This completes the proof of Claim~$2$.

By Lemma~\ref{DecompositionM} and 
$(5)$, 
we have\\
$\begin{array}{l}
|{\mathcal{W}^{{\rm Mid}}(E,m\pm 1)}|\\
= |\mathcal{W}({\rm Int}~E)|
+|\mathcal{W}_O^{{\rm Mid}}(C,m)|
+|\mathcal{W}_I(C,m)|
-|\mathcal{W}_I^{{\rm Mid}}(C,m)|\\
= 
(\mathcal V-|\mathcal{W}_O(C,m)|
-|\mathcal{W}_I(C,m)|)
+|\mathcal{W}_O^{{\rm Mid}}(C,m)|
+|\mathcal{W}_I(C,m)|
-|\mathcal{W}_I^{{\rm Mid}}(C,m)|\\
= 
\mathcal V
-|\mathcal{W}_O(C,m)|
+|\mathcal{W}_O^{{\rm Mid}}(C,m)|
-|\mathcal{W}_I^{{\rm Mid}}(C,m)|.
\end{array}$\\
By Claim~$1$ and Claim~$2$, we have\\
$\begin{array}{l}
0\le |\mathcal{W}^{{\rm Mid}}(E,m\pm 1)|-2\mathcal F\\ 
= (\mathcal V-|\mathcal{W}_O(C,m)|
+|\mathcal{W}_O^{{\rm Mid}}(C,m)|
-|\mathcal{W}_I^{{\rm Mid}}(C,m)|)
- (2+\mathcal V-|\mathcal{W}_O(C,m)|)\\
=|\mathcal{W}_O^{{\rm Mid}}(C,m)|
-|\mathcal{W}_I^{{\rm Mid}}(C,m)|-2.
\end{array}$\\
Therefore
$|\mathcal{W}_I^{{\rm Mid}}(C,m)|+2
\le |\mathcal{W}_O^{{\rm Mid}}(C,m)|.$
This proves Theorem~\ref{Inequation}.
{\hfill {$\square$}\vspace{1.5em}}

\begin{cor}
\label{CorInequation}
Let $\Gamma$ be a minimal chart, and 
$m$ a label of $\Gamma$.
Let $E$ be a $3$-color disk with 
$\partial E\subset \Gamma_m$
but without free edges nor simple hoops. 
If $\Gamma_m\cap E$ is connected, 
then we have 
$$ |\mathcal{W}_O^{{\rm Mid}}(\partial E,m)|\ge 2.$$
\end{cor}


\section{Dichromatic one-side pseudo paths and $2$-color disks}
\label{s:Dichromatic}

In this section,
we investigate an arc 
in the boundary of a $2$-color disk $E$  with $\partial E\subset\Gamma_m$
for a minimal chart $\Gamma$. 
We shall show that
if $\Gamma_m\cap E$ is connected and
if $E$ contains neither free edge nor simple hoop,
then $|\mathcal{W}_O(\partial E,m)-\mathcal{W}_O^{\rm Mid}(\partial E,m)|\ge2$.

Let $\Gamma$ be a chart, and 
$m$ a label of $\Gamma$.
A simple arc $P$ in $\Gamma_m$
is called a {\it path} of label $m$
provided that the endpoints of $P$ are vertices of $\Gamma$.
A simple arc $P^*$ in $\Gamma_m$
is called a {\it pseudo path} of label $m$
provided that
\begin{enumerate}
\item[(i)] $P^*$ contains at least one vertex of $\Gamma$, and
\item[(ii)]
the endpoints of $P^*$
are not black vertices, crossings, nor white vertices.
\end{enumerate}

\begin{figure}[h]
\begin{center}
\includegraphics{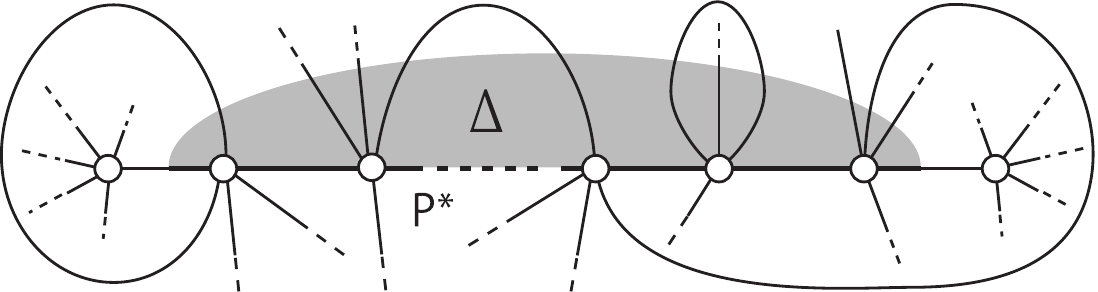}
\end{center}
\caption{ \label{fig14}
The thick line is a pseudo path $P^*$, and 
the gray area is a side-disk $\Delta$.}
\end{figure}

Let $P^*$ be a pseudo path of label $m$ in a chart $\Gamma$.
A disk $\Delta$ is called a 
{\it side-disk} of $P^*$
provided that $P^*\subset \partial \Delta$ 
(see Fig.~\ref{fig14}). 
Let $N$ be a regular neighborhood of $P^*$ 
in the side-disk $\Delta$. 
Let $e$ be an edge of $\Gamma$, and 
$\gamma$ the closure of a connected component of 
$e\cap $Int~$N$. 
If $\gamma$ contains a white vertex in $P^*$,
then $\gamma$ is called a {\it side-arc} of $P^*$
with respect to $\Delta$.
A side-arc is said to be
{\it at a vertex $v$}
if it contains the vertex $v$.
Similarly,
a side-arc is said to be
{\it middle} at a vertex $v$
if it contains a middle arc at $v$.

Let $P^*$ 
be a pseudo path of label $m$ in a chart
with a side-disk $\Delta$.
The pseudo path $P^*$ is said to be 
{\it inward $($resp. outward$)$} 
with respect to $\Delta$
provided that
\begin{enumerate}
\item[(i)] all the vertices in $P^*$ are white vertices, and
\item[(ii)]
for each vertex $v$ in $P^*$, 
any side-arc at $v$ with respect to $\Delta$ 
is oriented inward (resp. outward)
at $v$.
\end{enumerate}
An inward pseudo path and an outward pseudo path 
are called {\it I/O pseudo paths}.

Let $P^*$ be a pseudo path of label $m$ in a chart, and 
$v_1,v_2,\cdots,v_{s}$ all the vertices in $P^*$ 
which are situated in this order on $P^*$,
here some of $v_1,v_2,\cdots,v_{s}$ may be crossings.
For each $i=1,2,\cdots,s-1$, 
let $e_i$ be the edge of label $m$ in $P^*$ with $\partial e_i=\{v_i,v_{i+1}\}$.
Then the $s$-tuple $(v_1,v_2,\cdots,v_{s})$
is called 
the {\it associated vertex sequence} 
for the pseudo path $P^*$,
and the $(s-1)$-tuple $(e_1,e_2,\cdots,e_{s-1})$
is called 
the {\it associated edge sequence} 
for the pseudo path $P^*$.
The path $e_1\cup e_2\cup\cdots\cup e_{s-1}$
is denoted by $L(P^*)$.
The path $L(P^*)$ is the maximal path 
contained in the pseudo path $P^*$. 
Let $\gamma_0$ and $\gamma_s$ be arcs in edges of label $m$ with $\gamma_0\ni v_1,\gamma_s\ni v_s$ and
$P^*=\gamma_0\cup L(P^*)\cup \gamma_s$.
Then $\gamma_0,\gamma_s$ are called 
the {\it end-arcs} of the pseudo path $P^*$.

Let $P^*$ be a pseudo path of label $m$ in a chart, and
$(v_1,v_2,\cdots,v_{s})$ the associated vertex sequence for $P^*$.
The pseudo path $P^*$ 
is said to be {\it admissible} provided that 
\begin{enumerate}
\item[(i)] the vertices $v_1,v_s$ are white vertices, 
possibly $v_1=v_s$, 
\item[(ii)] 
there exists a side-disk $\Delta$ 
so that 
for each $i=1,s$,
there does not exist a 
side-arc of label $m$ at $v_i$
with respect to the side-disk $\Delta$.
\end{enumerate}
We also say that 
$P^*$ is admissible for the side-disk $\Delta$.\\

{\bf Important Notice.} 
Let $P^*$ be a pseudo path of label $m$ in a chart,
and $D$ a disk with $P^*\cap\partial D=\partial P^*$. Then $P^*$ splits $D$ into two disks $\Delta_1,\Delta_2$.
Both of the two disks are side-disks of $P^*$.
{\it 
If $P^*$ is an admissible pseudo path, 
then $P^*$ is admissible for 
one of the two side-disks $\Delta_1,\Delta_2$,
but NOT admissible for 
the other side-disk.
Thus for admissible pseudo paths,
we do not mention side-disks
unless the side-disks are 
needed to be mentioned.
Thus for an admissible pseudo path $P^*$,
'a side-arc' means a side-arc
with respect to a side-disk 
for which $P^*$ is admissible.}\\

A pseudo path $P^*$ of label $m$ in a chart 
is called a {\it dichromatic pseudo path} 
if there exists a label $k$
with $|m-k|=1$ such that
any vertex in $P^*$ is contained in an edge of label $k$.
The label $k$ is called the
{\it secondary label} 
of the dichromatic pseudo path $P^*$.

\begin{figure}
\begin{center}
\includegraphics{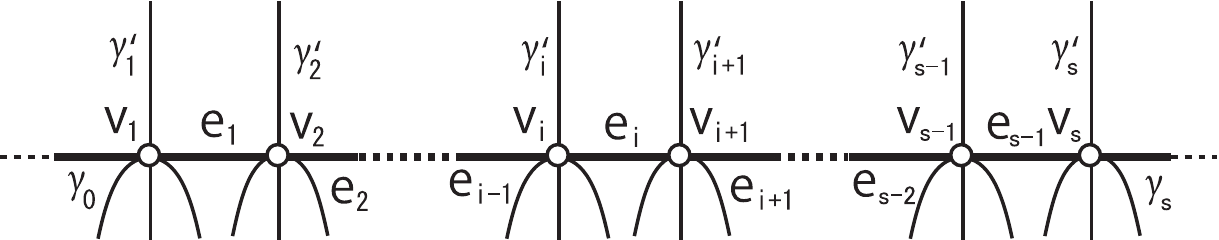}
\end{center}
\caption{ \label{fig15} The thick line is a one-side pseudo path $P^*$. 
}
\end{figure}

A pseudo path $P^*$ of label $m$ in a chart is called a {\it one-side pseudo path}
if there exists a side-disk $\Delta$ of $P^*$ such that 
(see Fig.~\ref{fig15})
\begin{enumerate}
\item[(i)] all the vertices in $P^*$ are white vertices, 
and
\item[(ii)]
for each vertex in $P^*$ 
the label of any side-arc with respect to $\Delta$
is different from $m$.
\end{enumerate}

\begin{remark}
\label{NoteOneSidePath}
{\rm
\begin{enumerate}
\item[$(1)$] 
Any one-side pseudo path is admissible.
Thus {\it we do not mention side-disks 
for 
one-side pseudo paths.}
\item[$(2)$]
If $P^*$ 
is a one-side (resp. dichromatic) pseudo path, 
then any pseudo path in $P^*$ 
is a one-side (resp. dichromatic) pseudo path.
\end{enumerate}
}
\end{remark}

Let $P^*$ be a one-side pseudo path of label $m$ in a chart, and $(v_1,v_2,\cdots,v_s)$ the associated vertex sequence for $P^*$.
Then for each vertex 
$v_i$ $(i=1,2,\cdots,s)$,
there exists a side-arc 
$\gamma_i'$ of label $m\pm1$ 
at $v_i$.
The $s$-tuple 
$(\gamma_1',\gamma_2',\cdots,\gamma_{s}')$ 
is called an 
{\it associated side-arc sequence} 
for the one-side pseudo path $P^*$.


\begin{lemma}
\label{OneWayPathLemma}
Let $\Gamma$ be a minimal chart.
Let $P^*$ be 
a dichromatic one-side pseudo path 
of label $m$ in $\Gamma$
with the associated vertex sequence $(v_1,v_2,\cdots,v_{s})$,
the associated edge sequence $(e_1,e_2,\cdots,e_{s-1})$ and
an associated side-arc sequence 
$(\gamma_1',\gamma_2',\cdots,\gamma_{s}')$.
Let $\gamma_0$ and $\gamma_s$ be the end-arcs of $P^*$ with
$\gamma_0\ni v_1$ and $\gamma_s\ni v_s$.
Suppose that each side-arc $\gamma_i'$
$(i=1,2,\cdots,s)$ 
is not middle at $v_i$. Then we have the following:
\begin{enumerate}
\item[{\rm $($a$)$}] 
If the end-arc $\gamma_0$ is oriented inward at $v_1$, then 
 each edge $e_i$ 
$(i=1,2,\cdots,s-1)$ is oriented inward at 
$v_{i+1}$, and
the end-arc $\gamma_s$ is oriented outward at $v_{s}$.
\item[{\rm $($b$)$}] 
If the end-arc $\gamma_0$ is oriented outward at $v_1$, 
then each edge $e_i$ 
$(i=1,2,\cdots,s-1)$ is oriented outward at 
$v_{i+1}$, and
the end-arc $\gamma_s$ is oriented inward at $v_{s}$. 
\item[{\rm $($c$)$}] The pseudo path $P^*$ is 
an I/O pseudo path.
\end{enumerate}
\end{lemma}

\begin{figure}
\begin{center}
\includegraphics{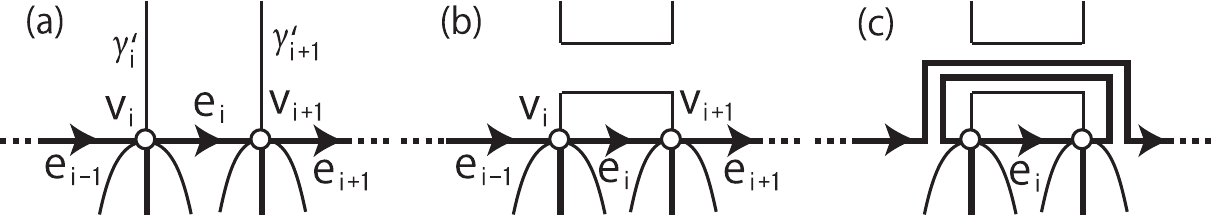}
\end{center}
\caption{ \label{fig16} 
}
\end{figure}

\begin{proof}
Set $e_0=\gamma_0,e_s=\gamma_s$.
Let $v_{s+1}$ be the endpoint of $\gamma_s$ different from $v_{s}$.

We show Statement (a). 
Now the end-arc $\gamma_0$ is oriented 
inward at $v_1$.
Suppose that  
 for some integer $i$ $(1\le i\le s)$
 $e_i$ is oriented 
outward at $v_{i+1}$.
Let $t=\min \{~j~|~e_j$ is oriented 
outward at $v_{j+1}~\}$.
Then $t\ge 1$ and 
 $e_t$ is oriented 
outward at $v_{t+1}$.
Thus  $e_{t}$ is oriented 
inward at $v_{t}$.
But  $e_{t-1}$ is oriented 
inward at $v_{t}$.
Hence the side-arc $\gamma_{t}'$ is 
middle at $v_{t}$.
This contradicts the assumption.
Hence Statement (a) holds.

Similarly
we can  show Statement (b).

We show Statement (c). 
We only show the case that
the end-arc $\gamma_0$ is oriented inward at $v_1$.
Then by Statement (a),
each edge $e_i$ 
$(i=1,2,\cdots,s-1)$ is oriented inward at 
$v_{i+1}$, and
the end-arc $\gamma_s$ is oriented outward at $v_{s}$.
Suppose that 
the pseudo path $P^*$ is not an I/O pseudo path.
Then for some two integers $i,j$ 
with $1\le i<j \le s$,
one of the following occurs.
\begin{enumerate}
\item[$(1)$]
The side-arc $\gamma_i'$ is oriented inward 
at $v_i$, and
the side-arc $\gamma_j'$ is oriented outward 
at $v_j$.
\item[$(2)$]
The side-arc $\gamma_i'$ is oriented outward 
at $v_i$, and 
the side-arc $\gamma_j'$ is oriented inward 
at $v_j$.
\end{enumerate}
Without loss of generality we can assume
$j=i+1$.
We show that we can eliminate 
the two white vertices 
$v_i$ and $v_{i+1}$ by C-moves.
For, either Case $(1)$ or Case $(2)$, 
we can apply a C-I-M2 move between 
the two side-arcs 
$\gamma_{i}'$ and $\gamma_{i+1}'$ (see Fig.~\ref{fig16}(a) and (b)).
Since $e_{i-1}$ is oriented 
inward at $v_{i}$ and 
since $e_{i+1}$ is oriented 
inward at $v_{i+2}$,
we can apply a C-I-M2 move between 
 $e_{i-1}$ and $e_{i+1}$ (see Fig.~\ref{fig16}(c)).
Finally
by applying a C-I-M3 move,
we can eliminate the two white vertices 
$v_{i}$ and $v_{i+1}$.
This contradicts the fact that 
the chart $\Gamma$ is minimal.
Hence the pseudo path $P^*$ 
is an I/O pseudo path.
Hence Statement (c) holds.
\end{proof}

\begin{lemma}
\label{OneSideLemma0}
Let $\Gamma$ be a minimal chart.
Let $P^*$ be 
a dichromatic one-side pseudo path 
of label $m$ in $\Gamma$
with the associated vertex sequence 
$(v_1,v_2,\cdots,v_{s})$ and
an associated side-arc sequence 
$(\gamma_1',\gamma_2',\cdots,\gamma_{s}')$.
Let $\gamma_0$ and $\gamma_s$ be 
the end-arcs of $P^*$ with
$\gamma_0\ni v_1$ and $\gamma_s\ni v_s$.
Suppose that 
\begin{enumerate}
\item[{\rm $($a$)$}] 
the end-arc $\gamma_0$ is middle at $v_1$ or 
the end-arc $\gamma_s$ is middle at $v_{s}$, and
\item[{\rm $($b$)$}] 
for each $i=2,3,\cdots,s-1$,
the side-arc $\gamma_i'$ is not middle at $v_i$.
\end{enumerate}
Then the pseudo path $P^*$ is 
an I/O pseudo path.
\end{lemma}

\begin{proof}
By (a),
without loss of generality
we can assume that 
\begin{enumerate}
\item[$(1)$] the end-arc $\gamma_s$ 
is middle at $v_s$.
\end{enumerate}
The side-arc $\gamma_1'$ is middle at $v_1$
or
not middle at $v_1$.

First, suppose that
the side-arc $\gamma_1'$ is not middle at $v_1$.
Then by (b) and $(1)$, 
we have that
for each $i=1,2,\cdots,s$,
the side-arc $\gamma_i'$ is not middle at $v_i$.
Thus by Lemma~\ref{OneWayPathLemma}(c),
the pseudo path $P^*$ is an I/O pseudo path.

Next, suppose that 
the side-arc $\gamma_1'$ is middle at $v_1$.
If $s=1$ then our lemma is true.
Thus we assume $s\ge 2$.
Let $(e_1,e_2,\cdots,e_{s-1})$ be the associated edge sequence for $P^*$.
Without loss of generality we can assume that
\begin{enumerate}
\item[$(2)$] the side-arc $\gamma_1'$ 
is oriented inward at $v_1$.
\end{enumerate}
Then the edge $e_1$ 
is oriented inward at $v_1$.
Thus the edge $e_1$ 
is oriented outward at $v_2$.
Since the end-arc $\gamma_s$ is middle at $v_{s}$
by $(1)$,
the side-arc $\gamma_{s}'$ 
is not middle at $v_{s}$.
Hence by (b), 
we have that
for each $i=2,3,\cdots,s$,
the side-arc $\gamma_i'$ is not middle at $v_i$.
Let $\gamma$ be a short arc containing $v_2$ in the edge $e_1$.
We can apply Lemma~\ref{OneWayPathLemma} 
to the dichromatic one-side pseudo path 
$P_1^*=\gamma\cup e_2\cup\cdots\cup e_{s-1}\cup\gamma_s$.
Then
the pseudo path $P^*_1$ is an I/O pseudo path,
and
the end-arc $\gamma_s$ is 
oriented inward at $v_{s}$
because $e_1$ is oriented outward at $v_2$
(i.e. the arc $\gamma$ is oriented outward at $v_2$).
Thus 
$\gamma_{s}'$ is oriented inward at $v_{s}$
because the end-arc $\gamma_s$ is middle at $v_{s}$
by $(1)$.
Hence 
the pseudo path $P^*_1$ is an inward pseudo path.
Thus for each $i=2,3,\cdots,s$,
the side-arc $\gamma_{i}'$ is 
oriented inward at $v_{i}$.
Considering $(2)$, we have 
for each $i=1,2,\cdots,s$,
the side-arc $\gamma_{i}'$ is 
oriented inward at $v_{i}$.
Therefore 
the pseudo path $P^*$ is an inward pseudo path.
This completes the proof of 
Lemma~\ref{OneSideLemma0}.
\end{proof}


\begin{lemma}
\label{OneSideLemma}
Let $\Gamma$ be a minimal chart.
Let 
$P^*$ 
be a dichromatic one-side pseudo path 
of label $m$ in $\Gamma$ with
the associated vertex sequence $(v_1,v_2,\cdots,v_{s})$ and
an associated side-arc sequence 
$(\gamma_1',\gamma_2',\cdots,\gamma_{s}')$. Suppose that 
the end-arcs $\gamma_0$ and $\gamma_{s}$ of $P^*$ are 
middle at $v_1$ and $v_{s}$ 
respectively.
Then we have the following:
\begin{enumerate}
\item[{\rm $($a$)$}] $s\ge 3$.
\item[{\rm $($b$)$}] 
For some integer $t$ with
$2\le t\le s-1$,
the side-arc $\gamma_t'$ is middle at $v_t$.
\end{enumerate}
\end{lemma}

\begin{proof}
First, we shall prove Statement (b).
Without loss of generality
we can assume that 
the end-arc $\gamma_0$ is 
oriented inward at $v_1$.
Since $\gamma_0$ is middle at $v_1$, we have
\begin{enumerate}
\item[$(1)$]
the side-arc $\gamma_1'$ is oriented inward 
at $v_1$.
\end{enumerate}
Since the end-arcs $\gamma_0,\gamma_s$ are 
middle at $v_1,v_s$ respectively,
Remark~\ref{EdgesAroundVertex}$(1)$ implies
\begin{enumerate}
\item[$(2)$]
the side-arc $\gamma_1'$ (resp. $\gamma_{s}'$) 
is not
middle at $v_1$ (resp. $v_{s}$).
\end{enumerate}
Now suppose that for each $i=2,\cdots,s-1$
the side-arc $\gamma_i'$ 
is not middle at $v_i$.
Then by $(2)$, we have that 
for $i=1,2,\cdots,s$
the side-arc $\gamma_i'$ 
is not middle at $v_i$.
Hence by Lemma~\ref{OneWayPathLemma}(c),
the dichromatic one-side pseudo path $P^*$ is 
an I/O pseudo path.
Further since the end-arc $\gamma_0$ is 
oriented inward at $v_1$, 
by Lemma~\ref{OneWayPathLemma}(a)
the end-arc $\gamma_s$ is oriented outward 
at $v_{s}$. 
Since $\gamma_s$ is middle at $v_{s}$, 
we have
\begin{enumerate}
\item[$(3)$]
the side-arc $\gamma_{s}'$ is oriented
outward at $v_{s}$.
\end{enumerate}
Hence the side-arc $\gamma_1'$ is
oriented outward at $v_1$
because $P^*$ is an I/O pseudo path.
This contradicts $(1)$.
Hence for some integer $t$ with 
$2\le t\le s-1$, 
the side-arc $\gamma_t'$ is middle at $v_t$. 
Thus Statement (b) holds.

We show Statement (a).
By Lemma~\ref{OneSideLemma}(b), 
we have $s-1\ge2$.
Thus $s\ge3$.
\end{proof}

\begin{figure}
\begin{center}
\includegraphics{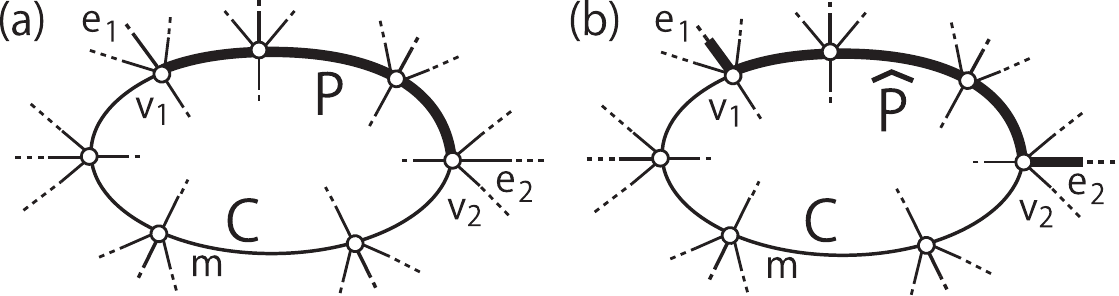}
\end{center}
\caption{ \label{fig17}
The vertices $v_1,v_2$ are 
in $\mathcal{W}_O(C,m)$.
(a) The thick line is a path $P$. 
(b) The thick line is 
the extended pseudo path $\widehat{P}$.}
\end{figure}

Let $C$ be a cycle of label $m$ 
in a chart $\Gamma$, and 
$P$ a path in $C$ with 
$\partial P=\{v_1,v_2\}\subset\mathcal{W}_O(C,m)$
(see Fig.~\ref{fig17}(a)).
For $i=1,2$,
let $e_i$ be the outside edge of label $m$ 
for $C$ containing $v_i$,
and $\gamma_i$ 
an arc in $e_i$
containing $v_i$ 
(see Fig.~\ref{fig17}(b)).
Set $\widehat{P}=\gamma_1\cup P\cup \gamma_2$.
Then the union $\widehat{P}$ is called 
an {\it extended pseudo path} of $P$.

\begin{remark}
\label{ExtendedPath}
{\rm Any extended pseudo path is admissible.}
\end{remark}

\begin{remark}
\label{TwoColorDisk}
{\rm 
By the condition of a $2$-color 
disk,
any $2$-color disk does not contain 
a crossing.
Hence 
any $2$-color disk is also a 
$3$-color disk.
}
\end{remark}

Let $\Gamma$ be a chart.
Let $C$ be a cycle or 
a path  of label $m$ 
in $\Gamma$,
and $\mathcal{S}$ a set of 
white vertices in $C$, here
we assume that
$|\mathcal{S}|\ge 2$ if $C$ is a cycle.
By cutting $C$ at all the white vertices 
in $\mathcal{S}$,
the set $C$ splits into paths.
Then the set of all the paths is called the 
{\it 
path decomposition of $C$ by $\mathcal{S}$}, 
denoted by $\mathcal{P}(C;\mathcal{S})$.

Let $m$ be a label of a chart $\Gamma$, and
$E$ a $2$-color disk with 
$\partial E\subset\Gamma_m$ 
and
$\Gamma\cap E\subset 
\Gamma_{m}\cup\Gamma_k$ for some label $k$. 
The label $k$ is called 
the {\it secondary label} of 
the $2$-color disk $E$.

\begin{lemma}
\label{DecompositionPaths}
Let $\Gamma$ be a chart and 
$m$ a label of $\Gamma$.
Let $C$ be a cycle of label $m$ 
in $\Gamma$ bounding a $2$-color disk.
Let $P$ be a path in 
$\mathcal{P}(C;\mathcal{W}_O(C,m))$.
Then any extended pseudo path of $P$
is a dichromatic one-side pseudo path.
\end{lemma}

\begin{proof}
Let $E$ be the $2$-color disk with $\partial E=C$, and 
$k$ the secondary label 
of the $2$-color disk $E$.
By Remark~\ref{TwoColorDisk},
the $2$-color disk does not contain 
any crossing, neither does $P$.
Let $v$ be a white vertex in $P-\partial P$. 
Then $v$ is not in $\mathcal{W}_O(C,m)$
by the assumption.
Hence the vertex $v$ is 
in $\mathcal{W}_I(C,m)$.
Thus there exists exactly one edge
 of label $m\pm 1$ containing the vertex $v$ in $Cl(E^c)$.
Since the cycle $C$ bounds a $2$-color
disk, 
the label of the edge is $k$.
Therefore the extended pseudo path of $P$ is 
a dichromatic one-side pseudo path.
\end{proof}


\begin{lemma}
\label{qPrickle}
Let $\Gamma$ be a minimal chart.
Let $C$ be a cycle of label $m$ in $\Gamma$ 
bounding a $2$-color disk $E$ 
without free edges nor simple hoops. 
If $\Gamma_m\cap E$ is connected, 
then 
\begin{enumerate}
\item[{\rm $($a$)$}]
in the path decomposition 
$\mathcal{P}(C;\mathcal{W}_O^{\rm Mid}(C,m))$
there exist at least two paths each of which
contains a white vertex in 
$\mathcal{W}_O(C,m)-\mathcal{W}_O^{\rm Mid}(C,m)$,
\item[{\rm $($b$)$}] 
$|\mathcal{W}_O(C,m)-\mathcal{W}_O^{\rm Mid}(C,m)|\ge2$, 
and 
\item[{\rm $($c$)$}]
in the path decomposition 
$\mathcal{P}(C;\mathcal{W}_O(C,m)
-\mathcal{W}_O^{\rm Mid}(C,m))$
there exist at least two paths 
each of which contains 
a white vertex in $\mathcal{W}_O^{\rm Mid}(C,m)$.
\end{enumerate}
\end{lemma}

\begin{figure}[h]
\begin{center}
\includegraphics{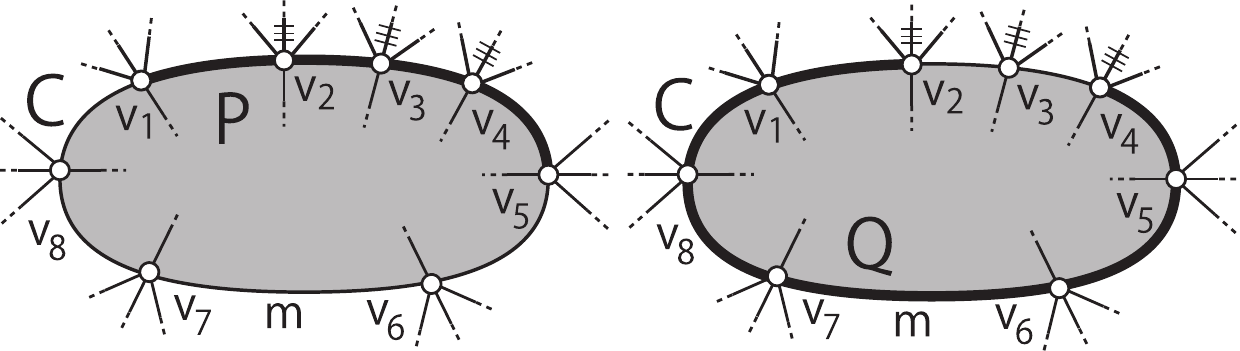}
\end{center}
\caption{ \label{fig18} 
An example of a path $P$ 
in a cycle $C$ of label $m$ with 
$\mathcal{W}_O^{\rm Mid}(C,m)\subset P$.
Each arc with three transversal short arcs 
is a middle arc. 
The thick lines are paths $P$ and $Q$ 
with $\partial P=\{v_1,v_5\}$ and 
$\partial Q=\{v_2,v_4\}$, 
$\mathcal{W}_O(C,m)=\{v_1,v_2,\cdots,v_8\}$ and
$\mathcal{W}_O^{\rm Mid}(C,m)=\{v_2,v_3,w_4\}$.
}
\end{figure}

\begin{proof}
By Remark~\ref{TwoColorDisk},
the $2$-color disk is also a 
$3$-color disk.
Let $s=|\mathcal{W}_O^{\rm Mid}(C,m)|$.
Then $s\ge2$ by Corollary~\ref{CorInequation}.
There are $s$ paths in the path decomposition 
$\mathcal{P}(C;\mathcal{W}_O^{\rm Mid}(C,m))$.

We show Statement (a).
Suppose that 
in $\mathcal{P}(C;\mathcal{W}_O^{\rm Mid}(C,m))$ 
there exists at most one path containing
a white vertex in 
$\mathcal{W}_O(C,m)-\mathcal{W}_O^{\rm Mid}(C,m)$.
Then 
in $\mathcal{P}(C;\mathcal{W}_O^{\rm Mid}(C,m))$ 
there are $s-1$ paths 
 each of which does not 
contain any vertices in $\mathcal{W}_O(C,m)$ 
except its end vertices. 
Hence the $s-1$ paths are in
$\mathcal{P}(C;\mathcal{W}_O(C,m))$.
By Lemma~\ref{DecompositionPaths},
the extended pseudo paths of the $s-1$ paths are 
dichromatic one-side pseudo paths.
Let $k$ be the secondary label of 
the $2$-color disk $E$.
By Lemma~\ref{OneSideLemma}(b), 
each of the pseudo paths has 
 a side-arc of label $k$ middle at 
 a vertex in the pseudo path,
namely each of the paths
contains a vertex in $\mathcal{W}_I^{\rm Mid}(C,m)$
by Remark~\ref{O(C)I(C)}$(3)$.
Thus $|\mathcal{W}_I^{\rm Mid}(C,m)|\ge s-1$.

By Remark~\ref{TwoColorDisk},
the $2$-color disk is a $3$-color disk.
Thus by Theorem~\ref{Inequation}
$$s + 1 = (s - 1) + 2 \le 
|\mathcal{W}_I^{\rm Mid}(C,m)| 
+ 2 \le |\mathcal{W}_O^{\rm Mid}(C,m)| = s.$$
This is a contradiction. 
Therefore Statement (a) holds.

Now Statement (b) is a direct result 
of Statement (a). 

We show Statement (c).
By Statement (b), 
we have 
$|\mathcal{W}_O(C,m)-\mathcal{W}_O^{\rm Mid}(C,m)|\ge2$.
Thus we can consider 
the path decomposition
$\mathcal{P}(C;\mathcal{W}_O(C,m)
-\mathcal{W}_O^{\rm Mid}(C,m))$.
Suppose that in the path decomposition
$\mathcal{P}(C;\mathcal{W}_O(C,m)
-\mathcal{W}_O^{\rm Mid}(C,m))$ 
there exists a path $P$
containing all the white vertices 
in $\mathcal{W}_O^{\rm Mid}(C,m)$ 
(see Fig.~\ref{fig18}).
Then
all the paths in 
$\mathcal{P}(C;\mathcal{W}_O^{\rm Mid}(C,m))$ 
are contained in $P$ except one, 
say $Q$. 
Thus $P\cup Q=C$.
Since $P\in \mathcal{P}(C;\mathcal{W}_O(C,m)
-\mathcal{W}_O^{\rm Mid}(C,m))$ implies
Int $P\cap (\mathcal{W}_O(C,m)-
\mathcal{W}_O^{\rm Mid}(C,m))=\emptyset$, 
we have 
$\mathcal{W}_O(C,m)-\mathcal{W}_O^{\rm Mid}(C,m)\subset Q$.
This contradicts Statement (a).
\end{proof}

\section{Suspicious cycles}
\label{s:SuspiciousCycle}

In this section,
for a tangle $(\Gamma\cap D,D)$ and 
each label $m$ with 
$\Gamma_m\cap D\neq\emptyset$,
we obtain an equation and 
search for conditions for the existence 
of a special cycle which 
never bounds a $2$-color disk 
in a minimal chart.
This special cycle is crucial in the proof of 
Theorem~\ref{NoNS-tangle}.

Let $G$ be a graph. For each vertex $v$ of $G$,
we denote by deg$_G v$ the degree of the vertex $v$
in $G$.

Let $\Gamma$ be a chart, 
and $m$ a label of $\Gamma$.
A tree $T$ in $\Gamma$ is called 
a {\it reducible tree} of label $m$
provided that
\begin{enumerate}
\item[(i)] 
each edge of $T$ is of label $m$,
\item[(ii)] 
the tree $T$ contains a white vertex,
\item[(iii)] 
if $v$ is a  crossing in $T$, 
then deg$_Tv=2$,
\item[(iv)]
if $T$ contains exactly one white vertex $v$, 
then deg$_Tv=3$,
\item[(v)] 
if $T$ contains at least two white vertices, 
then
\begin{enumerate}
\item[(a)]
for each white vertex $v$ of $T$, we have
deg$_T v=1$ or $3$, 
\item[(b)] 
there exists at most one 
white vertex $v_0$ in $T$
with deg$_T v_0=1$.
\end{enumerate}
\end{enumerate}
The white vertex $v_0$ with deg$_T v_0=1$
is called 
{\it the special vertex} 
of the reducible tree.

\begin{lemma}
\label{NoReducibleTree}
Any minimal chart
does not contain a reducible tree 
of any label.
\end{lemma}

\begin{proof}
Suppose that there exists a reducible tree $T$ 
of label $m$ in a minimal chart.
If $T$ contains exactly one white vertex $w$,
then $T$ contains three terminal edges
by Remark~\ref{Assumption0}$(1)$.
By Remark~\ref{EdgesAroundVertex}$(2)$,
two of the three terminal edges are 
not middle at $w$.
This contradicts Remark~\ref{Assumption0}$(2)$.
Hence $T$ contains at least two white vertices.

Let $T^*$ be the subtree 
obtained from the reducible tree $T$
by taking out all the terminal edges.
Since $T^*$ is a tree containing
at least two white vertices, 
there exist two white vertices 
each of whose degree in $T^*$ is $1$.
Let $w$ be one of them 
different from the special vertex of $T$.
By Condition (v) for a reducible tree, 
we have deg$_T w=3$. 
Hence in the reducible tree $T$
there exist two terminal edges 
of label $m$ at $w$.
By Remark~\ref{EdgesAroundVertex}$(2)$,
one of the two terminal edges is 
not middle at $w$.
This contradicts Remark~\ref{Assumption0}$(2)$.
\end{proof}

Let $C$ be a cycle of label $m$ 
in a chart $\Gamma$.
Then the cycle $C$ is said to 
be {\it suspicious}
provided that
\begin{enumerate}
\item[(i)]
the outside edges of label $m$ for $C$ 
are terminal edges except one,
\item[(ii)] 
the cycle contains a white vertex, and
\item[(iii)]
the cycle bounds a disk $E$ with
$\Gamma_m\cap E$ connected.
\end{enumerate}

The following lemma is an easy consequence of 
the definition of a suspicious cycle.

\begin{lemma}\label{No2ColorC1-tangle}
In a minimal chart $\Gamma$, 
for any label $m$ of $\Gamma$
there does not exist
a suspicious cycle of label $m$ 
bounding a $2$-color disk.
\end{lemma}

\begin{proof}
Suppose that there exists 
a suspicious cycle $C$
of label $m$ bounding 
a $2$-color disk $E$.
By Assumption~\ref{AssumptionFreeEdge} and 
Assumption~\ref{AssumptionFreeEdgeSimpleHoop},
the disk $E$ contains neither free edge nor simple hoop.
By Condition (i) of a suspicious cycle,
we have
$|\mathcal{W}_O(C,m)-\mathcal{W}_O^{{\rm Mid}}(C,m)|\leq 1$.
On the other hand, 
by 
Condition (iii) of a suspicious cycle,
Lemma~\ref{qPrickle}(b) implies that 
$|\mathcal{W}_O(C,m)-\mathcal{W}_O^{{\rm Mid}}(C,m)|\geq 2$.
This is a contradiction.
\end{proof}

Let $X$ be a subset of a chart.
A cycle in $X$ is said to be
{\it maximal in $X$}
if it is not contained in the disk bounded by 
another cycle in $X$.

Let $\Gamma$ be a chart, 
and $m$ a label of $\Gamma$.
Let $(\Gamma\cap D,D)$ be a tangle.
Let $G$ be a connected component 
of $\Gamma_m\cap D$ containing 
a white vertex, 
and
$E_1,E_2\cdots,E_d$ all the disks 
bounded by the maximal cycles in $G$. 
For each $i=1,2,\cdots,d$, let
\begin{enumerate}
\item[]
$\widehat E_i=E_i\cup(\cup\{e~|~e 
\text{ is a terminal edge in $G$ 
intersecting }E_i\})$.
\end{enumerate}
Let $P_1,P_2\cdots,P_{p}$ be 
all the closures of connected components of 
$G-\cup_{i=1}^d\widehat E_i$
not intersecting $\partial D$,
and 
$Q_1,Q_2,\cdots,Q_q$
all the closures of connected components of 
$G-\cup_{i=1}^d\widehat E_i$ intersecting 
$\partial D$.
Let
\begin{enumerate}
\item[]
$
\begin{array}{ll}
\mathcal{H}=(\cup_{i=1}^d \widehat E_i)
\cap ((\cup_{j=1}^p P_j)
\cup (\cup_{k=1}^q Q_k)),&
h=|\mathcal{H}|,\\
s=|(\cup_{i=1}^d \widehat E_i)
\cap (\cup_{j=1}^p P_j)|,&
t=|(\cup_{i=1}^d \widehat E_i)
\cap (\cup_{k=1}^q Q_k)|.
\end{array}
$
\end{enumerate}
For each $k=0,1,\cdots,h$, let

\begin{figure}[h]
\begin{center}
\includegraphics{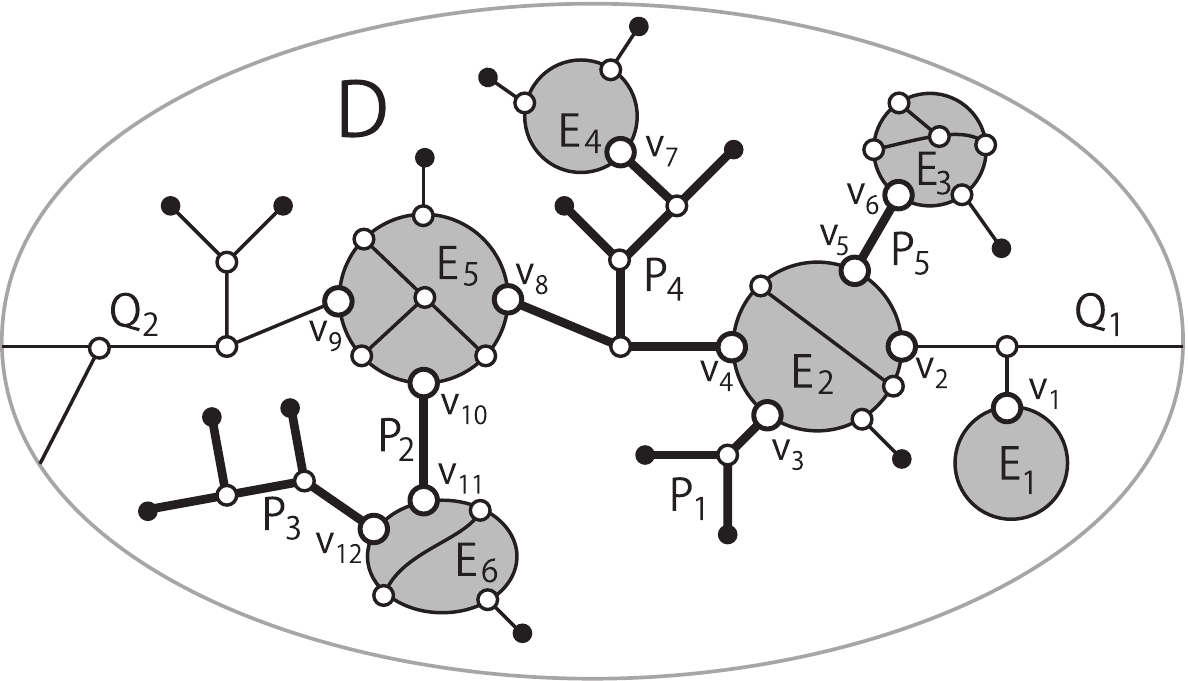}
\end{center}
\caption{\label{fig19}}
\end{figure}

$x_k=$ the number of $E_i$'s containing 
exactly $k$ points in $\mathcal{H}$,

$y_k=$ the number of $P_j$'s containing 
exactly $k$ points in $\mathcal{H}$.\\
Then 
$(E_1,E_2,\cdots,E_d;
P_1,P_2,\cdots,P_p;
Q_1,Q_2,\cdots,Q_q)$ and
$(\mathcal{H},h,s,t;
x_0,x_1,\cdots,x_h,$
$y_0,y_1,\cdots,y_h)$
are called the 
{\it primary fundamental information} 
and the
{\it secondary fundamental information} 
of $G$ respectively.
These will be used in 
Lemma~\ref{EQofMinimalC1-tangle}
and Lemma~\ref{MinimalC1-tangle}.
For $G$ in Fig.~\ref{fig19}, 
there are six disks $E_1,E_2,E_3,E_4,E_5,E_6$ 
and seven trees $P_1,P_2,P_3,P_4,P_5,Q_1,Q_2$. 
Thus we have
\begin{enumerate}
\item[]
$\mathcal{H}=\{v_1,v_2,v_3,v_4,v_5,v_6,v_7,v_8,
v_9,v_{10},,v_{11},v_{12}\},$
\item[]
$h=12,~~~s=9,~~~t=3,$
\item[]
$x_0=0,x_1=3,x_2=1,x_3=1, x_4=1,
x_5=\cdots=x_{12}=0,$
\item[]
$y_0=0,y_1=2,y_2=2,
y_3=1,y_4=\cdots=y_{12}=0$.
\end{enumerate}

Let $X$ be a subset of a chart.
A connected component $G$ of $X$ is called 
a {\it small component} of $X$ 
if any finite complementary domain of $G$
does not intersect $X$.


\begin{lemma}
\label{EQofMinimalC1-tangle}
Let $\Gamma$ be a minimal chart, 
and $m$ a label of $\Gamma$.
Let $(\Gamma\cap D,D)$ be a tangle, and 
$G$ a connected component 
of $\Gamma_m\cap D$ containing a white vertex with
the primary fundamental information
$(E_1,E_2,\cdots,E_d;
P_1,P_2,\cdots,P_p;
Q_1,Q_2,\cdots,Q_q)$ and 
the secondary fundamental information
$(\mathcal{H},h,s,t;
x_0,x_1,\cdots,x_h,
y_0,y_1,\cdots,y_h)$.
Then we have the following:
\begin{enumerate}
\item[{\rm $($a$)$}]
$y_0=y_1=0$.
\item[{\rm $($b$)$}]
$2x_0+x_1=
2-2q+t+
(x_3+y_3)+2(x_4+y_4)+\cdots+
(h-2)(x_{h}+y_{h})$.
\item[{\rm $($c$)$}]
Suppose that $G$ is a small component of $\Gamma_m\cap D$. 
If $1\le x_0$ or $1\le x_1$, 
then $D$ contains a suspicious cycle 
of label $m$.
\end{enumerate}
\end{lemma}

\begin{proof}
We show Statement (a).
Suppose $1\le y_0$ or $1\le y_1$.
Then
there exists a tree $P_j$ $(1\le j\le p)$ containing
at most one point of $\mathcal{H}$.
Hence $P_j$ is a reducible tree.
This contradicts Lemma~\ref{NoReducibleTree}.

We show Statement (b).
For each $i=1,2,\cdots,p$ 
and $j=1,2,\cdots,q$
we have $P_i\cap Q_j=\emptyset$.
Hence\\
$(1)~h=s+t,$\\
$(2)~x_0+x_1+\cdots+x_{h}=d,
y_0+y_1+\cdots+y_{h}=p,$ and \\
$(3)~1\times x_1+2\times x_2+\cdots+
h\times x_{h}=h,
~1\times y_1+2\times y_2+\cdots+
h\times y_{h}=s=h-t.$\\
By adding the two equations in $(3)$ we have\\
$(4)~(1\times x_1+2\times x_2+\cdots+
h\times x_{h})+
(1\times y_1+2\times y_2+\cdots+
h\times y_{h})=2\times h-t.
$\\
For each $i=1,2,\cdots,d$, let
\begin{enumerate}
\item[]
$\widehat E_i=E_i\cup(\cup\{e~|~e 
\text{ is a terminal edge in $G$ 
intersecting }E_i\})$.
\end{enumerate}
Let $E_*=\cup_{i=1}^d \widehat E_i,
P_*=\cup_{j=1}^p P_j,
Q_*=\cup_{k=1}^q Q_k$ and
$X=E_*\cup P_*\cup Q_*$.
Then we have $X=G\cup (\cup_{i=1}^d E_i)$.
Hence
the set $X$ is connected and 
any cycle in $X$ bounds a disk in $X$,
i.e. $X$ is simply connected.
Hence by Euler formula, we have
$\chi(X)=1$.
On the other hand, considering $(1)$ 
we have \\
$\chi(X)=\chi(E_*\cup P_*\cup Q_*)\\
=\chi(E_*)+
\chi(P_*)+
\chi(Q_*)
-(\chi(E_*\cap P_*)+\chi(E_*\cap Q_*)+\chi(P_*\cap Q_*))
+\chi(E_*\cap P_*\cap Q_*)\\
=d+p+q-(s+t+0)+0=d+p+q-h.$\\
Hence\\
$(5)~d+p+q-h=1$.\\
By doubling both sides of the equation $(5)$,
and
eliminating $d,p,h$ using $(2)$ and $(4)$, 
we have\\
$2(x_0+x_1+\cdots+x_{h})+
2(y_0+y_1+\cdots+y_{h})+2q\\
-((1\times x_1+2\times x_2+\cdots+
h\times x_{h})+
(1\times y_1+2\times y_2+\cdots+
h\times y_{h})+t)=2.$\\
Since $y_0=y_1=0$ by 
Statement (a), 
we have\\
$2x_0+x_1=
2-2q+t+
(x_3+y_3)+2(x_4+y_4)+\cdots+
(h-2)(x_{h}+y_{h})$.

We show Statement (c).
Suppose $1\le x_0$ or $1\le x_1$. 
Then
there exists an integer $i$ 
with $1\le i\le d$ 
such that 
the disk $E_i$ intersects 
$\mathcal{H}$ 
by at most one point.
Thus the outside edges of label $m$
for $\partial E_i$ 
are terminal edges except one.
Since $G$ is a small component 
of $\Gamma_m\cap D$,
the intersection $\Gamma_m\cap E_i$ is connected.
Since $G$ contains a white vertex, 
$G$ is not a hoop nor a ring.
Thus the cycle $\partial E_i$ contains 
a white vertex.
Therefore $\partial E_i$ is 
a suspicious cycle of label $m$.
\end{proof}

\begin{lemma}
\label{MinimalC1-tangle}
Let $\Gamma$ be a minimal chart, 
and $m$ a label of $\Gamma$.
Let $(\Gamma\cap D,D)$ be 
a tangle. 
Suppose that 
there exists a small component $G$ of 
$\Gamma_m\cap D$ 
containing a white vertex.
If  
$|G \cap \partial D|\le 1$, 
then 
$D$ contains 
a suspicious cycle of label $m$. 
\end{lemma}

\begin{proof}
First, we claim that
\begin{enumerate}
\item[$(1)$] 
$G$ contains a cycle.
\end{enumerate}
Suppose that 
$G$ does not contain a cycle.
If $|G\cap \partial D|=0$,
then $G$ is a reducible tree 
without the special vertex. 
This contradicts Lemma~\ref{NoReducibleTree}.
If $|G\cap \partial D|=1$,
let $\overline e$ be a terminal edge or 
an internal edge 
with $\overline e \ni G\cap\partial D$. 
Let $T=G\cup \overline e$.
Now for any white vertex $w$ of $G$,
we have deg$_T w=3$.
Hence $T$ is a reducible tree.
This contradicts Lemma~\ref{NoReducibleTree}.
Thus Statement $(1)$ holds. 

Let 
$(E_1,E_2,\cdots,E_d;
P_1,P_2,\cdots,P_p;
Q_1,Q_2,\cdots,Q_q)$ and
$(\mathcal{H},h,s,t;
x_0,x_1,\cdots,x_h,$
$y_0,y_1,\cdots,y_h)$
be the primary fundamental information  and 
the secondary fundamental information of 
the small component $G$.
By Lemma~\ref{EQofMinimalC1-tangle}(b), 
we have
\begin{enumerate}
\item[$(2)$] 
$2x_0+x_1=
2-2q+t+
(x_3+y_3)+2(x_4+y_4)+\cdots+
(h-2)(x_{h}+y_{h})$.
\end{enumerate}
Now $|G\cap\partial D|\le 1$
implies that $q=0$ or $1$.
If $q=0$ then $2-2q+t\ge 2$.
Hence $2x_0+x_1\ge 2$ by $(2)$.
Thus $x_0\ge 1$ or $x_1\ge 1$.
Hence there exists 
a suspicious cycle of label $m$
by Lemma~\ref{EQofMinimalC1-tangle}(c).

Suppose $q=1$.
By Statement $(1)$, we have $d\ge1$.
Then we have $t\ge 1$,
because $G\cup (\cup_{i=1}^d E_i)$ is connected.
Thus $2-2q+t\ge 1$.
Hence $2x_0+x_1\ge 1$ by Statement $(2)$.
Thus $x_0\ge 1$ or $x_1\ge 1$.
Therefore there exists 
a suspicious cycle of label $m$
by Lemma~\ref{EQofMinimalC1-tangle}(c).
This proves Lemma~\ref{MinimalC1-tangle}. 
\end{proof}


\section{Proof of Theorem~1.2}
\label{s:PfTheorem2}

In this section
we shall prove Theorem~\ref{NoNS-tangle}
by using Lemma~\ref{No2ColorC1-tangle}, Lemma~\ref{MinimalC1-tangle} 
and the three lemmata in this sections.

Let $\Gamma$ be a chart, and
$D$ a disk.
We define

$\omega(D)=$ 
the number of white vertices in $\Gamma\cap D$,

x$(D)=$ 
the number of crossings in $\Gamma\cap D$,

$n(\partial D)=$ 
the number of points in $\Gamma\cap\partial D$.

For a tangle $(\Gamma\cap D,D)$ 
in a chart $\Gamma$, 
let
$\tau(D)=(\omega(D),$ x$(D), n(\partial D))$.
We call $\tau(D)$ the {\it $\tau$-complexity} of the tangle.

Let $\Gamma$ be a chart. 
An NS-tangle $(\Gamma\cap D,D)$ 
of label $m$ 
is said to be 
{\it minimal} 
if its $\tau$-complexity of the tangle 
is minimal amongst the NS-tangles
of all the labels
with respect to 
the lexicographical order of
the triplet of integers.


\begin{lemma}
\label{MinimalNS-tangle}
Let $\Gamma$ be a minimal chart.
Let $(\Gamma\cap D,D)$ be 
a minimal NS-tangle of label $m$.
Then $D$ does not contain any ring.
\end{lemma}

\begin{proof}
Suppose that $D$ contains 
a ring $C$ of label $k$.
Then $C$ bounds a disk $E$.
Let $A$ be a regular neighbourhood of 
$\partial E$, and
$D'=Cl(E-A)$.
Then by Condition (iii) of an NS-tangle, 
the intersection $\Gamma_k\cap D$ 
contains at most one crossing.
Thus the ring $C=\partial E$ contains at most one crossing. Hence
\begin{enumerate}
\item[$(1)$] 
the intersection $\Gamma\cap\partial D'$ 
is at most one point.
\end{enumerate}
By Remark~\ref{Assumption0}$(3)$,
the disk $E$
contains a white vertex. Thus we have 

$(2)$ the disk $D'$ contains a white vertex.\\
By Condition (iii) of an NS-tangle,
for each label $i$ 
the intersection $\Gamma_i\cap D$ contains at most one crossing.
Since $D'\subset D$, we have

$(3)$ for each label $i$, 
the intersection $\Gamma_i\cap D'$ contains at most one crossing.\\
Hence $(\Gamma\cap D',D')$ is 
an NS-tangle.
Since $C$ is a ring with $C\cap D'=\emptyset$,
we have
$\omega(D')\le \omega(D)$, and
x$(D')<$ x$(D)$.
Hence we have $\tau(D')<\tau(D)$. 
This contradicts the fact that 
the NS-tangle $(\Gamma\cap D,D)$
is minimal.
\end{proof}

\begin{lemma}
\label{MinimalAdmissibleNS-tangle}
Let $\Gamma$ be a minimal chart.
If $(\Gamma\cap D,D)$ is 
a minimal NS-tangle of label $m$,
then $\partial D$ does not intersect any terminal edge.
\end{lemma}

\begin{proof}
Suppose that 
 $\partial D$ intersects a terminal edge $e$ 
with a black vertex $v$. 
Then there are two cases: $v\in D$ or 
$v\not\in D$.

{\bf Case 1.} $v\in D$.

Let $\ell$ be the connected component of $e\cap D$
with $v\in \ell$.
And let $N$ be a regular neighbourhood of $\ell$ in $D$,
and $D'=Cl(D-N)$.
Then $D'$ is a disk.
Since the terminal edge $e$ does not contain any crossings, we have $\Gamma\cap D'=(\Gamma\cap D)-\ell$
and $\Gamma\cap \partial D'=(\Gamma\cap \partial D)-\ell$.
Thus $(\Gamma\cap D',D')$ is an NS-tangle with
$\omega(D')=\omega(D)$, x$(D')=$x$(D)$, and
$n(\partial D')=n(\partial D)-1$.
Hence $\tau(D')<\tau(D)$.
This contradicts the fact that 
the NS-tangle $(\Gamma\cap D,D)$ is minimal.

{\bf Case 2.} $v\not\in D$.

Let $\ell^*$ be the arc in the terminal edge $e$
connecting $v$ and a point $p$ in $\partial D$
with $\ell^*\cap D=\{ p\}$. 
Let $N^*$ be a regular neighbourhood of $\ell^*$,
$D^*=D \cup N^*$. 
Then $D^*$ is a disk.
Since the terminal edge $e$ does not contain any crossings,
we have $\Gamma\cap D^*=(\Gamma\cap D)\cup \ell^*$
and $\Gamma\cap \partial D^*=(\Gamma\cap \partial D)-\{p\}$.
Thus 
$(\Gamma\cap D^*, D^*)$ is an NS-tangle with
$\omega(D^*)=\omega(D)$, x$(D^*)=$x$(D)$, and
$n(\partial D^*)=n(\partial D)-1$.
Hence $\tau(D^*)<\tau(D)$.
This contradicts the fact that 
the NS-tangle $(\Gamma\cap D,D)$ is minimal.

Therefore $\partial D$ does not intersect any terminal edge.
\end{proof}

\begin{lemma}
\label{MinimalNS-tangleNoHoop}
If there exists an NS-tangle 
in a minimal chart, 
then there exist a minimal chart $\Gamma$ and 
 a minimal NS-tangle 
$(\Gamma\cap D,D)$ such that 
$D$ does not contain any hoop.
\end{lemma}

\begin{proof}
Suppose that there exists an NS-tangle 
in a minimal chart $\Gamma$.
Let $(\Gamma\cap D^*,D^*)$ be 
a minimal NS-tangle of $\Gamma$. 
By Assumption~\ref{AssumptionFreeEdge} and
Assumption~\ref{AssumptionFreeEdgeSimpleHoop},
we can assume that 

(1) the disk $D^*$ does not contain any simple hoop.

Suppose that $D^*$ contains a hoop.
Let $C$ be an innermost hoop in $D^*$.
Let $E$ be the disk bounded by the hoop $C$, 
and 
$A$ a regular neighbourhood of $C$.
Set $D=Cl(E-A)$. 
Then $\Gamma\cap \partial D=\emptyset$.
Since $C$ is an innermost hoop in $D^*$,
the disk $D$ does not contain any hoop.
Since $C$ is not simple by (1),
the disk $E$ contains a white vertex
and so does $D$.
Thus
 $(\Gamma\cap D,D)$ is 
an NS-tangle. 
Since $(\Gamma\cap D^*,D^*)$ is 
a minimal NS-tangle, so is 
$(\Gamma\cap D,D)$.
Hence the tangle $(\Gamma\cap D,D)$ is 
a desired tangle.
\end{proof}

Let $C$ be a cycle of label $m$ in a chart $\Gamma$, and
$E$ the disk bounded by $C$. Let

$E^*=E\cup(\cup
\{e~|~\text{$e$ is a terminal edge intersecting $E$}
\})$.\\
Let $D^*$ be 
a regular neighbourhood of $E^*$.
Since any terminal edge does not contain a crossing,
the set $D^*$ is a disk.
Thus $(\Gamma\cap D^*,D^*)$ is a tangle.
The tangle $(\Gamma\cap D^*,D^*)$ is called 
a {\it  tangle induced from the cycle $C$}.

\begin{remark}
\label{RemTangleInducedFromC}
{\rm
\begin{enumerate}
\item[$(1)$] Int $D^*$ contains all terminal edges intersecting $C$, and
$\partial D^*$ does not intersect any terminal edge. 
\item[$(2)$]  If an edge intersects $\partial D^*$, then it must intersect $C$, 
because $D^*$ is 
a regular neighborhood of $E^*$.
\end{enumerate}
}
\end{remark}

For a subset $X$ of a chart $\Gamma$, we define

$
\overline \alpha(X)=
\min\{~i~|~\Gamma_i\cap X\neq\emptyset~\},~~~~
\overline{\beta}(X)=
\max\{~i~|~\Gamma_i\cap X\neq\emptyset~\}.
$\\

{\it Proof of Theorem~\ref{NoNS-tangle}.}
Suppose that
there exists an NS-tangle 
in a minimal chart.
By Lemma~\ref{MinimalNS-tangle},  
Lemma~\ref{MinimalAdmissibleNS-tangle}
and Lemma~\ref{MinimalNS-tangleNoHoop},
there exist a minimal chart $\Gamma$
and a minimal NS-tangle 
$(\Gamma\cap D,D)$ of label $m$ such that

$(1)$ $D$ contains neither hoop nor ring,

$(2)$ $\partial D$ does not intersect any terminal edge.\\
By Assumption~\ref{AssumptionFreeEdge}
and Assumption~\ref{AssumptionFreeEdgeSimpleHoop},

$(3)$ $D$ does not intersect any free edge.\\
Let $\alpha=\overline\alpha(\Gamma\cap D),
\beta=\overline\beta(\Gamma\cap D)$.
Then $\alpha\neq m$ or $\beta\neq m$.
Without loss of generality we can assume 
that $\alpha\neq m$.
By Condition (i) of 
an NS-tangle, we have 
$|\Gamma_\alpha\cap \partial D|\le 1$.
Let $G_\alpha$ be a small component of $\Gamma_\alpha\cap D$.
Then $|G_\alpha\cap \partial D|\le 1$.
Thus by (1), (2) and (3),
the set $G_\alpha$ contains a white vertex.
Hence by Lemma~\ref{MinimalC1-tangle},
there exists a suspicious cycle $C^*$ 
of label $\alpha$ in $D$.
Let $(\Gamma\cap D^*,D^*)$ 
be a tangle induced from $C^*$.
By Remark~\ref{RemTangleInducedFromC}$(1)$,

$(4)$ $\partial D^*$ does not intersect any terminal edge.

{\bf Claim~1.}
The tangle 
$(\Gamma\cap D^*,D^*)$ 
is an NS-tangle of label $\alpha+1$.

{\it Proof of Claim~$1$}. 
Since $C^*$ is a suspicious cycle, 
we have

$(5)$ the cycle $C^*$ contains 
a white vertex.\\
Let $i$ be a label with $i\ge \alpha+2$.
Since $C^*\subset D$,
the intersection $\Gamma_i\cap C^*$ 
consists of
at most one crossing 
by Condition (iii) of an NS-tangle.
If an edge of label $i$ 
intersects $\partial D^*$,
then it must intersect the cycle $C^*$
by Remark~\ref{RemTangleInducedFromC}$(2)$.
Hence there exists at most one edge 
of label $i$ intersecting $\partial D^*$.
Thus

$(6)$ for any label $i\ge\alpha+2$ 
the intersection 
$\Gamma_i\cap\partial D^*$ is 
at most one point.\\
Further, since $C^*$ is 
a suspicious cycle 
of label $\alpha$, 
by Remark~\ref{RemTangleInducedFromC}$(1)$ 
and Condition (i) 
of a suspicious cycle we have

$(7)$ the intersection 
$\Gamma_\alpha\cap \partial D^*$ 
is at most one point.\\
Since $\alpha$ is the lowest label in $D$, 

$(8)$ for any label $j<\alpha$, 
we have 
$\Gamma_j\cap \partial D^*=\emptyset$.\\
Since $(\Gamma\cap D,D)$ is an NS-tangle,
for each label $i$,
$\Gamma_i\cap D$ contains at most one crossing, and since $D^*\subset D$, we have

$(9)$ for each label $i$,
the intersection $\Gamma_i\cap D^*$ contains 
at most one crossing.\\
The tangle
$(\Gamma\cap D^*,D^*)$ is an NS-tangle
of label $\alpha+1$ in $D$
by $(5)$, $(6)$, $(7)$, $(8)$ and $(9)$. 
Hence Claim~$1$ holds.

{\bf Claim~2.} 
$\Gamma\cap D^{*}\subset \Gamma_\alpha\cup \Gamma_{\alpha+1}$.

{\it Proof of Claim~$2$.}
Suppose that there exists an integer 
$r\neq \alpha,\alpha+1$
with $\Gamma_r\cap D^{*}\neq\emptyset$.
Then $r>\alpha+1$,
because $\alpha$ is the lowest label in $D$. 
Let $\beta^*=\overline\beta(\Gamma\cap D^{*})$.
Then $\beta^*\ge r> \alpha+1$.
Thus $\Gamma_{\beta^*}\cap \partial D^{*}$
consists of at most one point by $(6)$.
Let $G_{\beta^*}$ be a small component of  $\Gamma_{\beta^*}\cap D^*$. 
Then $|G_{\beta^*}\cap \partial D^*|\le 1$.
Thus by $(1)$, $(3)$ and $(4)$, 
the set $G_{\beta^*}$ contains a white vertex.
Hence Lemma~\ref{MinimalC1-tangle} assures that
there exists a suspicious cycle $C^{**}$
of label $\beta^*$ in $D^*$.
Let $(\Gamma\cap D^{**},D^{**})$ 
be a tangle induced from $C^{**}$.
Since $\beta^{*}\neq\alpha,\alpha+1$,
the cycle $C^{**}$ is contained
in the interior of the disk 
bounded by $C^*$.
Hence $(5)$ implies that 

$(10)$ $\omega(D^{**})<
\omega(D^*)\le\omega(D)$.\\
Now we can show that 
the tangle $(\Gamma\cap D^{**},D^{**})$ 
is an NS-tangle of label $\beta^*-1$ 
by a similar way as the one used to 
show Claim~$1$.
Now by $(10)$ we have $\tau(D^{**})<\tau(D)$.
This contradicts the fact that
the NS-tangle $(\Gamma\cap D,D)$ is 
minimal.
Hence Claim~$2$ holds.

Therefore the suspicious cycle $C^*$
bounds a $2$-color disk.
This contradicts 
Lemma~\ref{No2ColorC1-tangle}.
This proves Theorem~\ref{NoNS-tangle}.
{\hfill {$\square$}\vspace{1.5em}}

Let $\Gamma$ be a chart, and
$(\Gamma\cap D,D)$ a tangle. 
Let $m$ be a label of $\Gamma$ 
with $\Gamma_m\cap D\neq\emptyset$,
and 
$G$ a connected closed subset of 
$\Gamma_m\cap D$. Set 

$X=\cup \{E~|~
\text{$E$ is a disk bounded by a cycle 
in $G$}\}$.\\
Let
$G^*=G\cup X$. 
Then 
$G^*$ is simply connected.
Let
$D^*$ be a regular neighbourhood of $G^*$ in $D$.
Then $D^*$ is a disk.
Thus  
$(\Gamma\cap D^*,D^*)$ is 
a tangle called 
a {\it tangle induced from $G$
with respect to $D$.}

\begin{remark}
\label{RemTangleInducedFromG}
{\rm 
\begin{enumerate}
\item[$(1)$] 
If an edge intersects $\partial D^*$,
then the edge must intersect $G$, 
because $D^*$ is 
a regular neighborhood of $G^*=G\cup X$ in $D$.
\item[$(2)$] 
Suppose that $G$ is a connected component of $\Gamma_m\cap D$.
Since $D^*$ is a regular neighbourhood 
of $G^*$ in $D$, we have 
$|\Gamma_m\cap\partial D^*|=
|G\cap\partial D^*|=
|G\cap\partial D|.$
\end{enumerate}
}
\end{remark}

Let $\Gamma$ be a chart. 
A tangle $(\Gamma\cap D,D)$ 
is  said to be 
{\it $2$-color} if 
there exist two labels $m,k$ 
with $|m-k|=1$ and
$\Gamma\cap D\subset \Gamma_m\cup\Gamma_k$.

The following lemma will be used in the proof of Theorem~\ref{Prickles} and Theorem~\ref{TwoColorGateTangle}.

\begin{lemma}
\label{ConnectedCycle}
Let $\Gamma$ be a minimal chart,
and $m$ a label of $\Gamma$. 
Then we have the following:
\begin{enumerate} 
\item[{\rm (a)}] If $(\Gamma\cap D,D)$ is a $2$-color admissible tangle with $|\Gamma_m\cap \partial D|=2$,
then $\Gamma_m\cap D$ is connected. 
\item[{\rm (b)}] If $E$ is a $2$-color disk 
with $\partial E\subset\Gamma_m$
but without free edges nor simple hoops, 
then $\Gamma_m\cap E$ is connected.
\end{enumerate}
\end{lemma}

\begin{proof}
We show Statement (a).
Suppose that $\Gamma_m\cap D$ is not connected.
Then $|\Gamma_m\cap \partial D|=2$ implies that
 
(1) there exists a connected component $G$ of 
$\Gamma_m\cap D$ with 
$|G\cap \partial D|\le 1$.\\
Let
$(\Gamma\cap D',D')$ be 
a tangle induced from $G$ with respect to $D$.

{\bf Claim.}
The disk $D'$ contains 
at least one white vertex.

{\it Proof of Claim.}
If $|G\cap \partial D|=1$,
then 
the point $G\cap \partial D$
is contained in 
an internal edge $\overline e$
by Condition (ii) of an admissible tangle.
Hence $G$ contains a white vertex 
by Condition (iii) of an admissible tangle.
Thus $D'$ contains a white vertex.

Suppose that $|G\cap \partial D|=0$.
If $G$ contains a white vertex, then
$D'$ contains a white vertex.
If $G$ does not contain any white vertex, then
$G$ is a hoop or a ring or a free edge. 
By Condition (i) of an admissible tangle,
the set $G$ is neither free edge nor simple hoop.
Hence $G$ is a non-simple hoop or a ring.
By the definition of a simple hoop and 
Remark~\ref{Assumption0}(3), 
the curve $G$ bounds a disk containing a white vertex. 
Hence $D'$ contains a white vertex.
Thus Claim holds.

By Remark~\ref{RemTangleInducedFromG}$(2)$ 
and Statement (1),
we have $|\Gamma_m\cap \partial D'|=|G\cap \partial D|\le 1$. Thus

(2) $\Gamma_m\cap \partial D'$ is 
at most one point.\\
Since $(\Gamma\cap D,D)$ is a $2$-color tangle,
there exists a label $k$ with $|m-k|=1$
and $\Gamma\cap D\subset\Gamma_m\cup\Gamma_k$.
Thus

(3) for any label $i$ with 
$i\not=m$ and $i\not=k$,
we have 
$\Gamma_i\cap \partial D'=\emptyset$.\\
Since $(\Gamma\cap D,D)$ is a 2-color tangle,
the disk $D$ does not contain any crossing.
Hence

(4) for each label $i$,
the intersection $\Gamma_i\cap D'$ does not contain any crossing.\\
Thus 
by Claim, Statement (2), (3) and (4),
we have that $(\Gamma\cap D',D')$ is an NS-tangle of label $k$.
This contradicts Theorem~\ref{NoNS-tangle}.
Hence Statement (a) holds.

We show  Statement (b).
Suppose that 
$\Gamma_m\cap E$ is not connected. 
Let $G$ be a connected component 
of $\Gamma_m\cap E$ 
with $G\cap \partial E=\emptyset$.
Let $D$ be a regular neighbourhood of $E$.
Then $(\Gamma\cap D,D)$ is a 2-color tangle
with $G\cap \partial D=\emptyset$.
Let $(\Gamma\cap D',D')$ be a tangle induced from $G$
with respect to $D$.
In a similar way to (a),
we can show that
$(\Gamma\cap D',D')$ is an NS-tangle.
This contradicts Theorem~\ref{NoNS-tangle}.
\end{proof}


\section{Inward pseudo paths and outward pseudo paths}
\label{s:IO-Path}

In this section
we investigate a $2$-color disk $E$
with $\partial E\subset \Gamma_m$ 
and an arc on the boundary of $E$ for a minimal chart $\Gamma$.

\begin{figure}
\begin{center}
\includegraphics{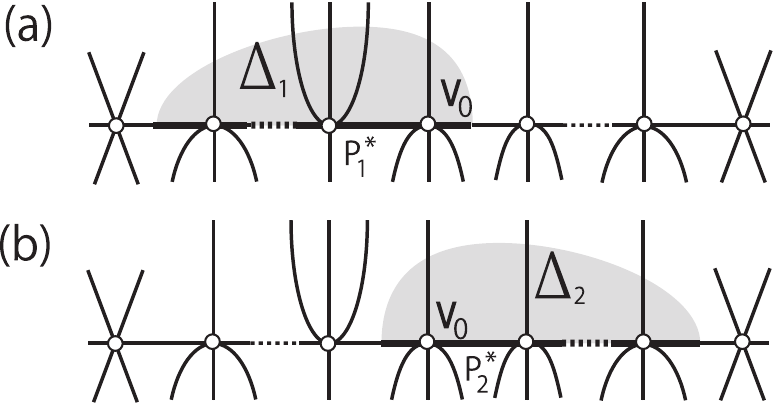}
\end{center}
\caption{\label{fig20} 
The thick lines are pseudo paths $P_1^*,P_2^*$.
}
\end{figure}

Let $\Gamma$ be a chart, and
$P_1^*,P_2^*$ 
pseudo paths of label $m$ in $\Gamma$.
The pair $(P_1^*,P_2^*)$ is called an
{\it I/O pair of type I}
if 
there exist
side-disks $\Delta_1,\Delta_2$ of $P_1^*,P_2^*$ respectively 
(see Fig.~\ref{fig20}) such that
\begin{enumerate}
\item[(i)]
$P_1^*,P_2^*$ are I/O pseudo paths with respect to 
the side-disks $\Delta_1,\Delta_2$ respectively,
\item[(ii)]
$P_1^*\cup P_2^*$ is a pseudo path, 
\item[(iii)] 
$\Delta_1\cup \Delta_2$ is a side-disk of $P_1^*\cup P_2^*$, and
\item[(iv)]
the intersection $P_1^*\cap P_2^*$ contains exactly one white vertex.
\end{enumerate}
We denote the pseudo path $P_1^*\cup P_2^*$ 
by $P_1^* *P_2^*$.
The union $\Delta_1\cup \Delta_2$ is called an 
{\it associated side-disk} of $P_1^* *P_2^*$.
The pair $(\Delta_1,\Delta_2)$ is called 
an {\it associated side-disk pair} of 
the I/O pair $(P_1^*,P_2^*)$.

\begin{remark}
{\rm 
Any vertex in an I/O pseudo path is a white vertex
by the definition of I/O pseudo paths.}
\end{remark}

\begin{lemma}\label{I/OPath1}
Let $\Gamma$ be a chart, and
$P^*_1,P^*_2$ pseudo paths of label $m$ in $\Gamma$.
If $(P^*_1,P^*_2)$ is an I/O pair of type I,
then $P^*_1 * P^*_2$ is an I/O pseudo path with respect to its associated side-disk.
\end{lemma}

\begin{proof}
We use the notations in the definition of 
an I/O pair of type I.
Without loss of generality we can assume 
that the pseudo path 
$P_1^*$ is inward with respect to 
the side-disk $\Delta_1$.
Then 
for each vertex $v$ in $P^*_1$,
any side-arc at $v$ of $P^*_1$ 
with respect to $\Delta_1$ 
is oriented inward at $v$. 
Thus
we have 
\begin{enumerate}
\item[$(1)$] 
for each vertex $v$ in $P^*_1\subset P_1^* *P_2^*$,
any side-arc at $v$ of $P_1^* *P_2^*$ 
with respect to $\Delta_1\cup \Delta_2$ 
is oriented inward at $v$. 
\end{enumerate}
By Condition (iv) for 
the definition of an I/O pair of type I,
there exists exactly one vertex $v_0$ in 
$P_1^*\cap P_2^*$.
Since 
the white vertex $v_0$ is 
a common vertex of $P_1^*$ and $P_2^*$, 
any side-arc at $v_0$ of $P^*_2$ 
with respect to $\Delta_2$
is oriented inward at $v_0$, too.
Hence the I/O pseudo path $P^*_2$ is 
inward with respect to $\Delta_2$.
Thus for each vertex $v$ in $P_2^*$,
any side-arc at $v$ of $P^*_2$ 
with respect to $\Delta_2$ 
is oriented inward at $v$. 
Hence 
\begin{enumerate}
\item[$(2)$] for each vertex $v$ in 
$P^*_2\subset P_1^* *P_2^*$,
any side-arc at $v$ of $P_1^* *P_2^*$ 
with respect to $\Delta_1\cup \Delta_2$ 
is oriented inward at $v$.
\end{enumerate}
Thus by $(1)$ and $(2)$, 
the pseudo path $P^*_1 * P^*_2$ is 
inward with respect to $\Delta_1\cup \Delta_2$.
\end{proof}



\begin{figure}
\begin{center}
\includegraphics{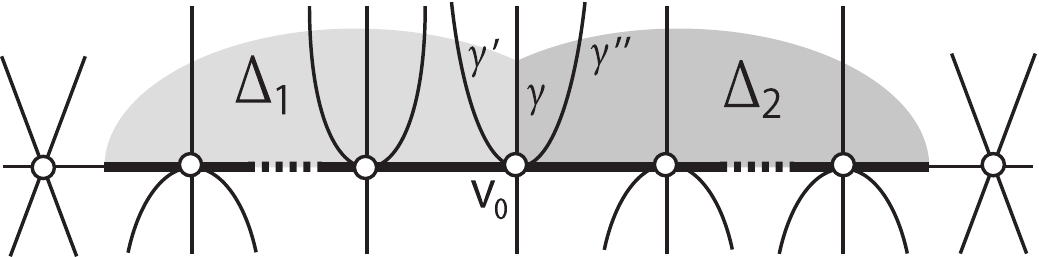}
\end{center}
\caption{ \label{fig21} 
The thick line is a pseudo path $P^*$. 
The end-arc $\gamma$ is a middle arc at $v_0$. 
And $\Delta_1,\Delta_2$ are side-disks of $P_1^*,P_2^*$ respectively.
}
\end{figure}

Let $\Gamma$ be a chart, and
$P_1^*,P_2^*$ pseudo paths of label $m$
in $\Gamma$.
The pair $(P_1^*,P_2^*)$ is called 
an {\it I/O pair of type II}
if there exist side-disks $\Delta_1,\Delta_2$ of 
$P_1^*,P_2^*$ respectively 
(see Fig.~\ref{fig21}) such that   
\begin{enumerate}
\item[(i)] 
the pseudo paths $P^*_1$ and $P^*_2$ 
are I/O pseudo paths with respect to 
the side-disks $\Delta_1,\Delta_2$ respectively,
\item[(ii)]
the intersection 
$\gamma=\Delta_1\cap \Delta_2$ is an end-arc of 
$P_1^*$ and $P_2^*$ respectively,
\item[(iii)] 
$\gamma$ is middle at 
the white vertex $v_0$ in $\gamma$,
\item[(iv)] 
$P^*=Cl((P_1^*\cup P_2^*)-\gamma)$ is a pseudo path, 
and $P^*\cap\gamma=v_0$ and
\item[(v)]
the union $\Delta_1\cup \Delta_2$ is a side-disk of $P^*$.
\end{enumerate}
We denote $P^*$ by $P^*_1 * P^*_2$.
The union $\Delta_1\cup \Delta_2$ is called 
an {\it associated side-disk} of $P_1^* *P_2^*$.
The pair $(\Delta_1,\Delta_2)$ is called 
an {\it associated side-disk pair} of 
the I/O pair $(P_1^*, P_2^*)$.

\begin{lemma}\label{I/OPath2}
Let $\Gamma$ be a chart.
Let $P^*_1,P^*_2$ be pseudo paths of label $m$ 
in $\Gamma$.
If $(P^*_1,P^*_2)$ is an I/O pair of type II,
then $P^*_1 * P^*_2$ is an I/O pseudo path
with respect to its associated side-disk.
\end{lemma}

\begin{proof}
We use the notations in the definition of 
an I/O pair of type II.
The end-arc $\gamma$ is a side-arc of $P^*_1*P_2^*$ 
middle at the vertex $v_0$.
Without loss of generality we can assume that
the side-arc $\gamma$ is oriented 
inward at $v_0$.
Let $\gamma'$ be a side-arc at $v_0$ of $P^*_1$ with respect to $\Delta_1$,
and $\gamma''$ a side-arc at $v_0$ of $P^*_2$ with respect to $\Delta_2$ (see Fig.~\ref{fig21}).
Since $\gamma$ is middle at $v_0$,
the side-arcs $\gamma',\gamma''$ are oriented inward at $v_0$.
Thus the pseudo paths $P^*_1$ and $P^*_2$ are 
inward with respect to 
$\Delta_1$ and $\Delta_2$ respectively.
Hence  
\begin{enumerate}
\item[$(1)$] for each vertex $v$ in $P_1^*$,
if $v\not=v_0$,
then any side-arc at $v$ of $P^*_1*P_2^*$ 
with respect to $\Delta_1\cup \Delta_2$ 
is oriented inward at $v$, and
\item[$(2)$] for each vertex $v$ in $P_2^*$,
if $v\not=v_0$,
then any side-arc at $v$ of $P^*_1*P_2^*$ 
with respect to $\Delta_1\cup \Delta_2$ 
is oriented inward at $v$.
\end{enumerate}
The three side-arcs $\gamma',\gamma,\gamma''$ are 
oriented inward at $v_0$.
Therefore by $(1)$ and $(2)$, 
the pseudo path $P^*_1 * P^*_2$ is inward
with respect to $\Delta_1\cup \Delta_2$.
\end{proof}

Let $\Gamma$ be a chart, 
$m$ a label of $\Gamma$, 
and $s\in\Bbb N$ with $s\ge 2$.
Let
$P^*$ be a pseudo path of label $m$ in $\Gamma$ with a side-disk $\Delta$, and
$P^*_1,P^*_2,\cdots,P^*_s$ pseudo paths 
of label $m$ in $\Gamma$ such that 
\begin{enumerate}
\item[(i)] for $i,j\in\{1,2,\cdots,s\}$,
if $|i-j|>1$, then
$P^*_i\cap P^*_j=\emptyset$,
\item[(ii)] for each $k=1,2,\cdots,s$,
$L_k=\cup_{i=1}^k L(P^*_i)$ is a path in $P^*$, 
here $L(P^*_i)$ is the maximal path 
contained in the pseudo path $P^*_i$ 
(see Section~\ref{s:Dichromatic})
and
\item[(iii)]
$L_s=L(P^*)$, here $L(P^*)$ is the maximal path 
contained in $P^*$.
\end{enumerate}
The s-tuple $(P^*_1,P^*_2,\cdots,P^*_s)$
is called an {\it I/O sequence} 
for $(P^*,\Delta)$
if there exist
side-disks $\Delta_1,\Delta_2,\cdots,\Delta_s$ of 
$P_1^*,P_2^*,\cdots,P_s^*$ respectively
such that
\begin{enumerate}
\item[(i)] 
for each $k=1,2,\cdots,s$,
$\Delta'_k=\cup_{i=1}^k \Delta_i$ is a disk in $\Delta$,
\item[(ii)]
for each $i=1,2,\cdots,s-1$,
the pair $(P^*_i,P^*_{i+1})$ is an I/O pair of 
type I or type II with 
an associated side-disk pair $(\Delta_i,\Delta_{i+1})$, and
\item[(iii)] $\Delta'_s$ is a side-disk of $P^*$.
\end{enumerate}
The s-tuple $(\Delta_1,\Delta_2,\cdots,\Delta_s)$
is called an {\it associated side-disk sequence} of 
the I/O sequence $(P^*_1,P^*_2,\cdots,P^*_s)$.
\begin{remark}
\label{NoteI/Osequence}
{\rm
Let $(P^*_1,P^*_2,\cdots,P^*_s)$ be 
an I/O sequence for $(P^*,\Delta)$
with an associated side-disk sequence
$(\Delta_1,\Delta_2,\cdots,\Delta_s)$.
Let $Q_1^*=P_1^*$.
Since the side-disk $\Delta$ of 
the pseudo path $P^*$ contains
all the side-disks $\Delta_1,\Delta_2,\cdots,\Delta_s$ of 
$P^*_1,P^*_2,\cdots,P^*_s$,
we can inductively show that
for each $i=2,3,\cdots,s$,
the pair $(Q_{i-1}^*,P^*_i)$ is 
an I/O pair of type I 
(resp. type II) with an associated side-disk pair $(\Delta'_{i-1},\Delta_{i})$
if the pair $(P^*_{i-1},P^*_{i})$ is 
an I/O pair of type I (resp. type II),
and that
$Q_i^*=Q_{i-1}^* * P^*_i$ is an I/O pseudo path
with respect to $\Delta'_i$
by Lemma~\ref{I/OPath1}
(resp. Lemma~\ref{I/OPath2}).
Hence $Q_s^*$ is an I/O pseudo path.
Therefore 
 $P^*$ is an I/O pseudo path.
We denote the I/O pseudo path $Q_s^*$ by 
$P^*_1 * P^*_2* \cdots* P^*_s$.
}
\end{remark}

\begin{lemma}
\label{OneSideLemma2}
Let $\Gamma$ be a minimal chart.
Let 
$P^*$ 
be a dichromatic one-side pseudo path 
of label $m$ in $\Gamma$ with
the associated vertex sequence $(v_1,v_2,\cdots,v_{s})$ and
an associated side-arc sequence 
$(\gamma_1',\gamma_2',\cdots,\gamma_{s}')$. Suppose that 
the end-arcs $\gamma_0$ and $\gamma_{s}$ of $P^*$ are 
middle at $v_1$ and $v_{s}$ 
respectively.
Suppose that there exists an integer $t$
with $2\le t\le s-1$ such that
\begin{enumerate}
\item[{\rm $($i$)$}]
the side-arc $\gamma_t'$ is middle at $v_t$, and
\item[{\rm $($ii$)$}]
for each $i=1,2,\cdots,s$ with 
$i\neq t$,
the side-arc $\gamma_i'$ is not middle at $v_i$.
\end{enumerate}
Then the pseudo path $P^*$ is 
an I/O pseudo path.
\end{lemma}

\begin{proof}
Let $(e_1,e_2,\cdots,e_{s-1})$ be the associated edge sequence for $P^*$.
Let $\gamma_{t-1}$ and $\gamma_{t}$ be short arcs containing $v_t$ in the edges $e_{t-1}$ and $e_{t}$
respectively.
Set $P^*_1=\gamma_0\cup e_1\cup e_2\cdots\cup e_{t-1}\cup \gamma_{t}$, and
$P^*_2=\gamma_{t-1}\cup e_{t}\cup \cdots\cup e_{s-1}\cup \gamma_s$.
Then $P^*_1$ and $P^*_2$ are dichromatic one-side pseudo paths
by Remark~\ref{NoteOneSidePath}$(2)$.
Thus $P^*_1$ and $P^*_2$ are I/O pseudo paths
by Lemma~\ref{OneSideLemma0}.
Thus $(P^*_1,P^*_2)$ is an I/O pair of type I. 
Therefore 
the pseudo path $P^*$ is an I/O pseudo path 
by Lemma~\ref{I/OPath1}.
\end{proof}


\begin{theorem}
\label{Prickles} 
Let $\Gamma$ be a minimal chart.
Let $C$ be a cycle of label $m$ in $\Gamma$
bounding a $2$-color disk $E$ such that 
the disk $E$ contains neither free edge nor simple hoop. 
If there exist two paths $S,T$ in
the path decomposition $\mathcal{P}
(C;\mathcal{W}_O(C,m)-\mathcal{W}_O^{{\rm Mid}}(C,m))$
with 
$\mathcal{W}_O^{{\rm Mid}}(C,m)\subset S\cup T$, 
then
\begin{enumerate}
\item[{\rm $($a$)$}] each of the paths $S$ and $T$ contains at least one white vertex in $\mathcal{W}_O^{{\rm Mid}}(C,m)$,
\item[{\rm $($b$)$}]
the extended pseudo paths 
$\widehat{S},\widehat{T}$
of ${S},{T}$
are I/O pseudo paths.
\end{enumerate}
\end{theorem}

\begin{figure}
\begin{center}
\includegraphics{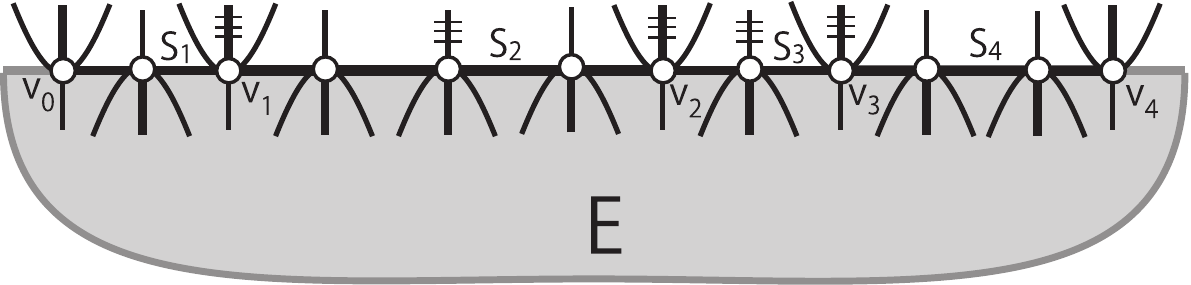}
\end{center}
\caption{ \label{fig22} 
The thicken lines are of label $m$. 
Each arc with three transversal short 
arcs is a middle arc at the white vertex.
For paths 
$S_1,S_2,S_3,S_4$
of label $m$, 
$\partial S_1=\{v_0,v_1\},
\partial S_2=\{v_1,v_2\},
\partial S_3=\{v_2,v_3\},
\partial S_4=\{v_3,v_4\}$. }
\end{figure}

\begin{proof}
By Lemma~\ref{ConnectedCycle}(b),
the intersection $\Gamma_m\cap E$ is connected.
Let $s=|S\cap\mathcal{W}_O^{{\rm Mid}}(C,m)|$, and
$t=|T\cap\mathcal{W}_O^{{\rm Mid}}(C,m)|$.
Suppose that $s=0$ or $t=0$. 
Then 
$\mathcal{W}_O^{{\rm Mid}}(C,m)\subset S\cup T$ 
implies that 
$\mathcal{W}_O^{{\rm Mid}}(C,m)$ 
is contained in 
one of the two paths $S,T$ in 
$\mathcal{P}(C;
\mathcal{W}_O(C,m)-
\mathcal{W}_O^{{\rm Mid}}(C,m))$.
This contradicts Lemma~\ref{qPrickle}(c).
Thus we have $s\ge 1,t\ge 1$.
Namely Statement (a) holds.

We show Statement (b).
Let
\begin{enumerate}
\item[]
$\mathcal{S}=\mathcal{P}(S;\mathcal{W}_O^{{\rm Mid}}(C,m)\cap S)
=
\{S_1,S_2,\cdots,S_{s+1}\}$, and 
\item[]
$\mathcal{T}=\mathcal{P}(T;\mathcal{W}_O^{{\rm Mid}}(C,m)\cap T)
=
\{T_1,T_2,\cdots,T_{t+1}\}$.
\end{enumerate}
Since
$S,T\in \mathcal{P}(C;\mathcal{W}_O(C,m)-\mathcal{W}_O^{{\rm Mid}}(C,m))$, we have
\begin{enumerate}
\item[$(1)$] 
$
\{S_1,S_2,\cdots,S_{s+1},T_1,T_2,\cdots,T_{t+1}\}
\subset \mathcal{P}(C;\mathcal{W}_O(C,m))$.
\end{enumerate}
Further,
without loss of generality we can assume that
\begin{enumerate}
\item[$(2)$] 
each of $\partial S_1,\partial S_{s+1},\partial T_1,\partial T_{t+1}$
consists of a vertex in 
$\mathcal{W}_O(C,m)-\mathcal{W}_O^{{\rm Mid}}(C,m)$
and a vertex in 
$\mathcal{W}_O^{{\rm Mid}}(C,m)$, and
\item[$(3)$]
for $i=2,3,\cdots,s$ and $j=2,3,\cdots,t$,
we have 
$\partial S_i\subset \mathcal{W}_O^{{\rm Mid}}(C,m)$ and
$\partial T_j\subset \mathcal{W}_O^{{\rm Mid}}(C,m)$.
\end{enumerate}
Let

$X=\{~P~|~P\in \mathcal{S}
\cup \mathcal{T}
,~\partial P\subset \mathcal{W}_O^{{\rm Mid}}(C,m)\}
=\{S_2,S_3,\cdots,S_{s},
T_2,T_3,\cdots,T_{t}
\}$.\\
Since $\mathcal{W}_O^{{\rm Mid}}(C,m)\subset S\cup T$ by the assumption,
we have
\begin{enumerate}
\item[$(4)$] 
$s+t=|\mathcal{W}_O^{{\rm Mid}}(C,m)|$.
\end{enumerate}
Let $k$ be the secondary label of 
the $2$-color disk $E$. 
By Remark~\ref{O(C)I(C)}$(3)$, we have
\begin{enumerate}
\item[$(5)$] for each path $P$ in 
$\mathcal{S}\cup \mathcal{T}$,
a side-arc of label $k$ of 
the extended pseudo path $\widehat{P}$ of 
$P$
is middle at a vertex $v$ in Int~$P$
if and only if $v\in \mathcal{W}_I^{{\rm Mid}}(C,m)$.
\end{enumerate}
By Statement $(1)$, 
Lemma~\ref{DecompositionPaths} implies that
\begin{enumerate}
\item[$(6)$] 
for each path 
$P\in \mathcal{S}\cup\mathcal{T}$,
the extended pseudo path $\widehat{P}$
of $P$
is a dichromatic one-side pseudo path.
\end{enumerate}
Further, Statement $(3)$ assures that
for each path 
$P\in X$ we can apply 
Lemma~\ref{OneSideLemma}(b) 
to $\widehat{P}$ so that 
the extended pseudo path $\widehat{P}$
has 
a side-arc of label $k$ middle 
at a vertex in Int~$P$.
Namely by Statement $(5)$ 
\begin{enumerate}
\item[$(7)$] 
for each $P\in X$,
the set Int~$P$ contains at least 
one vertex in $\mathcal{W}_I^{{\rm Mid}}(C,m)$.
\end{enumerate}
For each $i=1,2,\cdots,s+1$ 
and $j=1,2,\cdots,t+1$,
let

$
\sigma_i=|S_i\cap\mathcal{W}_I^{{\rm Mid}}(C,m)|,~~~~
\tau_j=|T_j\cap\mathcal{W}_I^{{\rm Mid}}(C,m)|
$.\\
Then Statement $(7)$ implies
\begin{enumerate}
\item[$(8)$] 
for each $i=2,3,\cdots,s$ 
and $j=2,3,\cdots,t$,
we have $1\le \sigma_i$ and $1\le \tau_j$.
\end{enumerate}

{\bf Claim~1.} 
$\sigma_1=\sigma_{s+1}=\tau_1=\tau_{t+1}=0,$ 
and\\
\hspace{19.5mm}
$\sigma_2=\sigma_3=\cdots=\sigma_{s}=
\tau_2=\tau_3=\cdots=\tau_{t}=1.$

{\it Proof of Claim $1$.}
If not, then by Statement $(8)$, 
we have
\begin{enumerate}
\item[]
$\sigma_1>0,$ or $\sigma_{s+1}>0,$ or 
$\tau_1>0,$ or $\tau_{t+1}>0,$ or
\item[]
$\sigma_i>1$ for some $2\le i\le s,$ or
$\tau_j>1$ for some $2\le j\le t$.
\end{enumerate}
Thus we have
\begin{enumerate}
\item[$(9)$]
$\begin{array}{ll}
|\mathcal{W}_I^{{\rm Mid}}(C,m)|&\ge 
\sigma_1+(\sigma_2+\cdots+\sigma_{s})+\sigma_{s+1}
+\tau_1+(\tau_2+\cdots+\tau_{t})+\tau_{t+1}\\
&> 0~+(1~+\cdots+~1)+~0~+~0~+(1~+\cdots+~1)+0\\
&=(s-1)+(t-1)=s+t-2.\\
\end{array}$
\end{enumerate}
Using the equation $(4)$ and 
the inequality $(9)$,
Theorem~\ref{Inequation} implies 
$$s+t=(s+t-2)+2<|\mathcal{W}_I^{{\rm Mid}}(C,m)|+2\le |\mathcal{W}_O^{{\rm Mid}}(C,m)|=s+t.$$
This is a contradiction.
Therefore Claim $1$ holds.

Now Claim $1$ implies that

(10) each path in $X$ contains 
exactly one white vertex 
in $\mathcal{W}_I^{{\rm Mid}}(C,m)$.

{\bf Claim 2.} 
$C-(\cup_{P\in X} P)$ 
does not contain any white vertex in 
$\mathcal{W}_I^{{\rm Mid}}(C,m)$.

{\it Proof of Claim $2$.} 
If $C-(\cup_{P\in X} P)$ 
contains a white vertex in 
$\mathcal{W}_I^{{\rm Mid}}(C,m)$,
then by Statement (10)
we have the inequality 
$|\mathcal{W}_I^{{\rm Mid}}(C,m)|> |X|=s+t-2$.
Thus we have the same contradiction 
as the one of Claim 1. 
Thus Claim~$2$ holds.

Without loss of generality
we can assume that 
$S_1,S_2,\cdots,S_{s+1}$ are situated on 
$S$ in this order.
By Statement $(3)$, $(6)$ and (10),
Lemma~\ref{OneSideLemma2} 
implies that
the extended pseudo paths 
$\widehat{S_2},\widehat{S_3},
\cdots,\widehat{S_{s}}$
are I/O pseudo paths 
(see $S_2$ and $S_3$ in Fig.~\ref{fig22}).
Further, Claim $2$ 
and Statement $(5)$ imply that 
for any vertex $v\in$Int~$S_1$ 
(resp. $v\in$Int~$S_{s+1}$)
any side-arc at $v$ of $\widehat{S_{1}}$ 
(resp. $\widehat{S_{s+1}}$) 
is not middle at $v$. Hence by 
Statement $(2)$ and Statement $(6)$,
Lemma~\ref{OneSideLemma0} implies that
the extended pseudo paths 
$\widehat{S_1},\widehat{S_{s+1}}$
are I/O pseudo paths 
(see $S_1$ and $S_4$ in Fig.~\ref{fig22}).
Now for each $i=1,2,\cdots,s$,
we have that
$(\widehat{S_i},\widehat{S_{i+1}})$
is an I/O pair of type II.
Thus
$(\widehat{S_1}, \widehat{S_2},
\cdots,\widehat{S_{s}},
\widehat{S_{s+1}})$
is an I/O sequence
for $(\widehat{S},\Delta)$ 
here $\Delta$ is a side-disk of $\widehat{S}$. 
Hence by Remark~\ref{NoteI/Osequence}
the extended pseudo path 
$\widehat{S_1}*
\widehat{S_2}*
\cdots*\widehat{S_{s}}*
\widehat{S_{s+1}}$
is an I/O pseudo path 
and so is $\widehat{S}$.
Similarly we can show that
$\widehat{T}$ is an I/O pseudo path. 
This completes the proof of
Theorem~\ref{Prickles}.
\end{proof}

\section{Bridges}
\label{s:Bridge}

In this section, 
we investigate a path of label $m$ such that 
each white vertex in the path is contained in 
a terminal edge.

\begin{figure}
\begin{center}
\includegraphics{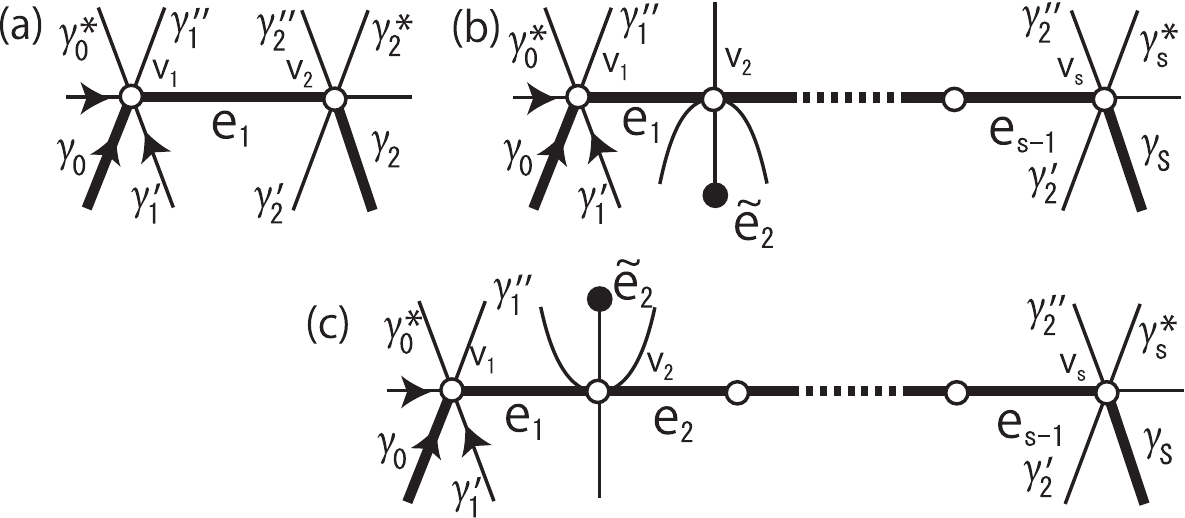}
\end{center}
\caption{ \label{fig23} The thicken lines are bridges $B$. }
\end{figure}

Let $\Gamma$ be a chart.
Let $B$
be an admissible pseudo path  of label $m$
 in $\Gamma$
 with the associated vertex sequence
$(v_1,v_{2},\cdots,v_s)$
 and the associated edge sequence
$(e_1,e_2,\cdots,e_{s-1})$.
The pseudo path $B$ is called a {\it bridge}
provided that
\begin{enumerate}
\item[(i)] $s\ge2$ and all the vertices in $B$ are white vertices,
\item[(ii)] 
the edge $e_1$ is not middle at $v_1$ and
the edge $e_{s-1}$ is not middle at $v_{s}$,
and
\item[(iii)] 
for each $i=2,3,\cdots,s-1$,
there exists a terminal edge of label $m$ 
at $v_i$.
\end{enumerate}
Let $\gamma_0$ and $\gamma_s$ be the end-arcs of $B$ with $\gamma_0\ni v_1$ and $\gamma_s\ni v_s$.
Let $\Delta$ be a side-disk for which $B$ is admissible.
Let $\gamma_0^*$ and $\gamma_s^*$ be short arcs
in edges of label $m$ in $\Gamma$
with $\gamma_0^*\cap \Delta=v_1$ and 
$\gamma_s^*\cap \Delta=v_{s}$. 
Then 
$B^*=\gamma_0^*\cup e_1\cup e_2\cup \cdots\cup e_{s-1}\cup\gamma_s^*$
is a bridge called the {\it co-bridge} 
of $B$ (see Fig.~\ref{fig23}). 
It is clear that the bridge $B$ is 
the co-bridge of $B^*$.
If there exists a label $k$ with
$|m-k|=1$ and
$v_1,\cdots,v_{s}\in \Gamma_m\cap\Gamma_k$,
then the bridge $B$ is called 
a {\it dichromatic bridge}.

\begin{lemma}
\label{2ColorBridge}
Let $\Gamma$ be a minimal chart.
Let $B$
be a dichromatic bridge 
of label $m$ in $\Gamma$, 
and $B^*$ a co-bridge of $B$.
Then
one of $B$ and $B^*$ 
is an inward pseudo path, and 
the other is an outward pseudo path.
\end{lemma}

\begin{proof}
We use the notations in the definition of 
a bridge and a co-bridge.
Since the edge $e_1$ is not middle at $v_1$
by Condition (ii) of a bridge,
one of the two end-arcs $\gamma_0$ and $\gamma_0^*$ is middle at $v_1$
by Remark~\ref{EdgesAroundVertex}$(2)$.
Since the co-bridge of $B^*$ is $B$,
without loss of generality we can assume that
\begin{enumerate}
\item[$(1)$] the end-arc $\gamma_0$ is middle at $v_1$.
\end{enumerate}
Further, we can assume that
\begin{enumerate}
\item[$(2)$] 
the end-arc $\gamma_0$ is oriented inward at $v_1$.
\end{enumerate}

First, we prove the case $s=2$.
Now $B=\gamma_0\cup e_1\cup\gamma_2$ and 
$B^*=\gamma_0^*\cup e_1\cup\gamma_2^*$ 
are dichromatic one-side pseudo paths 
(see Fig.~\ref{fig23}(a)).
Let $(\gamma_1',\gamma_2')$ and 
$(\gamma_1'',\gamma_2'')$ be 
 side-arc sequences 
of $B$ and $B^*$ respectively.
Then by $(1)$ and $(2)$, we have
\begin{enumerate}
\item[$(3)$] 
the side-arc $\gamma_1'$ is 
oriented inward at $v_1$,
and  the edge $e_1$ and 
the side-arc $\gamma_1''$ 
are oriented outward at $v_1$.
\end{enumerate}
Since the end-arc $\gamma_0$ 
is middle at $v_1$ by $(1)$, 
and
since $s=2$,
the bridge $B$ is an I/O pseudo path 
by Lemma~\ref{OneSideLemma0}.
Since the side-arc $\gamma_1'$ 
is oriented inward 
at $v_1$ by $(3)$, 
the bridge $B$ is an inward pseudo path.
Thus the side-arc $\gamma_2'$ is oriented inward at $v_2$.
Now $e_1$ is oriented 
inward at $v_2$ by $(3)$.
Since $e_1$ is not middle at $v_2$ 
by Condition (ii) of a bridge,
the side-arc $\gamma_2''$ is 
oriented outward at $v_2$.
Since the side-arc $\gamma_1''$ is 
oriented outward at $v_1$
by $(3)$,
the co-bridge $B^*$ is an outward pseudo path.

Suppose that $s\ge 3$.
Now for each $i=2,3,\cdots,s-1$
there exists a terminal edge 
$\widetilde{e}_i$ 
of label $m$ 
at $v_i$ by Condition (iii) of a bridge.
Since $\widetilde{e}_2$ is 
a terminal edge at $v_2$,
by Remark~\ref{Assumption0}$(2)$
we have
\begin{enumerate}
\item[$(4)$] 
the edge $\widetilde{e}_2$ is middle at $v_2$.
\end{enumerate}

If $\widetilde{e}_2$ contains a side-arc of $B$, say $\widetilde{\gamma}_2$,
then $P=\gamma_0\cup e_1\cup\widetilde{\gamma}_2$ is 
a dichromatic one-side pseudo path 
(see Fig.~\ref{fig23}(b)).
But the end-arc $\gamma_0$ of $P$
is middle at $v_1$ by (1), 
and the end-arc $\widetilde{\gamma}_2$ of $P$ is middle at $v_2$ by $(4)$.
Further, 
the dichromatic one-side pseudo path $P$ 
contains 
only two white vertices $v_1,v_2$.
This contradicts 
Lemma~\ref{OneSideLemma}(a).
Hence $\widetilde{e}_2$ contains 
a side-arc of $B^*$, say $\widetilde{\gamma}_2^*$ 
(see Fig.~\ref{fig23}(c)).
By $(4)$,
\begin{enumerate}
\item[$(5)$] 
the side-arc $\widetilde{\gamma}_2^*$ is middle at $v_2$.
\end{enumerate}
Let $\gamma_1$ (resp. $\gamma_2$) be 
a short arc in the edge $e_1$ (resp. $e_2$) 
containing $v_2$.
Let $B_1=\gamma_0\cup e_1\cup \gamma_2,
B_2=\gamma_1\cup e_2\cup\cdots\cup e_{s-1}\cup\gamma_s,
B_1^*=\gamma_0^*\cup e_1\cup\widetilde{\gamma}_2^*,
B_2^*=\widetilde{\gamma}_2^*\cup e_2\cup\cdots\cup e_{s-1}\cup\gamma_s^*$.
Since $\widetilde{\gamma}_2^*$ is middle at $v_2$ by $(5)$,
the edges $e_1$ and $e_2$ are 
not middle at $v_2$ by Remark~\ref{EdgesAroundVertex}$(2)$.
Thus $B_1,B_2,B_1^*,B_2^*$ are bridges.

Now $B_1^*,B_2^*$ are 
co-bridges of $B_1,B_2$ respectively.
Thus $B_1$ is inward 
and $B_1^*$ is outward by the case $s=2$.
By induction on 
the number of edges of a bridge,
we can show that 
$B_2$ and the co-bridge of $B_2$ 
are I/O pseudo paths.
Since 
$(B_1,B_2)$ is an I/O pair of type I,
the bridge $B$ is an I/O pseudo path
by Lemma~\ref{I/OPath1}.
Since $B_1$ is an inward pseudo path, 
so is $B$.
Since 
$(B_1^*,B_2^*)$ is an I/O pair of type II,
the bridge $B^*$ is an I/O pseudo path
by Lemma~\ref{I/OPath2}.
Since $B_1^*$ is an outward pseudo path, 
so is $B^*$.
This proves Lemma~\ref{2ColorBridge}.
\end{proof}

\begin{figure}
\begin{center}
\includegraphics{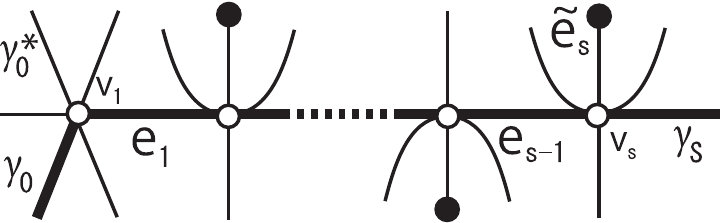}
\end{center}
\caption{ \label{fig24} The thick line is a pier $B$. }
\end{figure}

Let $\Gamma$ be a chart.
Let 
$B$
be a pseudo path of label $m$
 in $\Gamma$
with the associated vertex sequence
$(v_1,v_2,\cdots,v_s)$
and the associated edge sequence 
$(e_1,e_2,\cdots,e_{s-1})$.
Let $\gamma_0$ and $\gamma_s$ be the end-arcs of $B$ with $\gamma_0\ni v_1$ and $\gamma_s\ni v_s$.
The pseudo path $B$ is called a {\it pier}
provided that (see Fig.~\ref{fig24})
\begin{enumerate}
\item[(i)] all the vertices in $B$ are white vertices,
\item[(ii)]
the edge $e_1$ is not middle at $v_1$, and
\item[(iii)] 
for each $i=2,3,\cdots,s$,
there exists a terminal edge of label $m$ 
at $v_i$
which does not contain the end-arc $\gamma_s$.
\end{enumerate}
Let $\gamma_0^*$ be a short arc
in an edge of label $m$ in $\Gamma$
with $v_1=\gamma_0^*\cap e_{1}=\gamma_0^*\cap \gamma_{0}$. 
Then 
$B^*=\gamma_0^*\cup e_1\cup e_2\cup \cdots\cup e_{s-1}\cup\gamma_s$
is a pier called the {\it co-pier} 
of $B$. 
It is clear that the pier $B$ is 
the co-pier of $B^*$.
There exist side-disks $\Delta$ and $\Delta^*$  
of $B$ and $B^*$ respectively
with $\Delta\cap \Delta^*=B\cap B^*$.
The side-disk $\Delta$ (resp. $\Delta^*$) is called 
a {\it nice side-disk} for $B$ (resp. $B^*$).
If there exists a label $k$ with
$|m-k|=1$ and
$v_1,\cdots,v_{s}\in \Gamma_m\cap\Gamma_k$,
then the pier $B$ is called 
a {\it dichromatic pier}.

\begin{cor}
\label{2ColorPier}
Let $\Gamma$ be a minimal chart.
Let
$B$
be a dichromatic pier of label $m$ in $\Gamma$, and $B^*$
the co-pier of $B$.
Then 
one of $B$ and $B^*$ 
is an inward pseudo path 
with respect to a nice side-disk, and 
the other is an outward pseudo path 
with respect to a nice side-disk.
\end{cor}

\begin{proof}
We use the notations in the definition of a pier.
Let $\Delta$ and $\Delta^*$ be nice side-disks of
$B$ and $B^*$ respectively.
Our result is true for the case $s=1$,
because $e_1$ is not middle at $v_1$.
Suppose that $s\ge 2$.
Let $\widetilde{e}_s$ be the terminal edge of label $m$ at $v_s$ not containing the end-arc $\gamma_s$.
If $\widetilde{e}_{s}\cap$Int~$\Delta^*\not=\emptyset$, 
then $B$ is an admissible pseudo path.
If $\widetilde{e}_{s}\cap$Int~$\Delta\not=\emptyset$, 
then $B^*$ is an admissible pseudo path.
Since $B$ and $B^*$ are 
co-piers of each other,
without loss of generality we can assume that
$B$ is an admissible pseudo path
(see Fig.~\ref{fig24}).
Now 
the terminal edge $\widetilde{e}_{s}$ is middle at $v_{s}$ by Remark~\ref{Assumption0}$(2)$.
By Remark~\ref{EdgesAroundVertex}$(2)$,
the edge $e_{s-1}$ is not middle at $v_{s}$.
Thus $B$ is a bridge.
Let $\widetilde{\gamma}_{s}$ be a short arc in the terminal edge $\widetilde{e}_{s}$ of label $m$ with $v_{s}\in \widetilde{\gamma}_{s}$. 
Then $B^\dagger=\gamma_0^*\cup e_1\cup\cdots,
e_{s-1}\cup\widetilde{\gamma}_{s}$ 
is a co-bridge of $B$.
Thus $B^\dagger$ and $B$ are I/O pseudo paths
by Lemma~\ref{2ColorBridge}.
Without loss of generality we can assume that
$B$ is inward and $B^\dagger$ is outward.
Then a side-arc at $v_s$ of $B$ with respect to $\Delta$ is oriented inward at $v_s$.
Hence $\widetilde{e}_s$ is oriented outward at $v_s$.
Since the terminal edge $\widetilde{e}_{s}$ is middle at $v_{s}$,
the side-arcs at $v_s$ of $B^*$ with respect to $\Delta^*$ is oriented outward at $v_s$.
Thus the pier $B^*$ is an outward pseudo path with respect to $\Delta^*$.
This proves Corollary~\ref{2ColorPier}.
\end{proof}

\section{Proof of Theorem~1.3}

\label{s:PfTheorem3}


In this section
we shall prove Theorem~\ref{TwoColorGateTangle}
by using 
Theorem~\ref{Prickles}(b),
Lemma~\ref{2ColorBridge}, 
Corollary~\ref{2ColorPier}, and
Lemma~\ref{2ColorBeads} below.

\begin{lemma}
\label{2ColorBeads}
Let $\Gamma$ be a minimal chart, and
$m$ a label of $\Gamma$.
Let $(\Gamma\cap D,D)$ be
a $2$-color admissible tangle with $|\Gamma_m\cap\partial D|=2$.
If $\Gamma_m\cap D$ contains a cycle, 
then 
there exist disks $E_1,E_2,\cdots, E_d$ in {\rm Int}~$D$ and simple arcs $L_0,L_1,\cdots,L_d$ in 
$Cl(D-\cup^{d}_{i=1}E_i)$
 such that
\begin{enumerate}
\item[{\rm (a)}] 
$\partial E_i\subset \Gamma_m$ for each $i=1,2,\cdots,d$ and $L_j\subset \Gamma_m$ for each $j=0,1,\cdots,d$,
\item[{\rm (b)}] 
for each $j=1,\cdots,d-1$, $L_j$ connects a vertex in $\partial E_j$ and a vertex in $\partial E_{j+1}$, \\
the arc $L_0$ connects a point in $\partial D$ and a vertex in $\partial E_{1}$, and\\
the arc $L_d$ connects a vertex in $\partial E_{d}$ and a point in $\partial D$,
\item[{\rm (c)}]
if an edge of label $m$ intersects 
$D-((\cup_{i=1}^{d} E_i)\cup (\cup_{j=0}^{d} L_j))$, 
then it is a terminal edge.
\end{enumerate}
\end{lemma}

\begin{proof}
Since $(\Gamma\cap D,D)$ is a 2-color tangle,
\begin{enumerate}
\item[(1)]there exists a label $k$ with $|m-k|=1$ and
$\Gamma\cap D\subset\Gamma_m\cup\Gamma_k$, 
\item[(2)] the disk $D$ does not contain any crossing.
\end{enumerate}
By Lemma~\ref{ConnectedCycle}(a),
\begin{enumerate}
\item[(3)] $\Gamma_m\cap D$ is connected.
\end{enumerate}
Let 
$E_1,E_2,\cdots,E_d$ be all the disks 
bounded by the maximal cycles 
in $\Gamma_m\cap D$.
Since $\Gamma_m\cap D$ contains a cycle,
we have $d\ge1$.
For each $i=1,2,\cdots,d$, let
\begin{enumerate}
\item[]
$\widehat E_i=E_i\cup(\cup\{e~|~e 
\text{ is a terminal edge in $\Gamma_m\cap D$ 
intersecting }E_i\})$.
\end{enumerate}
Let $P_1,P_2,\cdots,P_{p}$ be 
all the closures of connected components of 
$(\Gamma_m\cap D)-\cup_{i=1}^d\widehat E_i$
not intersecting $\partial D$,
and 
$Q_1,Q_2,\cdots,Q_q$
all the closures of connected components of 
$(\Gamma_m\cap D)-\cup_{i=1}^d\widehat E_i$
intersecting $\partial D$. Then
\begin{enumerate}
\item[(4)] 
for each $i=1,2,\cdots,p$, 
$P_i\cap\partial D=\emptyset$, and 
the tree $P_i$ is not a terminal edge.
\end{enumerate}
By the condition $|\Gamma_m\cap\partial D|=2$, 
we have $q\le 2$.
Set
\begin{enumerate}
\item[]
$\mathcal{E}=\{E_1,E_2,\cdots,E_d\}$,
\item[]
$\mathcal{P}=\{P_1,P_2,\cdots,P_p\}$, and
\item[]
$\mathcal{Q}=\{Q_1,Q_2,\cdots,Q_q\}$.
\end{enumerate}
Since $\Gamma_m\cap D$ is connected by $(3)$
and since $|\Gamma_m\cap \partial D|=2$,
\begin{enumerate}
\item[(5)] 
for each $i=1,2,\cdots,d$, 
the cycle $\partial E_i$ contains 
a white vertex.
\end{enumerate}
Further, by (3) and (4), we have
\begin{enumerate}
\item[(6)] 
for each $i=1,2,\cdots,p$, 
the tree $P_i$ contains 
at least two white vertices.
\end{enumerate}

{\bf Claim 1.} For each $i=1,2,\cdots, d$, 
the disk $E_i$ intersects 
exactly two of $\mathcal{P}\cup\mathcal{Q}$.

{\it Proof of Claim $1$.}
If there exists a disk $E_i$ intersecting 
at most one of $\mathcal{P}\cup\mathcal{Q}$,
let $D'$ be a regular neighbourhood of $\widehat E_i$.
Then $\Gamma_m\cap\partial D'$ consists of 
at most one point.
Thus by (5),
the tangle $(\Gamma\cap D',D')$ is 
an NS-tangle of label $k$. 
This contradicts Theorem~\ref{NoNS-tangle}.

Suppose that 
there exists a disk $E_i$ intersecting 
at least three of $\mathcal{P}\cup\mathcal{Q}$.
Let $D'$ be a regular neighbourhood of $\widehat E_i$.
Then $Cl((\Gamma_m\cap D)-D')$ consists of 
at least three connected components.
Thus 
the condition $|\Gamma_m\cap \partial D|=2$
implies that
there exists a connected component $X$ of
$Cl((\Gamma_m\cap D)-D')$
with $X\cap\partial D=\emptyset$. 
Let $(\Gamma\cap D'',D'')$ be 
a tangle induced from $X$ with respect to $D$.
Then $\Gamma_m\cap\partial D''$ consists of 
one point.
Further, by (6),
the component $X$ contains 
at least one white vertex.
Thus $(\Gamma\cap D'',D'')$ is 
an NS-tangle of label $k$.
This contradicts Theorem~\ref{NoNS-tangle}.
Therefore Claim~$1$ holds.

{\bf Claim~2.} For each $i=1,2,\cdots,p,$ 
the tree $P_i$ intersects 
exactly two of $\mathcal{E}$.

{\it Proof of Claim~$2$.}
If there exists an element $P_i$ intersecting 
at most one of $\mathcal{E}$, 
then $P_i$ is a reducible tree. 
This contradicts 
Lemma~\ref{NoReducibleTree}.

Suppose that there exists an element $P_i$
intersecting at least three of $\mathcal{E}$.
Then $Cl((\Gamma_m\cap D)-P_i)$ consists of 
at least three connected components. 
Again, 
the condition $|\Gamma_m\cap \partial D|=2$ 
implies that
there exists a connected component $X$ 
of $Cl((\Gamma_m\cap D)-P_i)$
with $X\cap\partial D=\emptyset$.
Let $(\Gamma\cap D',D')$ be 
a tangle induced from $X$ with respect to $D$.
Then $\Gamma_m\cap\partial D'$ consists of 
one point.
Thus by (5),
the tangle $(\Gamma\cap D',D')$ is 
an NS-tangle of label $k$.
This contradicts Theorem~\ref{NoNS-tangle}.
Therefore Claim~$2$ holds.

{\bf Claim~3.}
The set $\mathcal{Q}$ consists of 
exactly two elements 
each of which
intersects exactly one of $\mathcal{E}$.

{\it Proof of Claim~$3$.}
If $\mathcal{Q}$ consists of 
only one element $Q_1$, 
then $Cl((\Gamma_m\cap D)-Q_1)$ 
does not intersect $\partial D$.
Let $(\Gamma\cap D',D')$ be a tangle induced from a connected component of 
$Cl((\Gamma_m\cap D)-Q_1)$ with respect to $D$.
Then $\Gamma_m\cap\partial D'$ consists of 
one point.
Thus by (5),
the tangle $(\Gamma\cap D',D')$ is 
an NS-tangle of label $k$.
This contradicts Theorem~\ref{NoNS-tangle}.
Hence 
the condition $|\Gamma_m\cap\partial D|=2$ 
implies that
$\mathcal{Q}$ consists of 
exactly two elements. 

Now by (3), 
each element of $\mathcal{Q}$ 
intersects at least one of $\mathcal{E}$.
Suppose that there exists an element $Q_i$
intersecting at least two of $\mathcal{E}$.
Thus $Cl((\Gamma_m\cap D)-Q_i)$ consists of 
at least two connected components. 
Since $Q_i\cap\partial D\not=\emptyset$,
the condition $|\Gamma_m\cap\partial D|=2$
implies that 
$Cl((\Gamma_m\cap D)-Q_i)\cap\partial D$ 
consists of at most one point.
Hence there exists a connected component 
$X$ of 
$Cl((\Gamma_m\cap D)-Q_i)$ 
with $X\cap\partial D=\emptyset$.
Let $(\Gamma\cap D'',D'')$ be 
a tangle induced from $X$ with respect to $D$.
Then 
$\Gamma_m\cap\partial D''$ 
consists of one point.
Hence by (5),
the tangle $(\Gamma\cap D'',D'')$ is an NS-tangle of label $k$.
This contradicts Theorem~\ref{NoNS-tangle}.
Therefore Claim~$3$ holds.

By Claim~$1$, Claim~$2$, and Claim~$3$,
we have $d=p+1$. 
Thus by renumbering 
$E_1, E_2,\cdots, E_d$ and
$P_1,P_2,\cdots,P_{d-1}$,
by setting 
$P_0=Q_1,P_d=Q_2$
we can assume that
\begin{enumerate}
\item[(7)] 
for each $i=1,2,\cdots,d$ and $j=0,1,\cdots,d$,\\
$E_i\cap P_j=
\biggl\{
\begin{array}{ll}
\text{one point} &\ \ {\rm if}\ j=i-1
\text{ or }i,\\
\emptyset &~~\text{otherwise.}
\end{array}$
\end{enumerate}
Set 
$w_0=P_0\cap\partial D$, and 
$v_{d+1}=P_d\cap\partial D$. 
For each $i=1,2,\cdots,d$ and 
$j=0,1,\cdots,d$,
let $v_i=E_i\cap P_{i-1}$,  
$w_i=E_i\cap P_i$, and
\begin{enumerate}
\item[]
$L_j=$the simple arc in $P_j$ 
connecting $w_j$ and $v_{j+1}$.
\end{enumerate}
Then $E_1,E_2,\cdots,E_d$ and 
$L_0,L_1,\cdots,L_d$ 
satisfy Condition (a) and Condition (b) 
in Lemma~\ref{2ColorBeads}.

Let
$Y=(\cup_{i=1}^dE_i)\cup(\cup_{j=0}^{d}L_j)$. 
Then $Y$ is simply connected and 
$Y\cap\partial D=\{w_0,v_{d+1}\}$.

{\bf Claim~4.} 
If an edge of label $m$ intersects $D-Y$, then 
it is a terminal edge.

{\it Proof of Claim~$4$.}
If not, then
there exists the closure $T$ of a connected component
of $\Gamma_m\cap(D-Y)$ 
such that 
$T$ is a tree 
containing at least two white vertices.
Now Statement $(3)$ implies that 
$|T\cap Y|\ge 1$.
Suppose that $|T\cap Y|>1$. 
Since $Y$ is connected,
we can find a new cycle of label $m$ 
not contained in $\cup_{i=1}^dE_i$.
This contradicts that
$E_1,E_2,\cdots,E_d$ are all the disks bounded by the maximal cycles in $D$.
Hence $|T\cap Y|=1$.
Thus $T$ is a reducible tree 
with the special vertex $T\cap Y$.
This contradicts Lemma~\ref{NoReducibleTree}.
Hence Claim~$4$ holds.

Therefore Lemma~\ref{2ColorBeads} holds.
\end{proof}


{\it Proof of 
Theorem~\ref{TwoColorGateTangle}.}
Neither terminal edge nor free edge is 
an internal edge.
By Condition (ii) of an admissible tangle,
the boundary $\partial D$ intersects 
neither terminal edge nor free edge.
Thus the tangle $(\Gamma\cap D, D)$ 
satisfies Condition (i) of an IO-tangle.

By Lemma~\ref{2ColorBeads},
there exist disks $E_1,E_2,\cdots, E_d$ in {\rm Int}~$D$ and simple arcs $L_0,L_1,\cdots,L_d$ in 
$Cl(D-\cup^{d}_{i=1}E_i)$
satisfying Conditions (a), (b), (c)
in Lemma~\ref{2ColorBeads}.

The condition $|\Gamma_m\cap \partial D|=2$ 
implies that
$\partial D-(\Gamma_m\cap\partial D)$ consists of 
two connected components.
Let $L_I,L_O$ be the closures of 
the connected components.
Then $L_I,L_O$ are arcs on $\partial D$
with 
\begin{enumerate}
\item[(1)] 
$L_I\cap L_O=\partial L_I
=\partial L_O=\Gamma_m\cap \partial D$. 
\end{enumerate}
Thus the two arcs $L_I,L_O$ satisfy
Condition (iii) of an IO-tangle.

Let $Y=(\cup_{i=1}^d E_i)\cup(\cup_{j=0}^{d}L_j)$.
Then $Y\cap \partial D$ consists 
a point in $\partial L_0$
and a point in $\partial L_d$.
Now $Y$ is simply connected. 
Hence $|Y\cap\partial D|=
|\Gamma_m\cap\partial D|=2$ 
implies that
the set $D-Y$ consists of 
two connected components.
Let $\Delta_I,\Delta_O$ be 
the closures of connected components
with $\Delta_I\supset L_I$ and 
$\Delta_O\supset L_O$.
Then $\Delta_I$ and $\Delta_O$ are disks.
Let $T_I=Cl(\partial \Delta_I-L_I)$
and $T_O=Cl(\partial \Delta_O-L_O)$.
Then $T_I$ and $T_O$ are 
pseudo paths of label $m$
with side-disks $\Delta_I$ and $\Delta_O$
respectively. 

\begin{figure}
\begin{center}
\includegraphics{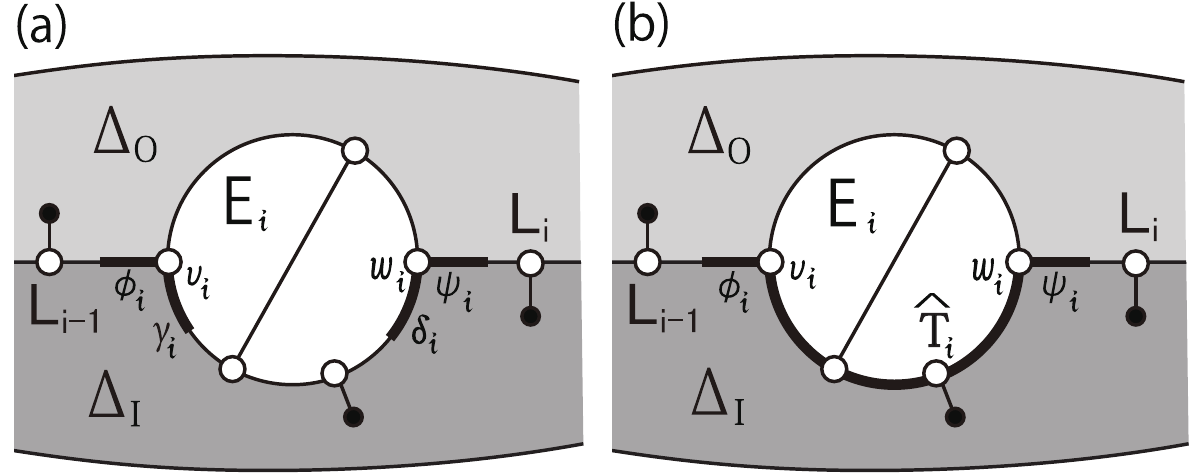}
\end{center}
\caption{ \label{fig25}
(a) The thick lines are arcs
$\gamma_i,\delta_i,\phi_i$ and $\psi_i$.
(b) The thick line is a pseudo path $\widehat T_i$.}
\end{figure}

For each $i=1,2,\cdots,d$, let 
(see Fig.~\ref{fig25})\\
$\begin{array}{l}
v_i=E_i\cap L_{i-1}, \ \ w_i=E_i\cap L_{i},\\ 
\gamma_i=
\text{an arc in the edge in $E_i\cap \Delta_I$ 
with $v_i\in \gamma_i$}, \\
\delta_i=
\text{an arc in the edge in $E_i\cap \Delta_I$ 
with $w_i\in \delta_i$},\\
\phi_i=
\text{an arc in the edge in $L_{i-1}$ 
with $v_i\in \phi_i$},\\
\psi_i=
\text{an arc in the edge in $L_{i}$ 
with $w_i\in \psi_i$},\\
T_i=E_i\cap \Delta_I, \ \ \ \widehat{T_i}=
\phi_i\cup T_i \cup \psi_i.\\
\end{array}
$\\
For each $i=0,1,\cdots,d$, let
$\widehat{L_i}=
\delta_{i}\cup L_i\cup\gamma_{i+1}$,
here
$\delta_0=\emptyset,\gamma_{d+1}=\emptyset$.

By Lemma~\ref{2ColorBeads}(c),
for each $i=1,2,\cdots,d$,
the outside edges of label $m$ for $\partial E_i$
are terminal edges except two edges 
containing one of $v_i,w_i$.
Thus by Remark~\ref{Assumption0}(2) 
we have 
\begin{enumerate}
\item[(2)]
$\mathcal{W}_O(\partial E_i,m)-
\mathcal{W}_O^{{\rm Mid}}(\partial E_i,m)
\subset \{ v_i,w_i\}$.
\end{enumerate}

Since $(\Gamma\cap D,D)$ is an admissible tangle,
the disk $D$ contains neither free edge nor simple hoop.
Since $(\Gamma\cap D,D)$ is a $2$-color tangle,
for each $i=1,2,\cdots,d$
the disk $E_i$ is a $2$-color disk 
with $\partial E_i\subset \Gamma_m$.
Thus by Lemma~\ref{ConnectedCycle}(b),

\begin{enumerate}
\item[$(3)$] 
for each $i=1,2,\cdots,d$,
the intersection
$\Gamma_m\cap E_i$ is connected.
\end{enumerate}
By $(2)$, $(3)$ and 
Lemma~\ref{qPrickle}(b), 
we have
$|\mathcal{W}_O(\partial E_i,m)-
\mathcal{W}_O^{{\rm Mid}}(\partial E_i,m)|= 2$.
Hence
\begin{enumerate}
\item[$(4)$]
$\{v_i,w_i\}=\mathcal{W}_O(\partial E_i,m)-
\mathcal{W}_O^{{\rm Mid}}(\partial E_i,m)$.
\end{enumerate}

{\bf Claim~1.}
The pseudo paths $\widehat{L_0},\widehat{L_1},\cdots,
\widehat{L_d}$ are I/O pseudo paths
with respect to $\Delta_I$.

{\it Proof of Claim~$1$.}
By (4), 
the arc $\phi_i$ is 
not middle at $v_i$,
nor
the arc $\psi_i$ is 
not middle at $w_i$. 
Hence 
by Lemma~\ref{2ColorBeads}(c),
for each $i=1,2,\cdots,d-1$
the pseudo path $\widehat{L_i}$ is 
a dichromatic bridge.
Further $\widehat{L_0}$ and 
$\widehat{L_d}$ are dichromatic piers
because $\phi_1$ and $\psi_d$
are not middle at $v_1$ and $w_d$ 
respectively.
Thus Claim~$1$ follows from 
Lemma~\ref{2ColorBridge} and
Corollary~\ref{2ColorPier}.
Hence Claim~$1$ holds.

By $(4)$,
Theorem~\ref{Prickles}(b) implies that
\begin{enumerate}
\item[$(5)$] for each $i=1,2,\cdots,d$,
the extended pseudo path $\widehat{T_i}$ 
of the path $T_i$ is an I/O pseudo path
with respect to $\Delta_I$.
\end{enumerate}

{\bf Claim~2.}
$T_I,T_O$ are I/O pseudo paths 
with respect to $\Delta_I,\Delta_O$ 
respectively.

{\it Proof of Claim~$2$.}
Let $i\in\{ 1,2,\cdots,d\}$. 
Since $\widehat{L_{i-1}}\cap\widehat{T_i}$ 
contains the white vertex $v_i$, and since 
$\widehat{T_i}\cap\widehat{L_i}$
contains the white vertex $w_i$,
the pairs 
$(\widehat{L_{i-1}},\widehat{T_i})$ 
and
$(\widehat{T_i},\widehat{L_i})$ are
I/O pairs of type I
by Claim~$1$ and $(5)$.
Thus $(\widehat{L_{0}},\widehat{T_1},
\widehat{L_{1}},\widehat{T_2},\cdots,
\widehat{L_{d-1}},\widehat{T_d},
\widehat{L_d})$
is an I/O sequence 
for $(T_I,\Delta_I)$.
Hence by Remark~\ref{NoteI/Osequence},
$T_I=\widehat{L_{0}}*\widehat{T_1}*
\widehat{L_{1}}*\widehat{T_2}*\cdots *
\widehat{L_{d-1}}*\widehat{T_d}*\widehat{L_d}$
is an I/O pseudo path 
with respect to $\Delta_I$.

Similarly $T_O$ is an I/O pseudo path 
with respect to $\Delta_O$.
Thus Claim~$2$ holds.

\begin{figure}
\begin{center}
\includegraphics{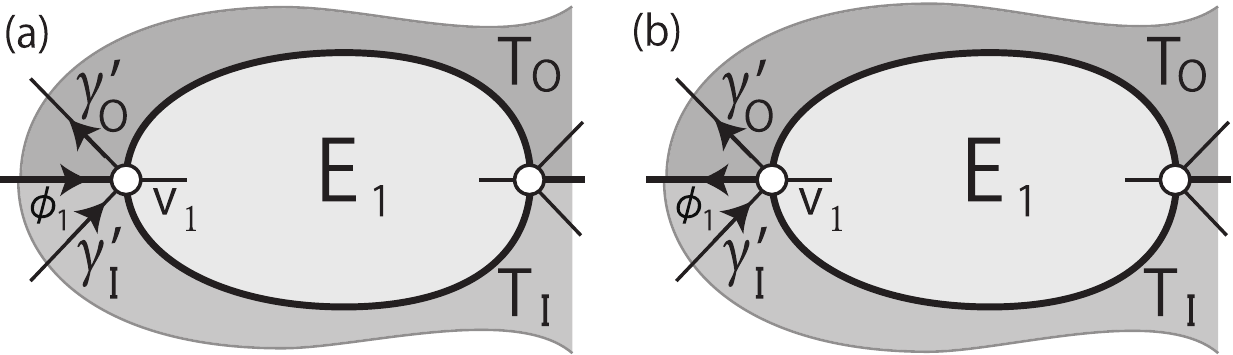}
\end{center}
\caption{ \label{fig26} 
The thick lines are of label $m$.}
\end{figure}

Let $\gamma_I'$  be 
a side-arc of $T_I$ at $v_1$ 
with respect to $\Delta_I$, 
and
$\gamma_O'$ 
a side-arc of $T_O$ at $v_1$ 
with respect to $\Delta_O$.
Without loss of generality
we can assume that 
the arc $\gamma_I'$ is oriented inward at $v_1$.
Suppose that $\phi_1$ is oriented inward at $v_1$.
Since $\phi_1$ is not middle at $v_1$ by $(4)$,
the arc $\gamma_O'$ 
is oriented outward at $v_1$ 
(see Fig.~\ref{fig26}(a)).
Suppose that $\phi_1$ is oriented outward at $v_1$.
Since $\gamma_I'$ is oriented inward at $v_1$, 
the arc $\gamma_O'$ 
is oriented outward at $v_1$ 
(see Fig.~\ref{fig26}(b)).
Namely, if $\gamma_I'$ is oriented inward at $v_1$, 
then the arc $\gamma_O'$ 
is oriented outward at $v_1$.
Since $\gamma_I'$ is a side-arc of $T_I$
with respect to $\Delta_I$
and 
since $\gamma_O'$ is a side-arc of $T_O$
with respect to $\Delta_O$,
we have that
the I/O pseudo path $T_I$ 
is inward with respect to $\Delta_I$
and the I/O pseudo path $T_O$ is outward
with respect to $\Delta_O$.

Let $p$ be a point in 
$\Gamma\cap{\rm Int}~L_I$, and
$e$ an edge of $\Gamma$ 
containing the point $p$.
Then by $(1)$, 
the label of $e$ is $k$.
Since 
the tangle $(\Gamma\cap D,D)$ is admissible, 
there exists an internal edge $\overline{e}$ of label $k$
containing $e$ such that 
each connected component of $\overline{e}\cap D$ contains a white vertex.
Let $\ell$ be the arc in $\overline{e}\cap D$ connecting the point $p$
and a white vertex $v$ of $\overline{e}$. 
Now $p\in\Delta_I$ implies $v\in\Delta_I$.
Since 
an edge of label $m$ intersecting $D-Y$ is a terminal edge by Lemma~\ref{2ColorBeads}(c),

(6) the set $D-Y$ does not contain any white vertex.\\
Hence $v\in T_I$.
Thus $\ell$ is oriented inward at $v$.
Hence the edge $e$ is locally inward at $p$.
Similarly
for any point $p\in \Gamma\cap~{\rm Int}~L_O$
there exists an edge of label $k$ locally outward at $p$.
Hence $(\Gamma\cap D,D)$ is an IO-tangle.

{\bf Claim~3.} 
$D$ does not contain any terminal edge 
of label $k$. 
 
{\it Proof of Claim~$3$.}
Suppose that $D$ contains a terminal edge $e$ 
of label $k$. 
Let $w$ be the white vertex in $e$.
Since for each $i=1,2,\cdots,d$
the $2$-color disk $E_i$ does not contain any terminal edge by Corollary~\ref{Cor3ColorDiskNoTerminalEdge}(b),
the terminal edge $e$ is not contained 
in $\cup^{d}_{i=1}E_i$. 
Thus $e$ intersects $D-\cup^{d}_{i=1}E_i$.
Hence by (6),
the white vertex $w$ is in $(\cup_{i=1}^d \partial E_i)\cup(\cup_{j=0}^{d}L_j)$.
Thus $w\in\partial E_i$ for some $i\in\{1,2,\cdots,d\}$ or $w\in {\rm Int}~L_i$ 
for some $i\in\{0,1,2,\cdots,d\}$.

Suppose that $w\in\partial E_i$ for some $i\in\{1,2,\cdots,d\}$.
There are exactly three edges of label $m$ at $w$.
Since $\partial E_i$ contains at least two white vertices
by (4), 
and since $w\in\partial E_i$,
two of the three edges are contained in $\partial E_i$. 
Thus one of $(e,w)$-edges is in $\partial E_i$.
Let $e'$ be an $(e,w)$-edge in $\partial E_i$.
Then $e'$ is an edge of label $m$ in $D$.
Since $\partial E_i$ contains at least two white vertices
by (4),
the edge $e'$ contains two white vertices.
Let $w'$ be the white vertex of $e'$ 
different from $w$. 
Since $(\Gamma\cap D,D)$ is a $2$-color tangle,
there exists an edge $e''$ of label $k$ 
at $w'$
such that $(e,e',e'')$ is 
a non-admissible consecutive triple.
This contradicts Consecutive Triplet Lemma 
(Lemma~\ref{ConsecutiveTriplet}).

Suppose that $w\in {\rm Int}~L_i$ 
for some $i\in\{0,1,2,\cdots,d\}$.
There are exactly three edges of label $m$ at $w$.
Since $L_i$ is an arc,
there exists exactly one edge $e^*$ of label $m$ 
with $L_i\cap e^*=w$.
Thus $e^*\cap (D-Y)\not=\emptyset$.
Hence by Lemma~\ref{2ColorBeads}(c),
the edge $e^*$ is a terminal edge.
Since $e$ is a terminal edge,
neither $(e,w)$-edges are terminal edges of label $m$.
Thus $e^*$ is not an $(e,w)$-edge.
Hence  both of the two $(e,w)$-edges contain
 short arcs contained in $L_i$.
Since $|L_i\cap \partial D|\le1$,
 one of the two $(e,w)$-edges is in $L_i$.
Let $e'$ be an $(e,w)$-edge in $L_i$,
and $w'$ the white vertex of $e'$ 
different from $w$. 
Since $(\Gamma\cap D,D)$ is a $2$-color tangle,
there exists an edge $e''$ of label $k$ 
at $w'$
such that $(e,e',e'')$ is 
a non-admissible consecutive triple.
This contradicts Consecutive Triplet Lemma 
(Lemma~\ref{ConsecutiveTriplet}).
Thus $D$ does not contain 
any terminal edge of label $k$. 
Hence Claim~$3$ holds.

Since $(\Gamma\cap D,D)$ is a $2$-color tangle,
by Claim~$3$
all the terminal edge in $D$ is of label $m$.
Thus $(\Gamma\cap D,D)$ is a simple IO-tangle.
Therefore
this proves Theorem~\ref{TwoColorGateTangle}.
{\hfill {$\square$}\vspace{1.5em}}

\section{Appendix}
\label{s:Appendix}

In this section, 
we give complementary lemmata 
to make our paper self-contained.\\

{\it Proof of Consecutive Triplet Lemma $($Lemma~\ref{ConsecutiveTriplet}$)$.}
Suppose that there exists 
a non-admissible consecutive triplet
$(e_1,e_2,e_3)$ in
a minimal chart.
Then $e_1$ is a terminal edge and 
the labels of $e_1$ and $e_3$ 
are same. 
Let $\partial e_2=\{w_1,w_2\}, 
w_1\in e_1$, and $w_2\in e_3$ 
(possibly $w_1=w_2$). 
Suppose that $\partial e_3=\{w_2,w_3\}$ 
(possibly $w_2=w_3$).

To make argument simple 
we assume that the three points 
$w_1,w_2$ and $w_3$ are 
mutually different.
The vertex $w_3$ may not be a white vertex.
Without loss of generality,
we can assume that 
the edge $e_1$ is oriented outward 
at the white vertex $w_1$.
In a minimal chart,
by Remark~\ref{Assumption0}(2)
any terminal edge
must contain 
a middle arc at its white vertex.
Hence 
$e_1$ contains a middle arc
at $w_1$.
Thus
the edge $e_2$ 
is oriented outward
at the white vertex $w_1$, too. 

\begin{figure}[h]
\begin{center}
\includegraphics{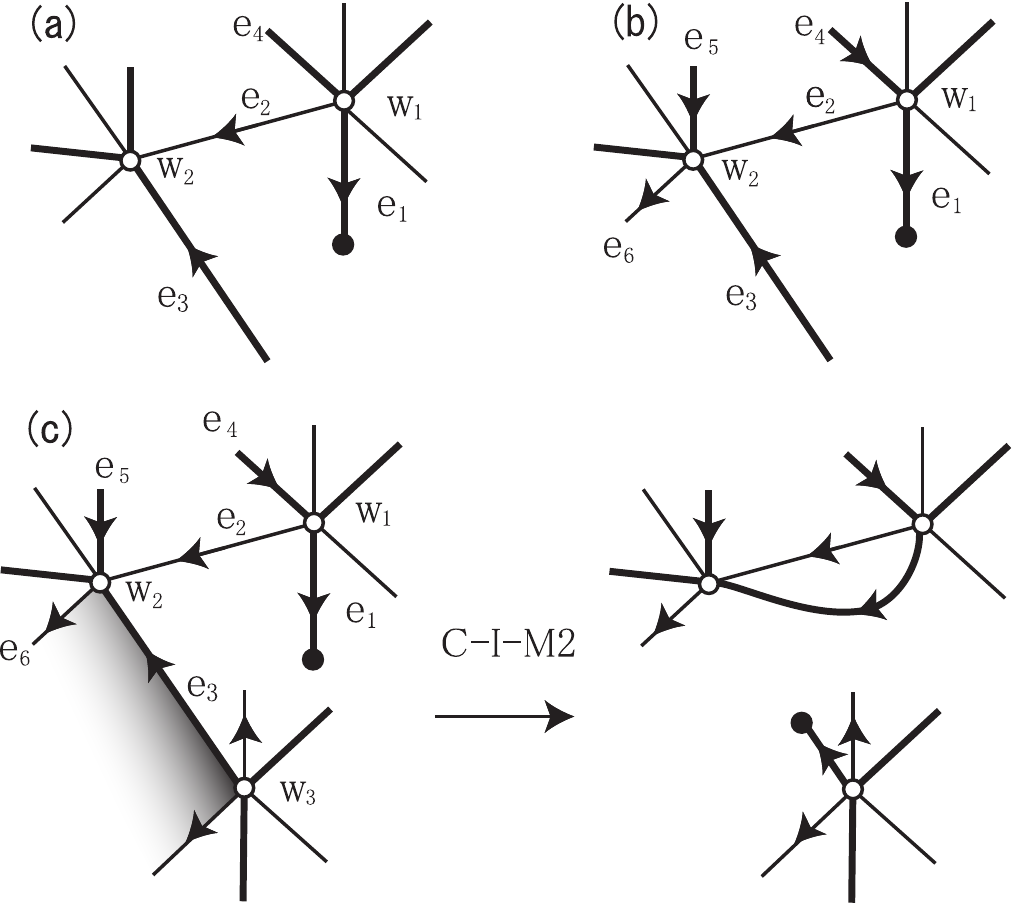}
\end{center}
\caption{ \label{fig27}}
\end{figure}

The edge $e_3$ is oriented inward 
at the white vertex $w_2$.
For, if the edge $e_3$ is oriented outward 
at $w_2$,
then the edge $e_3$ contains 
a non-middle arc at $w_2$.
By applying a C-I-M2 move
between $e_1$ and $e_3$
we can get a new terminal edge
which contains a non-middle arc
at the white vertex $w_2$.
Hence we can eliminate 
the white vertex $w_2$
by a C-III move.
This contradicts that
the chart is minimal.
Therefore
the edge $e_3$ is oriented inward 
at the white vertex $w_2$
(see Fig.~\ref{fig27}(a)).

Let $e_4$ be the $(e_2,w_1)$-edge
different from the edge $e_1$.
Since $e_1$ contains a middle arc and
oriented outward at $w_1$,
the edge $e_4$ is oriented inward
at the white vertex $w_1$.
Let $e_5$ be the $(e_2,w_2)$-edge
different from the edge $e_3$.
The edge $e_5$ is oriented inward 
at the white vertex $w_2$.
For, 
if the edge $e_5$ is oriented 
outward at the white vertex $w_2$,
then by applying a C-I-M2 move
between $e_1$ and $e_3$ and further
applying a C-I-M2 move
between $e_4$ and $e_5$,
we get three consecutive edges
connecting the two white vertices $w_1$ and $w_2$.
Hence by applying a C-I-M3 move
we can eliminate the two white vertices.
This contradicts that the chart is minimal.
Therefore the edge $e_5$ is oriented 
inward at the white vertex $w_2$
(see Fig.~\ref{fig27}(b)).

Let $e_6$ be the $(e_3,w_2)$-edge
different from $e_2$.
Since the three edges $e_2, e_3, e_5$ are oriented inward
at the white vertex $w_2$,
the edge $e_6$ is oriented outward
at the white vertex $w_2$
(see Fig.~\ref{fig27}(b)).

The vertex $w_3$ is a white vertex.
For, if $w_3$ is not a white vertex,
then $w_3$ is a crossing or a black vertex.
Hence the edge $e_3$ is not contained in a bigon.
By applying a C-I-M2 move between $e_1$ and $e_3$,
we can get a new bigon without destroying old bigons.
Thus the number of bigons increases.
This contradicts that the chart is minimal.
Hence the vertex $w_3$ must be a white vertex.

Since we can apply a C-I-M2 move between $e_1$ and $e_3$,
the edge $e_3$ must contain a middle
arc at the white vertex $w_3$.
Hence the $(e_3,w_3)$-edges are oriented 
outward at $w_3$.
Since the edge $e_6$ is oriented 
outward at $w_2$ (see Fig.~\ref{fig27}(c)),
neither of the $(e_3,w_3)$-edges is equal to $e_6$.
Hence the edge $e_3$ is not contained in
a bigon.

Now 
by applying a C-I-M2 move
between $e_1$ and $e_3$
we can get a new bigon
without destroying 
old bigons.
Thus the number of bigons
increases.
This contradicts 
that the chart is minimal. 
{\hfill {$\square$}\vspace{1.5em}}

As a consequence of Theorem~\ref{NoNS-tangle},
we can show the following lemma.

\begin{lemma}{\rm (Boundary Condition Lemma) 
(\cite[Lemma 4.1]{TwoCrossingI}) }
\label{LemmaBoundaryCondition}
Let $(\Gamma\cap D,D)$ be a tangle 
in a minimal chart $\Gamma$ 
such that 
$D$ does not contain any crossing,
free edge nor simple hoop.
Let $a=\overline\alpha(\Gamma\cap\partial D)$ and
$b=\overline\beta(\Gamma\cap\partial D)$.
Then
 $\Gamma_i\cap D=\emptyset$ 
 except for $a\le i \le b$.
\end{lemma}

\begin{proof}
We claim that $D$ does not contain any hoop.
If $D$ contains a hoop, say $C'$,
then $C'$ is not a simple hoop 
by the assumption.
Hence $C'$ bounds a disk $E'$ containing a white vertex.
Let $D'$ be a regular neighborhood of $E'$.
Then $\Gamma\cap \partial D'=\emptyset$.
Since $D$ does not contain any crossing,
neither does $D'$.
Thus $(\Gamma\cap D',D')$ is an NS-tangle.
This contradicts Theorem~\ref{NoNS-tangle}.
Hence 

(1) $D$ does not contain any hoop.

Let $i$ be a label with $\Gamma_i\cap D\not=\emptyset$.
Suppose that $i<a$ or $b<i$.
Without loss of generality
we can assume $i<a$.
Let $\alpha=\overline\alpha(\Gamma\cap D)$.
Since
$\alpha\le i<a=\overline\alpha(\Gamma\cap\partial D)$,
we have $\Gamma_\alpha\cap\partial D=\emptyset$.
Let $G$ be a connected component of $\Gamma_\alpha\cap D$.
Then $G\cap\partial D=\emptyset$.
Let $(\Gamma\cap D^*,D^*)$ be
a tangle induced from $G$ with respect to $D$.
Then $D^*\subset D$.
Since $G\cap\partial D=\emptyset$,
by Remark~\ref{RemTangleInducedFromG}$(2)$
we have $|\Gamma_\alpha\cap\partial D^*|=
|G\cap\partial D|=0$.
Hence 

(2) $\Gamma_\alpha\cap \partial D^{*}=\emptyset$.\\
Since $\alpha$ is the lowest label in $D$,

(3) for each label $j<\alpha$, we have 
$\Gamma_j\cap \partial D^{*}=\emptyset$.\\
Since $D$ does not contain any crossing by assumption, 
for each label $j>\alpha+1$
we have $G\cap\Gamma_j=\emptyset$. 
Hence 
by Remark~\ref{RemTangleInducedFromG}$(1)$

(4) for each label $j>\alpha+1$, 
we have $\Gamma_j\cap \partial D^{*}=\emptyset$ .\\
Since $D^*\subset D$, by assumption

(5) $D^{*}$ does not contain any crossing.\\
By (1), 
the disk $D^*$ does not contain any hoop.
By the assumption, 
the disk $D^*$ does not contain any free edge.
By (5), the disk $D^*$ does not contain any ring.
Hence $D^*$ contains at least one white vertex.
Thus 
by (2), (3), (4) and (5),
the tangle $(\Gamma\cap D^{*},D^{*})$ is 
an NS-tangle of label $\alpha+1$.
This contradicts Theorem~\ref{NoNS-tangle}.
Therefore $a\le i\le b$.
\end{proof}


In Section~\ref{s:Intro},
we have mentioned the following:
Let $m$ be any label of 
a chart $\Gamma$,
and $E$ a disk 
with $\partial E\subset\Gamma_m$ but
without crossings, free edges nor simple hoops.
If $\Gamma$ is a minimal chart,
then we can show that  
$E$ is a $3$-color disk.
Further, if $m$ is 
the minimal label or the maximal label 
of the chart, 
then 
$E$ is a $2$-color disk.

This can be shown by the following lemma.

\begin{lemma}
\label{Appendix3ColorDisk}
Let $\Gamma$ be a minimal chart.
Let $C$ be a cycle of label $m$ bounding a disk $E$ 
without crossings, free edges nor simple hoops. 
Then we have the following:
\begin{enumerate}
\item[{\rm $($a$)$}] $E$ is a $3$-color disk.
\item[{\rm $($b$)$}] If there exists a label $k$ with $|m-k|=1$ such that 
all the white vertices in $C$ are in $\Gamma_m\cap\Gamma_k$,
then $E$ is a $2$-color disk. 
\end{enumerate}
\end{lemma}

\begin{proof}
Let $(\Gamma\cap D,D)$ be a tangle induced from the cycle $C$.
Then $E\subset D$.
Since $\partial E$ is a cycle of label $m$,
we have $m-1\le \overline\alpha(\Gamma\cap\partial D)$ 
and $\overline\beta(\Gamma\cap\partial D)\le m+1$.
By Boundary Condition Lemma (Lemma~\ref{LemmaBoundaryCondition}),
we have 
$m-1\le \overline\alpha(\Gamma\cap D)$ 
and $\overline\beta(\Gamma\cap D)\le m+1$.
Thus $\Gamma\cap D\subset \Gamma_{m-1}\cup\Gamma_m\cup\Gamma_{m+1}$.
Since $E\subset D$,
the disk $E$ is a $3$-color disk.

Similarly we can show Statement (b).
\end{proof}


{\bf Index of Notations}
\vspace{2mm}

$
\begin{array}{ll||}
w(\Gamma) & p 2\\
f(\Gamma) & p 2\\
c(\Gamma) & p 3\\
b(\Gamma) & p 3\\
\Gamma_m & p 4  \\
{\mathcal W}_I^{{\rm Mid}}(C,m) & p 5, p 17 \\
{\mathcal W}_O^{{\rm Mid}}(C,m) & p 5, p 17\\
\end{array}
$
~~
$
\begin{array}{ll||}
{\rm Main}(\Gamma) & p 11 \\
{\mathcal W}_I(C,m)  & p 17 \\
{\mathcal W}_O(C,m) & p 17 \\
{\mathcal W}^{{\rm Mid}}(E,\pm1) & p 18 \\
L(P^*) & p 21 \\
{\mathcal P}(C;{\mathcal S}) & p 26 \\
{\rm deg}_Gv & p 27 \\
\end{array}
$
~~
$
\begin{array}{ll}
\tau(D) & p 32 \\
\overline\alpha(X) & p 33 \\
\overline\beta(X) & p 33 \\
P_1^* *P_2^* & p 36, p 38 \\
P_1^* *P_2^* *\cdots *P_s^* & p 39 \\
& \\
& \\
\end{array}
$

\newpage

\vspace{6mm}

{\bf Index of words}

\vspace{2mm}

{\small $
\begin{array}{ll||}
\text{admissible pseudo path} & p 21 \\
\text{admissible tangle} & p 5 \\
\text{associated edge sequence} & p 21\\
\text{associated side-arc sequence} & p 22 \\
\text{associated vertex sequence} & p 21 \\
\text{bigon} & p 3 \\
\text{bridge}  & p 42 \\
\text{chart} & p 3 \\
\text{C-move equivalent} & p 1 \\
\text{consecutive triplet} & p 13 \\
\text{cycle of label $m$} & p 4, p 17 \\
\text{dichromatic bridge} & p 42 \\
\text{dichromatic pier} & p 44 \\
\text{dichromatic pseudo path} & p 22\\
\text{end-arc} & p 21 \\
\text{$(e,w)$-edge} & p 14 \\
\text{extended pseudo path} & p 25 \\
\text{finite domain} & p 10 \\
\text{free edge} & p 2\\
\text{fundamental information} & p 29 \\
\text{hoop} & p 4 \\
\text{inside edge} & p 4, p 17 \\
\text{internal arc} & p 12 \\
\text{internal edge} & p 5 \\
\text{inward pseudo path} & p 21 \\
\text{locally inward} & p 6 \\
\text{locally outward} & p 6 \\
\text{loop} & p 13 \\
\text{I/O pair of type I} & p 36 \\
\text{I/O pair of type II} & p 37 \\
\text{I/O pseudo path} & p 21 \\
\end{array}
$}
~~
{\small $
\begin{array}{ll}
\text{I/O sequence for $(P^*,D)$} & p 39 \\
\text{IO-tangle} & p 6 \\
\text{maximal cycle in $X$} & p 29 \\
\text{middle arc} & p 2, p 3\\
\text{middle at $v$} & p 4, p 21 \\
\text{minimal chart} & p 4 \\
\text{minimal NS-tangle} & p 32 \\
\text{NS-tangle} & p 5 \\
\text{nice side-disk for a pier} & p 44\\
\text{one-side pseudo path} & p 22 \\
\text{outside edge} & p 4, p 17 \\
\text{outward pseudo path} & p 21 \\
\text{path} & p 20 \\
\text{path decomposition ${\mathcal P}(C;{\mathcal S})$} & p 26 \\
\text{pier} & p 44 \\
\text{pseudo path} & p 20 \\
\text{reducible tree} & p 27 \\
\text{ring} & p 12 \\
\text{secondary label} & p 22, p 26 \\
\text{side-arc} & p 21 \\
\text{side-disk} & p 20 \\
\text{simple hoop} & p 4 \\
\text{small component} & p 29 \\
\text{suspicious cycle} & p 28 \\
\text{tangle} & p 5 \\
\text{tangle induced from a cycle} & p 33 \\
\text{tangle induced from $G$} & p 35 \\
\text{terminal edge} & p 6 \\
\text{2-color disk} & p 4 \\
\text{2-color tangle} & p 35 \\
\text{3-color disk} & p 4 \\
\end{array}
$}

\vspace{10mm}

\begin{minipage}{65mm}
{Teruo NAGASE
\\
{\small Tokai University \\
4-1-1 Kitakaname, Hiratuka \\
Kanagawa, 259-1292 Japan\\
\\
nagase@keyaki.cc.u-tokai.ac.jp
}}
\end{minipage}
\begin{minipage}{65mm}
{Akiko SHIMA 
\\
{\small Department of Mathematics, 
\\
Tokai University
\\
4-1-1 Kitakaname, Hiratuka \\
Kanagawa, 259-1292 Japan\\
shima@keyaki.cc.u-tokai.ac.jp
}}
\end{minipage}

\end{document}